\documentclass[a4paper]{article}
\usepackage{graphics}
\usepackage[english]{babel}
\usepackage{geometry}
\usepackage[numbers]{natbib}
\usepackage{mathrsfs}
\usepackage{amsmath}
\usepackage{amssymb}
\usepackage{amsthm}
\usepackage{enumerate}
\usepackage{bm}
\usepackage{color}
\usepackage{stmaryrd}
\usepackage{graphicx}
\usepackage{manfnt}
\usepackage{hyperref}
\usepackage{authblk}
\usepackage{comment}
\usepackage{graphbox}

\newtheorem{theorem}    {Theorem}[section]

\newcommand{\ie}{\textit{i.e.,}}

\newcommand{\cf}{\textit{cf.}}
\newcommand{\etc}{\textit{etc.}}

\newcommand{\reals}{\mathbb{R}}
\newcommand{\metric}{\bm{\mathfrak{G}}}

\newcommand{\pd}[2]{\frac{\partial #1}{\partial #2}}
\newcommand{\ddd}[1]{\tilde{\Delta}#1}
\let\mbs=\bm

\newcommand\revone{\textcolor{black}}
\newcommand\revtwo{\textcolor{black}}

\begin{document}

\title{A new conservative/dissipative time integration scheme\\ for nonlinear mechanical systems}

\author[1]{Cristian G. Gebhardt}
\author[2,3]{Ignacio Romero}
\author[1]{Raimund Rolfes}

\affil[1]{Institute of Structural Analysis and ForWind Hannover, Leibniz Universit\"at Hannover, Appelstra\ss e 9 A, 30167 Hannover, Germany}

\affil[2]{IMDEA Materials Institute, C/ Eric Kandel 2, Tecnogetafe, Madrid 28906, Spain}

\affil[3]{Universidad Polit\'ecnica de Madrid, Jos\'e Guti\'errez Abascal 2, Madrid 29006, Spain}

\date{}

\maketitle    

\begin{abstract}
	We present a conservative/dissipative time integration scheme for nonlinear mechanical
	systems. Starting from a weak form, we derive algorithmic forces and velocities that guarantee the
	desired conservation/dissipation properties. Our approach relies on a collection of linearly
	constrained quadratic programs defining high order correction terms that modify, in the minimum
	possible way, the classical midpoint rule so as to guarantee the strict energy
	conservation/dissipation properties. The solution of these programs provides explicit formulas for
	the algorithmic forces and velocities which can be easily incorporated into existing
	implementations. Similarities and differences between our approach and well-established
	methods are discussed as well. The approach, suitable for reduced-order models, finite element
	models, or multibody systems, is tested and its capabilities are illustrated by means of several
	examples.
\end{abstract}

\section{Introduction}

A key feature in the numerical approximations of conservative mechanical systems is their ability to
exactly preserve the first integrals of their motion (energy, momenta, symplecticity, \ldots), replicating the
properties of the continuous counterparts (see, e.g., \citep{Arnold1989,Hairer:2002vg}). This
interest in \emph{structure preserving integrators} is hence justified by the qualitative similarity
between the dynamical behaviour of a mechanical system and the \emph{discrete dynamics} generated
by the time integration scheme \cite{Stuart:1996ve}. In addition, a wealth of evidence supports
the fact that this kind of time-stepping methods behaves extremely well for long-term simulations
\cite{Simo1992,Simo1994,Gonzalez1996,Kane1999,McLachlan1999,Armero2001_1}.

It is not easy to formulate numerical schemes that unconditionally preserve one or more invariants
of the discrete motion. Generally speaking, this goal is accomplished by ensuring
that some of the (abstract) geometric structures that appear in the continuous picture are
replicated in the discrete dynamics. Since it is well-known that, in general, all invariants cannot be
preserved for a fixed time step size scheme, different families of methods strive for the preservation
of specific subsets of the various symmetries of the continuous system. For example, some numerical methods
resemble discrete Hamiltonian systems \cite{Gonzalez1996}, based on discrete gradient operators,
and unconditionally preserve the energy and the (at most quadratic) momenta. Other methods emanate
from discrete variational principles \cite{Marsden:2001vi} and obtain the update formula from
the stationarity conditions of these principles. In fact, it is possible to formulate methods
that preserve energy, momenta, and the symplectic form of the system, if the time step size
is added as an unknown to the method's equations \cite{Kane1999}.

In the context of nonlinear elastodynamics, the first energy and momentum conserving algorithms were
developed by Simo and co-workers \cite{Simo1992}. This pioneering work showed that for Saint
Venant-Kirchhoff materials, such structure preserving methods can be easily obtained by a simple
modification of the midpoint rule in which the \revone{stress, instead of being evaluated at the midpoint
instant, should be taken as the average of the stresses at the endpoints of the time interval. Since the constitutive law is linear in the strain, this turns to be equivalent to compute the algorithmic stress with the average of the strains at the endpoints of the time interval.} This
simple idea was later applied to the conserving integration of shells \cite{Simo1994}, rods
\cite{Simo:1995tz,Romero:rod-2002}, contact mechanics \cite{Armero:1998th}, multibody systems
\cite{Goicolea:2000uf, Betsch2010}, \etc, and generalized to elastic materials of arbitrary type
\cite{Gonzalez:2000wx,Laursen:2001ue}. The key idea for such generalization is the definition of a
\emph{discrete gradient} operator, a consistent approximation of the gradient that guarantees the
strict conservation of energy in Hamiltonian systems \cite{Gotusso1985a, Itoh1988a,
Gonzalez:2000wx,Romero2012}. Alternatively, one might derive conserving methods by defining an
\emph{average vector field} \cite{Harten1983, McLachlan1999}. \revtwo{In the context of the continuous Galerkin method, an optimization
	approach was employed to systematically develop high-order energy conserving schemes \cite{French1990, Gross2005}. Very recently,
a new mixed variational framework that takes advantage of the structure	of polyconvex stored energy functions was proposed \cite{Betsch2018}, and the properties of several formulas for the \textit{discrete gradient} that are available in the literature were carefully analyzed in the context of multibody systems \cite{GarciaOrden2019}.}

Many Hamiltonian problems are modeled with stiff differential equations for which conserving
integration schemes might not be the most robust. For these problems, numerical methods with
controllable numerical dissipation in the high-frequency range provide often a practical solution
\cite{Kuhl:1996vu, Kuhl:1999tq, Bottasso:1997wm, Bottasso:2001te, Romero:2016td}. Based on
a modification of the discrete gradient operator, Armero and Romero
\cite{Armero2001_1,Armero2001_2} developed a family of schemes for nonlinear three-dimensional
elastodynamics that exhibits this kind of algorithmic dissipation, while preserving the momenta and
providing a strict control of the energy, applicable to elastodynamics, as well as to rods and shells
\cite{Romero:2002wb,Armero:2003vt}. Following an alternative path based on the
\emph{average vector field}, Gebhardt and co-workers have proposed similar conserving/dissipative methods
for general solid and structural problems \cite{Gebhardt2019a, Gebhardt2019b}.

This work considers the conservative/dissipative time integration of the equations of motion that
typically arise during the analysis of nonlinear mechanical systems. More specifically, we present a
novel approach that renders, by construction, methods with the desired conservation or dissipation
properties. These methods discretize the equations of motion and add some perturbations related to
the main field variables through a collection of ancillary linearly constrained quadratic programs
that guarantee the conservation/dissipation properties. This kind of programs are analytically
solvable and therefore, very attractive from the computational point of view. One possible
interpretation of the contributions in this article is that it results in conservative/dissipative
methods where the geometric arguments typically employed for their design have been replaced by
optimality conditions.

The perturbations proposed in the new methods are designed to correct some of the unwanted effects
coming from the discretization of the governing equations. From a geometric point of view, the idea
is to redesign the problem in such a way that the behavior of the system on the discrete constrained
sub-manifold remains unaltered, but acts as an attractor for trajectories outside of it. Since the
constrained programs can be solved in closed form, corrected formulas for the algorithmic internal
forces and generalized velocities can be provided, and thus easily incorporated in existing
simulation codes. The similarities and differences of the newly proposed method with respect to
existing ones are pointed out and discussed critically.

The remaining of the article is organized as follows: In Section~\ref{sec-mechanics}, we present the
basic framework for nonlinear mechanical systems. In Section~\ref{sec-time}, we address in a
comprehensive manner the new time discretization.  In Section~\ref{sec-results}, we present several
examples of increasing complexity for the verification of the method.  Finally, conclusions,
limitations and future work are presented in Section~\ref{sec-closure}.  Additionally, the Appendix
introduces the precision quotient, with which the correctness of an implementation can be tested.

\section{Mechanical framework}
\label{sec-mechanics}

\subsection{Statement}

In this work we consider mechanical systems whose configuration is completely
defined by a vector $\bm{q}\in Q$, where $Q\subseteq \reals^n$.
Denoting by $t$ the time, the state of the system at any instant is given
by the pair $(\bm{q}, \bm{s})\in W\equiv TQ$, where $\bm{s}=\dot{\bm{q}}$ is the velocity,
and in which we have employed the notation $\dot{(\cdot)} = \frac{\mathrm{d}(\cdot)}{\mathrm{d}t}$.
The dynamical behavior of this system, for $t\in[t_a,t_b]$ is governed by
the variational equation:
\begin{equation}
  \label{eq-ode}
  \int_{t_{a}}^{t_{b}}
  \left[ \bigl\langle
    \delta\bm{s},\bm{p}(\dot{\bm{q}})-\bm{\pi}(\bm{s})\bigr\rangle-
    \bigl\langle
    \delta\bm{q},\dot{\bm{\pi}}(\bm{s})+\bm{f}^{\textrm{int}}
    (\bm{q})-\bm{f}^{\textrm{ext}}(\bm{q})\bigr\rangle \right] \textrm{d}t \\
  = 0\:,
\end{equation}
where $(\delta\bm{q},\delta\bm{s})\in TW$ are admissible
variations of the generalized coordinates and velocities,
$\bm{p}(\dot{\bm{q}})\in T_{\bm{s}}^*S$
and $\bm{\pi}(\bm{s})\in T_{\bm{s}}^*S$
stand for the generalized-coordinate-based and generalized-velocity-based
momenta, respectively,
$\bm{f}^{\textrm{int}}\in T_{\bm{q}}^*Q$
is the vector of internal forces, $\bm{f}^{\textrm{ext}}\in T_{\bm{q}}^*Q$
is the vector of external loads that can be of conservative or non-conservative
nature, and finally, $\bigl\langle\cdot,\cdot\bigr\rangle$ represents
a suitable pairing. Additionally, we assume the following two conditions:
\emph{i}) the system possesses a positive-definite symmetric mass matrix~$\mbs{M}$
such that
\begin{equation}
  \mbs{\pi}(\mbs{s}) = \mbs{M} \mbs{s}\ ,
  \qquad
  \mbs{p}(\dot{\mbs{q}}) = \mbs{M} \dot{\mbs{q}}\ 
  \label{eq-mass-matrix}
\end{equation}
and, \emph{ii}) both the internal and the external forces derive from potential
functions depending only on the configuration $\bm{q}$, \ie
\begin{equation}
  \mbs{f}^{\mathrm{int}} = - \pd{V^{\mathrm{int}}}{\mbs{q}}
  \ ,\qquad
  \mbs{f}^{\mathrm{ext}} = - \pd{V^{\mathrm{ext}}}{\mbs{q}}\ ,
  \label{eq-potentials}
\end{equation}
and we define the total potential energy of the system as $V=V^{\textrm{int}}+V^{\textrm{ext}}$.

We would like to analyze next the implications that symmetry has on the form of the internal forces
and the appropriate notions of linear and angular momentum in the abstract space~$Q$. For that, the
relation between the configuration space $Q$ and the ambient space~$\reals^3$ has to be carefully
considered.  We start by defining $\Phi:\reals^3 \times Q\to Q$ to be a smooth action of $\reals^3$
on the configuration space such that $\Phi(\mbs{a},\mbs{q})$ is the configuration of the system
after all its points have been translated in space by constant vector $\mbs{a}$. The infinitesimal
generator of this translation at~$\mbs{q}$ is the vector $\mbs{\tau}_{\mbs{a}}(\mbs{q})\in
T_{\mbs{q}}Q$ defined as
\begin{equation}
  \mbs{\tau}_{\mbs{a}}(\mbs{q})
  =
  \left.
    \pd{}{\epsilon}
  \right|_{\epsilon=0}
  \Phi(\epsilon\, \mbs{a}, \mbs{q})\ ,
  \label{eq-inf-generator}
\end{equation}
with $\epsilon\in\reals$. Let us now assume internal potential energy is invariant under
translations, \ie
\begin{equation}
  V^{\mathrm{int}} = V^{\mathrm{int}} \circ \Phi\ .
  \label{eq-trans-invariant}
\end{equation}
Then, choosing a one parameter curve of translations $\Phi(\epsilon\,\mbs{a},\cdot)$ in
Eq.~\eqref{eq-trans-invariant} and differentiating with respect to $\epsilon$, it
follows that a translation invariant potential implies that the internal
forces satisfy
\begin{equation}
  0
  =
  \left.
    \pd{}{\epsilon}
  \right|_{\epsilon=0}
  V^{\mathrm{int}} (\Phi(\epsilon\,\mbs{a},\mbs{q}))
  =
  \langle
  \pd{}{\mbs{q}} V^{\mathrm{int}} (\mbs{q}), \mbs{\tau}_{\mbs{a}}(\mbs{q})
  \rangle
  =
  - \langle \mbs{f}^{\textrm{int}}(\mbs{q}), \mbs{\tau}_{\mbs{a}}(\mbs{q}) \rangle .
  \label{eq-trans-invariance2}
\end{equation}
To study the conservation of angular momentum, we must
repeat the same argument but considering now a second
smooth action $\Psi:\reals^3\times Q \to Q$ such that
$\Psi(\mbs{\theta},\mbs{q})$ is the configuration of the system
after all its points have rotated in ambient space by
the application of a rotation $\exp[\hat{\mbs{\theta}}]$. Defining,
as before, the infinitesimal generator of this action to be
the vector $\mbs{\rho}_{\mbs{\theta}}(\mbs{q})\in T_{\mbs{q}}Q$
calculated as
\begin{equation}
  \mbs{\rho}_{\mbs{\theta}}(\mbs{q})
  =
  \left.
    \pd{}{\epsilon}
  \right|_{\epsilon=0}
  \Psi(\epsilon\, \mbs{\theta}, \mbs{q})\ ,
  \label{eq-inf-generator-rot}
\end{equation}
again with $\epsilon\in\reals$. If the potential energy is now
rotation invariant, \ie
\begin{equation}
  V^{\textrm{int}} = V^{\textrm{int}} \circ \Psi\, .
\end{equation}
Then the internal force must satisfy
\begin{equation}
  0
  =
  \left.
    \pd{}{\epsilon}
  \right|_{\epsilon=0}
  V^{\textrm{int}} (\Psi(\epsilon\,\mbs{\theta},\mbs{q}))
  =
  \langle
  \pd{}{\mbs{q}} V^{\textrm{int}}(\mbs{q}), \mbs{\rho}_{\mbs{\theta}}(\mbs{q})
  \rangle
  =
  - \langle \mbs{f}^{\textrm{int}}(\mbs{q}), \mbs{\rho}_{\mbs{\theta}}(\mbs{q}) \rangle .
  \label{eq-rot-invariance2}
\end{equation}

The precise notion of linear and angular momentum for the
system defined in this section is provided by the following
result:

\begin{theorem} Consider a mechanical system with configuration space $Q\subseteq\reals^n$ and
vanishing external forces. Let $\Phi(\mbs{a},\cdot),\Psi(\mbs{\theta},\cdot)$ be the translation and
rotation actions on the configuration space with infinitesimal generators $\mbs{\tau}_{\mbs{a}}$ and
$\mbs{\rho}_{\mbs{\theta}}$, respectively, and define the linear momentum $\mbs{l}\in\reals^3$ and the
angular momentum $\mbs{j}\in\reals^3$ as the two quantities that verify
\begin{equation}
  \langle \mbs{l}, \mbs{a} \rangle = \langle \mbs{\tau}_{\mbs{a}}(\mbs{q}), \mbs{\pi} \rangle
  \ , \qquad
  \langle \mbs{j}, \mbs{\theta} \rangle = \langle \mbs{\rho}_{\mbs{\theta}}(\mbs{q}), \mbs{\pi} \rangle
  \:.
  \label{eq-momenta}
\end{equation}
Then, the linear momentum is conserved if the potential energy is invariant with respect to
translations.
Similarly, if the potential energy is invariant under rotations, the angular momentum is a
constant of the motion. Moreover, the total energy
\begin{equation}
  E = \frac12 \langle \mbs{s}, \mbs{M} \mbs{s} \rangle + V(\mbs{q})
  \label{eq-energy}
\end{equation}
is preserved by the motion, due to its time invariance.
\end{theorem}

\begin{proof}
  The proof of momenta conservation follows from taking the derivative of these quantities and using
\eqref{eq-ode} with admissible variations $(\delta \mbs{q}, \delta \mbs{s}) =
(\mbs{\tau}_{\mbs{a}}(\mbs{q}), \mbs{0})$ and $(\mbs{\rho}_{\mbs{\theta}}(\mbs{q}), \mbs{0})$, respectively. The
conservation of energy property follows similarly by choosing $(\delta \mbs{q}, \delta \mbs{s}) = (
\mbs{s}, \mbs{0})$.
\end{proof}

\section{Time discretization}
\label{sec-time}
The interest in the current work is in algorithms to approximate the solution of Eq.~\eqref{eq-ode}.
To define them, let us start by considering a partition of the interval $[t_a,t_b]$
into disjoint subintervals $(t_n,t_{n+1}]$ with $t_a=t_0<t_1<\ldots<t_N=t_b$, and
$\Delta t_n = t_{n+1}-t_n$. Then, the integration algorithms that we consider are based on
the midpoint approximation of Eq.~(\ref{eq-ode}) and of the form:
\begin{equation}
  \label{eq-general-integration}
    0
    =
    \langle
    \delta\bm{s},
    \mbs{M} \frac{\mbs{q}_{n+1} - \mbs{q}_n}{\Delta t_n} - \mbs{M} \mbs{s}_{n+1/2}
    \rangle
    \\
    -
    \langle
    \delta\bm{q},
    \frac{\mbs{\pi}_{n+1}-\mbs{\pi}_n}{\Delta t_n}
    +\mbs{\mathsf{f}}^{\mathrm{int}}
    (\bm{q}_{n},\bm{q}_{n+1})-\bm{f}^{\mathrm{ext}}(\bm{q}_{n+1/2})
    \rangle
    \:
\end{equation}
where the configuration and rate, respectively, at time $t_n$ are approximated by
$\mbs{q}_n,\mbs{s}_n$, we have defined $\mbs{\pi}_n = \mbs{M}\mbs{s}_n$, and we have used the
notation $(\cdot)_{n+1/2} = \frac12 (\cdot)_n + \frac12 (\cdot)_{n+1}$. The
update depends on the definition of an approximation to the internal force
at the midpoint $t_{n+1/2}$ that we have denoted as $\mbs{\mathsf{f}}^{\mathrm{int}}$.

Eq.~\eqref{eq-general-integration} provides an implicit or explicit update
$(\bm{q}_n,\bm{s}_n)\mapsto (\bm{q}_{n+1},\bm{s}_{n+1})$ that, together with the initial
conditions of the configuration and velocity, suffices to generate discrete trajectories. We note
that the approximation to the internal force in Eq.~\eqref{eq-general-integration} is a function of
two arguments that, by consistency, must satisfy
\begin{equation}
  \mbs{\mathsf{f}}^{\mathrm{int}}(\mbs{q},\mbs{q}) = \mbs{f}^{\mathrm{int}}(\mbs{q})\ ,
  \label{eq-f-discrete}
\end{equation}
for all configurations $\mbs{q}\in Q$.

We are interested, in particular, in formulating time integration schemes of
the form~\eqref{eq-general-integration} that preserve (some of) the invariants in
the motion of the system~\eqref{eq-ode}, while controlling the value of the energy at all times.
Let us first consider the update Eq.~\eqref{eq-general-integration} with
variations of the form $(\delta \bm{q}, \delta \bm{s})=(\bm{0},\bm{c})$, where
$\mbs{c}$ is an arbitrary but constant vector in $TQ$. Then, trivially, it
follows that
\begin{equation}
  \frac{\mbs{q}_{n+1} - \mbs{q}_n}{\Delta t_n}
  =
  \mbs{s}_{n+1/2}\ .
  \label{eq-velocity-update}
\end{equation}
Next, we would like to explore whether the proposed class of integration schemes preserves momenta
for mechanical systems defined in the configuration space~$Q$.
The result, for every configuration space, is given next.

\begin{theorem}
  Consider the time discretization \eqref{eq-general-integration} of a mechanical
  system with configuration space $Q\subseteq\reals^n$. Let $\Phi$ and $\Psi$ denote,
  as above, the actions of $\reals^3$ on $Q$ representing, respectively, translations
  and rotations. The integration scheme preserves linear momentum if
  
  %
  \begin{subequations}
    \label{eq-lm-conservation}
    \begin{gather}
      0 = \langle \mbs{\tau}_{\mbs{a}}(\mbs{q}_{n+1/2}) , \mbs{\mathsf{f}}^{\mathrm{int}}(\mbs{q}_n,\mbs{q}_{n+1}) \rangle
          \label{eq-lm-cons1} \\
      \langle \mbs{\tau}_{\mbs{a}}(\mbs{q}_{n+1/2}) , \mbs{\pi}_{n+1} - \mbs{\pi}_n \rangle
        =
          \langle \mbs{l}_{n+1} - \mbs{l}_n , \mbs{a}\rangle\ .
          \label{eq-lm-cons2}
    \end{gather}
  \end{subequations}

for every $\mbs{a}\in\reals^3$. Likewise, the integration algorithm preserves angular momentum
if for every $\mbs{\theta}\in\reals^3$

 \begin{subequations}
    \label{eq-am-conservation}
    \begin{gather}
      0 = \langle \mbs{\rho}_{\mbs{\theta}}(\mbs{q}_{n+1/2}) ,
      \mbs{\mathsf{f}}^{\mathrm{int}}(\mbs{q}_n,\mbs{q}_{n+1}) \rangle
      \label{eq-am-cons1} \\
      \langle \mbs{\rho}_{\mbs{\theta}}(\mbs{q}_{n+1/2}) ,
      \mbs{\pi}_{n+1} - \mbs{\pi}_n \rangle
      =
      \langle \mbs{j}_{n+1} - \mbs{j}_n , \mbs{\theta}\rangle\ .
      \label{eq-am-cons2}
    \end{gather}
  \end{subequations}
\end{theorem}

The verification of conditions~\eqref{eq-lm-conservation}-\eqref{eq-am-conservation} depends, first, on the
structure of $Q$. For example, if we consider $Q\equiv\reals^{3n}$, the configuration space
of $n$ particles in three-dimensional Euclidean space, conditions
\eqref{eq-lm-cons2}-\eqref{eq-am-cons2} are easily
verified. Conditions~\eqref{eq-lm-cons1}-\eqref{eq-am-cons1} depend not only on $Q$ but also on the
form of $\mbs{\mathsf{f}}^{\mathrm{int}}$ which has been, up to this point, left unspecified.  For
example, the canonical midpoint rule employs $\mbs{\mathsf{f}}^{\mathrm{int}}(\mbs{x},\mbs{y}) =
\mbs{f}^{\textrm{int}}((\mbs{x}+\mbs{y})/2)$, and preserves both linear and angular momenta, but not energy.  In
turn, the Energy-Momentum method \citep{Simo1992,Simo1994} provides an expression for this
force that guarantees strict energy conservation in the discrete update map, without upsetting the
preservation of momenta. Expanding on this idea, the Energy-Dissipative-Momentum-Conserving method
\citep{Armero2001_1,Armero2001_2} adds controllable energy dissipation to the solution, so small
that does not upset the accuracy of the solution, yet large enough that can damp out some of the
spurious oscillations in the high-frequency part of the solution.

\subsection{Discrete derivative}
\label{subs-dd}

As already mentioned, the direct evaluation of the internal forces
at the midpoint configuration does not guarantee, in general, the preservation
of energy. There exist however, consistent approximations of these
forces that strictly enforce this property of conservative equations.

To introduce the form of this ``conservative'' approximation of
the internal energy let us assume as in Section~\ref{sec-mechanics}
that the internal forces derive from a smooth potential $V:Q\to\reals$.
Be aware that from now on, we remove the superindex ${\textrm{int}}$, since no external force is longer considered along the derivation presented next.  
The type of approximations we search for are referred in the
literature as ``discrete derivatives'' \cite{Gonzalez1996} and are functions
$\mbs{\mathsf{f}}:Q\times Q\to\reals$ that satisfy, for
every $\mbs{x},\mbs{y}\in Q$, two properties, namely:

\begin{enumerate}[i.]
\item Directionality:
\begin{equation}
\langle\mbs{\mathsf{f}}(\bm{x},\bm{y}),\bm{y}-\bm{x}\rangle=V(\bm{y})-V(\bm{x})\:.
\end{equation}

\item Consistency:
\begin{equation}
  \mbs{\mathsf{f}}(\bm{x},\bm{x})
  =
  - DV(\mbs{x})
  =\bm{f}(\bm{x}) ,
  \label{eq-consistency}
\end{equation}
where $D$ denotes the standard derivative operator.
\end{enumerate}

When $Q\subset\reals$, there only exists one discrete derivative \citep{McLachlan1999}
and its closed form expression is given by
\begin{equation}
  \mathsf{f}(x,y) =
  \frac{V(y)-V(x)}{|y-x|}\ ,
  \label{eq-dd-oned}
\end{equation}
with the well-defined limit
\begin{equation}
  \lim_{y\to x} \mathsf{f}(x,y) = -DV(x) = f(x)\ .
  \label{eq-dd-limit}
\end{equation}
In higher dimensions, there are actually an infinite number of discrete derivatives
\cite{Romero2012,McLachlan1999} since only the component of $\mbs{\mathsf{f}}$ along the direction of
$\bm{y}-\bm{x}$ needs to have a precise value in order to guarantee energy conservation, and its
orthogonal complement is free to vary, as long as consistency of the approximation is preserved. This statement is formalized next:

\begin{theorem}
Any discrete derivative can be rewritten as
\begin{equation}
	\bm{\mathsf{f}}(\bm{x}, \bm{y}) = \frac{V(\bm{y})-V(\bm{x})}{\| \bm{y}-\bm{x} \|^2}(\bm{y}-\bm{x})+\bm{\mathsf{g}}(\bm{x}, \bm{y}),	
\label{eq-f-general}
\end{equation}
with $\bm{\mathsf{g}}(\bm{x}, \bm{y})$ a vector-valued function such that
\begin{equation}
	\left\langle \bm{\mathsf{g}}(\bm{x}, \bm{y}),  \bm{y}-\bm{x} \right\rangle = 0, \quad (\bm{y}\neq\bm{x}),  
\end{equation}
and
\begin{equation}
	\lim_{\bm{y}\to \bm{x}} \left( \bm{\mathsf{g}}(\bm{x}, \bm{y}) - \bm{\mathfrak{P}}_{(\bm{y}-\bm{x})}^\perp \bm{f}(\bm{x})\right) = \bm{0}, 
\end{equation}
where $\bm{\mathfrak{P}}_{(\bm{y}-\bm{x})}^\perp$ is the projection on the component perpendicular to $\bm{y}-\bm{x}$.
\end{theorem}
\begin{proof}
Let $\bm{\mathsf{g}}(\bm{x}, \bm{y})$ be defined as
\begin{equation}
\bm{\mathsf{g}}(\bm{x}, \bm{y}) = \bm{\mathsf{f}}(\bm{x}, \bm{y})-\frac{V(\bm{y})-V(\bm{x})}{\| \bm{y}-\bm{x} \|^2}(\bm{y}-\bm{x}),	
\end{equation}
It is apparent that $\bm{\mathsf{g}}(\bm{x}, \bm{y})$ is perpendicular to $\bm{y}-\bm{x}$ because
$\left\langle \bm{\mathsf{f}}(\bm{x}, \bm{y}),\bm{y}-\bm{x}\right\rangle =
V(\bm{y})-V(\bm{x})$, and
\begin{equation}
\begin{aligned}
	\bm{\mathsf{g}}(\bm{x}, \bm{y})\!-\!\bm{\mathfrak{P}}_{(\bm{y}-\bm{x})}^\perp \bm{f}(\bm{x}) 
	& = \bm{\mathsf{f}}(\bm{x}, \bm{y})\!-\!\frac{V(\bm{y})\!-\!V(\bm{x})}{\| \bm{y}\!-\!\bm{x} \|^2}(\bm{y}\!-\!\bm{x})\!-\!\bm{f}(\bm{x})+\frac{\left\langle \bm{f}(\bm{x}),\bm{y}\!-\!\bm{x}\right\rangle}{\| \bm{y}\!-\!\bm{x} \|^2}(\bm{y}\!-\!\bm{x})\\
	& = \bm{\mathsf{f}}(\bm{x}, \bm{y})\!-\!\bm{f}(\bm{x})-\frac{1}{\| \bm{y}\!-\!\bm{x} \|}\left(V(\bm{y})\!-\!V(\bm{x})\!-\!\left\langle \bm{f}(\bm{x}),\bm{y}\!-\!\bm{x}\right\rangle\right)\frac{(\bm{y}\!-\!\bm{x})}{\| \bm{y}\!-\!\bm{x} \|}, 	
\end{aligned}	
\end{equation}
which tends to zero as $\bm{y} \to \bm{x}$.
\end{proof}

\subsection{Conservative algorithmic force}

We explore next a type of discrete derivative that is different to the
one usually employed in nonlinear mechanics \cite{Simo1992,Gonzalez:2000wx}.
For that, we construct first a convex combination of the (exact) derivative
at two configurations, \ie
\begin{equation}
  \mbs{\mathsf{f}}^{\mathrm{cons}}
  =
  \frac{1}{2}(1-\alpha^{\mathrm{cons}})
  \bm{f}(\bm{x})
  +
  \frac{1}{2}(1+\alpha^{\mathrm{cons}})\bm{f}(\bm{y})
\end{equation}
or in a more compact form 
\begin{equation}
  \mbs{\mathsf{f}}^{\mathrm{cons}}
  =
  \bm{f}_{a} + \alpha^{\mathrm{cons}}
  \,
  \ddd{\mbs{f}}
  \,.
  \label{eq:conservative_part_current_0}
\end{equation}
In this expression, the scalar  $\alpha^{\mathrm{cons}}$ has to be determined in order to guarantee directionality,
$\bm{f}_{a}$ is the averaged
force
\begin{equation}
  \bm{f}_{a}
  =
  \frac{\bm{f}(\bm{x})+\bm{f}(\bm{y})}{2}\,,
\end{equation}
and $\ddd{\mbs{f}}$ is one
half of the force jump between the configurations at times $t_{n}$
and $t_{n+1}$, \ie
\begin{equation}
  \ddd{\mbs{f}}
  =
  \frac{\bm{f}(\bm{y})-\bm{f}(\bm{x})}{2}\,.
\end{equation} Notice that this conservative approximation satisfies, by construction, the
consistency condition~\eqref{eq-consistency}. The satisfaction of directionality depends, as
advanced, on the choice of the parameter $\alpha^{\mathrm{cons}}$. To enforce
it, we select $\alpha^{\mathrm{cons}}$ by means of an optimality condition~\citep{Romero2012},
namely, as the scalar that minimizes
\begin{equation}
  \begin{array}{ccc}
    & &
      \frac{1}{2}\|\mbs{\mathsf{f}}^{\mathrm{cons}}-\bm{f}_{m}\|_{\metric}^{2}\\
    \mathrm{subject~to} &
    & \langle\mbs{\mathsf{f}}^{\mathrm{cons}},\bm{y}
      -
      \bm{x}\rangle-V(\bm{y})+V(\bm{x})=0\,.
\end{array}
\end{equation}
Assuming $\mbs{f}(\mbs{x}) \ne \mbs{f}(\mbs{y})$, this
optimization problem is a linearly constrained quadratic program that can be solved
in closed form. Moreover, the only requirement for the optimization problem to be convex is that $\|\bm{f}(\bm{y})-\bm{f}(\bm{x})\|^2_{\metric} > 0$.
Its solution can be interpreted as
the discrete derivative that is closest to $\mbs{f}_m$, the continuous force at the midpoint
\textbf{$\bm{f}_{m}$}, namely,
\begin{equation}
\bm{f}_{m}=\bm{f}\left(\frac{\bm{x}+\bm{y}}{2}\right)\,.
\end{equation}
Here, $\metric$ is a metric tensor. The Lagrangian of the optimization problem is
\begin{equation}
  \mathcal{L}(\alpha^{\mathrm{cons}},\lambda^{\mathrm{cons}})
  =
  \frac{1}{2}\|\mbs{\mathsf{f}}^{\mathrm{cons}}-\bm{f}_{m}\|_{\metric}^{2}+
  \lambda^{\mathrm{cons}}
  \left(
    \langle\mbs{\mathsf{f}}^{\mathrm{cons}},\bm{y}-\bm{x} \rangle
    -
    V(\bm{y}) + V(\bm{x})
  \right)
  \,,
\end{equation}
where $\lambda^{\mathrm{cons}}$ is a Lagrange multiplier that enforces directionality.
To find the stationarity condition, the variation of $\mathcal{L}$ is calculated as:
\begin{equation}
\delta\mathcal{L}(\alpha^{\mathrm{cons}},\lambda^{\mathrm{cons}})={}
\langle\delta\bm{f}^{\mathrm{cons}},
\metric(\mbs{\mathsf{f}}^{\mathrm{cons}}-\bm{f}_{m})+
	\lambda^{\mathrm{cons}}(\bm{y}-\bm{x})\rangle
       + \delta\lambda^{\mathrm{cons}}
        \left(
        \langle\mbs{\mathsf{f}}^{\mathrm{cons}},\bm{y}-\bm{x}
        \rangle-
	V(\bm{y}) + V(\bm{x})
        \right) \,.
\end{equation}
Now for the sake of brevity, let us introduce a discrete function
defined as
\begin{equation}
  \tilde{\mathcal{C}}_{f}(\bm{x},\bm{y})
  =
  V(\bm{y})-V(\bm{x})
  -
  \left\langle \bm{f}_a , \bm{y}-\bm{x}\right\rangle \,.
\end{equation}
From now on, we refer to this function as a conservation function,
which allows the preservation energy in the discrete setting for a
given fixed time step $\Delta t$. This conservation function is not unique
and depends, in principle, on the shape of the approximated discrete
form.

The stationarity condition for the associated Lagrangian function can
be reformulated as the linear system that is explicitly given by
\begin{equation}
\left(\begin{array}{cc}
A_{11}^{f} & A_{12}^{f}\\
A_{12}^{f} & 0
\end{array}\right)\left[ \begin{array}{c}
\alpha^{\mathrm{cons}}\\
\lambda^{\mathrm{cons}}
\end{array}\right] =
\left[\begin{array}{c}
b_1^{f,\:\mathrm{cons}}\\
b_2^{f,\:\mathrm{cons}}
\end{array}\right]\, ,
\end{equation}
with
\begin{equation}
A_{11}^{f}=\frac{1}{2}\|\bm{f}(\bm{y})-\bm{f}(\bm{x})\|_{\metric}^{2}\:,
\end{equation}
\begin{equation}
A_{12}^{f}=\langle\bm{f}(\bm{y})-\bm{f}(\bm{x}),\bm{y}-\bm{x}\rangle\:,
\end{equation}
\begin{equation}
b_1^{f,\:\mathrm{cons}}=
\langle\bm{f}(\bm{y})-\bm{f}(\bm{x}),\metric(\bm{f}_{m}-\bm{f}_{a})\rangle
\ ,
\end{equation}
and
\begin{equation}
b_2^{f,\:\mathrm{cons}}=
2\tilde{\mathcal{C}}_{f}(\bm{x},\bm{y})\:.
\end{equation}
The solution of this program is
\begin{equation}
  \alpha^{\mathrm{cons}}=\frac{2\tilde{\mathcal{C}}_{f}(\bm{x},\bm{y})}{\langle\bm{f}
    (\bm{y})-\bm{f}(\bm{x}),\bm{y}-\bm{x}\rangle}
  \ ,
\end{equation}
and
\begin{equation}
  \lambda^{\mathrm{cons}} =
-\frac{\langle\bm{f}
  (\bm{y})-\bm{f}(\bm{x}),\metric(\bm{f}_{m}-\bm{f}_{a})
  \rangle}{\langle\bm{f}(\bm{y})-\bm{f}(\bm{x}),\bm{y}-\bm{x}\rangle}
-\frac{\|\bm{f}(\bm{y})-\bm{f}(\bm{x})\|
  _{\metric}^{2}\tilde{\mathcal{C}}_{f}
  (\bm{x},\bm{y})}{\langle\bm{f}(\bm{y})-\bm{f}(\bm{x}),\bm{y}-\bm{x}
  \rangle^{2}}\,.
\end{equation}
Notice that $\alpha^{\mathrm{cons}}$ does not depend on the chosen metric.
Then, we can claim that the adopted construction affords a unique definition. This feature
represents a main innovation of the current work.  However and up to this point, it is not clear to
which extent the current formula approaches the commonly used formulas like the one due to
Gonzalez \citep{Gonzalez1996} or the one due to Harten \emph{et al.}~\citep{Harten1983};
a comparison of that second method with the current one is beyond the scope of
this work.

In contrast, the Lagrange multiplier $\lambda^{\mathrm{cons}}$ depends on the chosen metric.
Finally, the conservative part of the discrete force takes the following explicit form:
\begin{equation}
  \mbs{\mathsf{f}}^{\mathrm{cons}}(\bm{x},\bm{y}) =  
  \frac{\bm{f}(\bm{x})+\bm{f}(\bm{y})}{2}+\\
  \frac{\tilde{\mathcal{C}}_{f}(\bm{x},\bm{y})}{\langle
    \bm{f}(\bm{y})-\bm{f}(\bm{x}),\bm{y}-\bm{x}\rangle}(\bm{f}(\bm{y})-\bm{f}(\bm{x}))\,.	
\label{eq-f-cons}
\end{equation}
In the context of nonlinear elastodynamics, the first
term of the formula is equivalent to the definition proposed
by Simo and Tarnow~\citep{Simo1992} that was derived in
the context of Saint Venant-Kirchhoff materials. The second term
can be interpreted as a correction for the most
general hyperelastic case. The formula proposed by Gonzalez~\citep{Gonzalez1996}
cannot be algebraically reduced to the proposed expression,
because the former is basically a correction for $\bm{f}_{m}$ and
the latter,  for $\bm{f}_{a}$.

The conserving force given by Eq.~\eqref{eq-f-cons} can be rewritten in
the form of Eq.~\eqref{eq-f-general} and therefore, is a discrete
derivative with
\begin{equation}
  \bm{\mathfrak{P}}_{\bm{y}-\bm{x}}^\parallel \bm{\mathsf{f}}^{\mathrm{cons}}(\bm{x}, \bm{y})
  = \frac{V(\bm{y})-V(\bm{x})}{\| \bm{y}-\bm{x} \|^2}(\bm{y}-\bm{x})\, ,
\end{equation}
where $\bm{\mathfrak{P}}_{\bm{y}-\bm{x}}^\parallel$ is the projection parallel to $\bm
{y}-\bm{x}$, and
\begin{equation}
  \bm{\mathsf{g}}^{\mathrm{cons}}(\bm{x}, \bm{y})
  := \bm{\mathfrak{P}}_{\bm{y}-\bm{x}}^\perp \bm{\mathsf{f}}^{\mathrm{cons}}(\bm{x}, \bm{y}).
\end{equation}

\subsection{Dissipative algorithmic force}
To account for dissipation, let us assume the existence of a dissipative
part of the algorithmic internal force that is proportional to
$\ddd{\bm{f}}$, this is
\begin{equation}
\bm{\mathsf{f}}^{\mathrm{diss}}
=
\alpha^{\text{diss}}
\ddd{\bm{f}}\,,
\end{equation}
where $\alpha^{\mathrm{diss}}$ is a scalar whose precise definition is still open. This construction
is supported by the analysis done by Romero \citep{Romero2012}, which showed that other choices may
destroy the accuracy of the approximation.  We will see later that this expression is very
attractive since it provides an unifying treatment of both conservative and dissipative parts of the
algorithmic internal force. 

To find the value of $\alpha^{\mathrm{diss}}$,
we define a discrete dissipation function $\tilde{\mathcal{D}}_{f}(\bm{x},\bm{y})$,
which must be positive semi-definite,
at least second order in $\| \bm{y}-\bm{x}\|$
to avoid spoiling the accuracy of the algorithm, and tend to $0$ as
$\bm{x}$ tends to $\bm{y}$. Then $\alpha^{\mathrm{diss}}$ can
be obtained as the scalar that minimizes
\begin{equation}
\begin{array}{ccc}
&  & \frac{1}{2}\|\bm{\mathsf{f}}^{\mathrm{diss}}\|_{\metric}^{2}\\
\mathrm{subject~to} &  & \langle\bm{\mathsf{f}}^{\mathrm{diss}},
\bm{y}-\bm{x}\rangle-\tilde{\mathcal{D}}_{f}(\bm{x},\bm{y})=0\,.
\end{array}
\end{equation}
This is also a linearly constrained quadratic program.
The solution of this optimization problem can be interpreted as the
smallest perturbation force that satisfies the dissipation
relation $\langle \mbs{y}- \mbs{x},\bm{\mathsf{f}}^{\mathrm{diss}}\rangle
=\mathcal{\tilde{D}}_{f}(\bm{x},\bm{y})$.
Once again, $\metric$ is a given metric tensor and the associated Lagrangian is simply
\begin{equation}
\mathcal{L}(\alpha^{\mathrm{diss}},\lambda^{\mathrm{diss}})
= 
\frac{1}{2}\|\bm{\mathsf{f}}^{\mathrm{diss}}\|_{\metric}^{2}
+
\lambda^{\mathrm{diss}}
(
\langle\bm{\mathsf{f}}^{\mathrm{diss}},\bm{y}-\bm{x}\rangle-
\tilde{\mathcal{D}}_{f}(\bm{x},\bm{y})
)\,,
\end{equation}
where $\lambda^{\mathrm{diss}}$ is a Lagrange multiplier that enforces
the dissipation constraint. To formulate the stationarity condition,
the variation of the associated Lagrangian has to be computed.
This procedure yields
\begin{equation}
\delta\mathcal{L}(\alpha^{\mathrm{diss}},\lambda^{\mathrm{diss}})
= 
\langle
\delta\bm{\mathsf{f}}^{\mathrm{diss}},\metric\bm{\mathsf{f}}^{\mathrm{diss}}
+
\lambda^{\mathrm{diss}}(\bm{y}-\bm{x})
\rangle
+
\delta\lambda^{\mathrm{diss}}
(
\langle\bm{\mathsf{f}}^{\mathrm{diss}},\bm{y}-\bm{x}\rangle-
\tilde{\mathcal{D}}_{f}(\bm{x},\bm{y})
)\,.
\end{equation}
Noting that $\delta \bm{\mathsf{f}}^{\mathrm{diss}} = \delta\alpha^{\mathrm{diss}}
\ddd{\mbs{f}}$, the stationarity condition of the Lagrangian
can be written explicitly as
\begin{equation}
\left(\begin{array}{cc}
A_{11}^{f} & A_{12}^{f}\\
A_{12}^{f} & 0
\end{array}\right)\left[ \begin{array}{c}
\alpha^{\mathrm{diss}}\\
\lambda^{\mathrm{diss}}
\end{array}\right] =\left[\begin{array}{c}
0\\
b_{2}^{f,\:\mathrm{diss}}
\end{array}\right] \,,
\end{equation}
with
\begin{equation}
b_{2}^{f,\:\mathrm{diss}}=2\tilde{\mathcal{D}}_{f}(\bm{x},\bm{y})\:.
\end{equation}
The solution of this linearly constrained quadratic program is
\begin{equation}
\alpha^{\mathrm{diss}}=
\frac{2\;\tilde{\mathcal{D}}_{f}(\bm{x},\bm{y})}{
	\langle\bm{f}(\bm{y})-\bm{f}(\bm{x}),\bm{y}-\bm{x}\rangle}
\ ,
\end{equation}
and
\begin{equation}
\lambda^{\mathrm{diss}}=-\frac{\| \bm{f}(\bm{y})-\bm{f}(\bm{x})\|
	_{\metric}^{2}\tilde{\mathcal{D}}_{f}(\bm{x},\bm{y})}{\langle\bm{f}(\bm{y})-\bm{f}(\bm{x}),\bm{y}-\bm{x}\rangle^{2}}\,.
\end{equation}
As in the case of the conservative part of the algorithmic force,
the parameter $\alpha^{\mathrm{diss}}$ does not depend on the chosen
metric and therefore, it is unique. However and up to this point, it is
not clear to which extent the current formula approaches already established
formulas, especially the one due to {Armero and Romero} \citep{Armero2001_1,Armero2001_2}.
As in the case of the conservative part of the approximation,
the multiplier $\lambda^{\mathrm{diss}}$
depends on the chosen metric. Finally, the formula for the dissipative
part of the algorithmic internal force takes the following explicit
form:
\begin{equation}
\bm{\mathsf{f}}^{\mathrm{diss}}=
\frac{\tilde{\mathcal{D}}_{f}(\bm{x},\bm{y})}
{\langle
	\bm{f}(\bm{y})-\bm{f}(\bm{x}),\bm{y}-\bm{x}\rangle}(\bm{f}(\bm{y})-\bm{f}(\bm{x}))\,.
\label{eq-f-diss}
\end{equation}

Notice that this formula has the same structure as the second term
of the conservative part of the algorithmic force. However, instead of
the conservation function, the dissipation function appears in the numerator.
This fact suggests that a unifying formula containing both conservative
and dissipative parts is possible, which makes the approach very attractive. 

\subsection{Generalization and preservation of momenta}
Assuming an additive composition of the algorithmic approximation
of the force, that is,
\begin{equation}
\bm{\mathsf{f}}
=
\bm{\mathsf{f}}^{\mathrm{cons}} + \bm{\mathsf{f}}^{\mathrm{diss}}\,,
\end{equation}
we can write, using Eqs.~\eqref{eq-f-cons} and~\eqref{eq-f-diss},
\begin{equation}
\bm{\mathsf{f}}(\bm{x},\bm{y})
=
\frac{\bm{f}(\bm{x})+\bm{f}(\bm{y})}{2}+\\
\frac{\tilde{\mathcal{C}}_{f}(\bm{x},\bm{y})+\tilde{\mathcal{D}}_{f}
	(\bm{x},\bm{y})}{\langle\bm{f}(\bm{y})-\bm{f}(\bm{x}),\bm{y}-\bm{x}\rangle}
(\bm{f}(\bm{y})-\bm{f}(\bm{x}))
\,.
\label{eq:f-sum}
\end{equation}

This formula is very compact and a simple inspection confirms that when the dissipation is zero, the
directionality condition is exactly verified.

To accommodate the preservation of linear and angular momenta as discussed in
Eqs. \eqref{eq-lm-conservation} and \eqref{eq-am-conservation}, Eq. \eqref{eq:f-sum} must be
modified as indicated next: Let $G$ be a Lie group with algebra $\mathfrak{g}$ and coalgebra
$\mathfrak{g}^{*}$, which acts on the configuration space $Q\subseteq\reals^{3n}$ by means of the
action $\bm{\chi}:G\times Q\to Q$. For every $\bm{\xi}\in\mathfrak{g}$, let $\bm{\xi}_Q: Q\to TQ$ denote
the infinitesimal generator of the action. Following again Gonzalez \citep{Gonzalez1996}, we can
define $G$-equivariant derivatives.  If $V:Q\to\reals$ is a $G$-invariant function, its
$G$-invariant discrete derivative is a smooth map $\bm{\mathsf{f}}^G:Q\times Q\to\reals$ that
satisfies the requirements of discrete derivatives and, moreover, the equivariance and orthogonality
condition, namely,
\begin{equation}
  \bm{\mathsf{f}}^G(\bm{\chi}_g(\bm{x}),\bm{\chi}_g(\bm{y}))
  =\left(\mathsf{D}\bm{\chi}_g\left(\frac{\bm{x}+\bm{y}}{2}\right)\right)^{-T}\bm{\mathsf{f}}^G(\bm{x},\bm{y})\, ,
\end{equation}
for all $\bm{x}, \bm{y} \in Q$, $g\in G$, and
\begin{equation}
\bm{\mathsf{f}}^G(\bm{x},\bm{y})\cdot\bm{\xi}_{\bm{Q}}\left(\frac{\bm{x}+\bm{y}}{2}\right)=0\, .
\end{equation}

To construct a $G$-equivariant discrete derivative, consider invariant functions under the symmetry
action denoted by $\pi_i$, $i=1,2,....,q$ where $q$ is the dimension of the quotient space
$Q/G$. Let $\bm{\Pi}=(\pi_1, \pi_2,...,\pi_q)$. If $V:Q\to\reals$ is $G$-invariant, a reduced
function $\tilde{V}$ can be defined by the relation $V=\tilde{V}\circ\bm{\Pi}$. If each of the
invariants is at most of degree two, then a $G$-equivariant discrete derivative for $V$ can be
constructed as
\begin{equation}
\bm{\mathsf{f}}^G(\bm{x}, \bm{y}) = \tilde{\bm{\mathsf{f}}}(\bm{\Pi}(\bm{x}), \bm{\Pi}(\bm{y}))\circ D\bm{\Pi}\left(\frac{\bm{x}+\bm{y}}{2}\right)\, ,
\label{eq-f-G-equi}
\end{equation}
where, as before $\bm{x}, \bm{y} \in Q$.
In particular, the formulation of a $G$-equivariant discrete derivative that
preserves linear and angular momenta (\cf~Eqs.~\eqref{eq-lm-conservation}
and \eqref{eq-am-conservation}) is straightforward.

The expression of the $G$-equivariant force in the current context is given by
\begin{equation}
\bm{\mathsf{f}}
(\bm{x},\bm{y})=
D \mbs{\Pi}^T(\mbs{z})
\left( \frac{\bm{f}(\bm{\Pi}(\bm{y}))+\bm{f}(\bm{\Pi}(\bm{x}))}{2}+\\
\alpha
(\bm{\Pi}(\bm{x}),(\bm{\Pi}(\bm{y}))
(\bm{f}(\bm{\Pi}(\bm{y}))-\bm{f}(\bm{\Pi}(\bm{x})))\right)\ ,
\label{eq-f-objective}
\end{equation}
with $\bm{z}=(\bm{x}+\bm{y})/2$ and
\begin{equation}
\alpha
(\bm{\Pi}(\bm{x}),(\bm{\Pi}(\bm{y})) = 
\frac{\tilde{\mathcal{C}}
_{f}(\bm{\Pi}(\bm{y}),\bm{\Pi}(\bm{x}))+\tilde{\mathcal{D}}
_{f}(\bm{\Pi}(\bm{y}),\bm{\Pi}(\bm{x}))}
{\langle\bm{f}(\bm{\Pi}(\bm{y}))-\bm{f}(\bm{\Pi}(\bm{x})),\bm{\Pi}(\bm{y})-\bm{\Pi}(\bm{x})\rangle}\, .
\end{equation}

\subsection{Interpretation of the conservative algorithmic force}
There exist infinite second order accurate approximations
of the midpoint force that 
that lead to energy and momentum conserving discretizations \citep{Romero2012}.
A general expression for the algorithmic force that is in agreement with the definition of the
discrete derivative is given by
\begin{equation}
  \bm{\mathsf{f}}_{\metric}(\bm{x},\bm{y})
  =\bm{f}\left(\frac{\bm{x}+\bm{y}}{2}\right)+\frac{\hat{\mathcal{C}}_{f}(\bm{x},\bm{y})}{\|\bm{y}-\bm{x}\| _{\metric^{-1}}^{2}}\metric^{-1}(\bm{y}-\bm{x})\,,
\label{eq:conservative_part_romero}
\end{equation}
where $\metric$ is the matrix representation of a suitable metric tensor and its associated
conservation function is
\begin{equation}
\hat{\mathcal{C}}_{f}(\bm{x},\bm{y})=\left[V(\bm{y})-V(\bm{x})\right]-\left\langle 
\bm{f}_m
,\bm{y}-\bm{x}\right\rangle \,.
\end{equation}

As shown in \citep{Romero2012}, the choice of the metric tensor in
Eq.~\eqref{eq:conservative_part_romero} is crucial, since a wrong choice can destroy the accuracy of the
solution even when the directionality and consistency properties are verified. The general
expression can be reduced to the original formula proposed in \citep{Gonzalez1996} just by adopting
the standard Euclidean metric tensor, \ie
\begin{equation}
\bm{\mathsf{f}}_{\bm{I}}(\bm{x},\bm{y})=\bm{f}\left(\frac{\bm{x}+\bm{y}}{2}\right)+\frac{\hat{\mathcal{C}}_{f}(\bm{x},\bm{y})}{\| \bm{y}-\bm{x}\| ^{2}_{\bm{I}}}(\bm{y}-\bm{x})\,.
\label{eq:conservative_part_gonzalez}
\end{equation}

Until now, this original formula has been regarded as the optimal
one from the implementation point of view. By visual inspection of
Eq.~\eqref{eq-f-cons}, it is also possible
to claim that the formula derived in this work is as easy to implement
as the original one given in Eq.~\eqref{eq:conservative_part_gonzalez}.
Moreover, the new formula requires only evaluations at the endpoints of the time
interval and not at the midpoint.

Next, we would like to analyze the formula~\eqref{eq-f-cons}
in terms of the general expression provided by Eq.~\eqref{eq:conservative_part_romero}.
First, defining $\mbs{z}$ to be the average $\mbs{z}=(\mbs{x}+\mbs{y})/2$,
we make use of Taylor's theorem to compute
\begin{equation}
\begin{aligned}
\bm{f}(\bm{x}) = 
& 
\bm{f}(\bm{z})-
\frac{1}{2}D\bm{f}\cdot(\bm{y}-\bm{x})+
\frac{1}{4}D^{2}\bm{f}\cdot((\bm{y}-\bm{x}),(\bm{y}-\bm{x}))\\
& 
-\frac{1}{8}D^{3}\bm{f}\cdot((\bm{y}-\bm{x}),(\bm{y}-\bm{x}),(\bm{y}-\bm{x}))
+\mathcal{O}(\| \bm{y}-\bm{x}\| ^{4})
\end{aligned}
\end{equation}
and
\begin{equation}
\begin{aligned}
\bm{f}(\bm{y}) = 
&
\bm{f}(\bm{z})+
\frac{1}{2}D\bm{f}\cdot(\bm{y}-\bm{x})+
\frac{1}{4}D^{2}\bm{f}\cdot((\bm{y}-\bm{x}),(\bm{y}-\bm{x}))\\
&
+\frac{1}{8}D^{3}\bm{f}\cdot((\bm{y}-\bm{x}),(\bm{y}-\bm{x}),(\bm{y}-\bm{x}))
+\mathcal{O}(\| \bm{y}-\bm{x}\| ^{4})\,.
\end{aligned}
\end{equation}
The averaged force can be expressed as
\begin{equation}
\frac{\bm{f}(\bm{x})+\bm{f}(\bm{y})}{2}=\bm{f}(\bm{z})+\mathcal{O}(\| \bm{y}-\bm{x}\| ^{2})
\end{equation}
and
\begin{equation}
\left\langle \frac{\bm{f}(\bm{x})+\bm{f}(\bm{y})}{2},\bm{y}-\bm{x}\right\rangle ={}
\langle\bm{f}(\bm{z}),\bm{y}-\bm{x}\rangle+
\mathcal{O}(\| \bm{y}-\bm{x}\| ^{3})\,.
\end{equation}
The force jump can be written as
\begin{equation}
\bm{f}(\bm{y})-\bm{f}(\bm{x})=\langle D\bm{f},\bm{y}-\bm{x}\rangle+\mathcal{O}(\| \bm{y}-\bm{x}\| ^{3}),
\end{equation}
and
\begin{equation}
\begin{aligned}
\langle\bm{f}(\bm{y})-\bm{f}(\bm{x}),\bm{y}-\bm{x}\rangle
&=\left\langle \bm{y}-\bm{x},
  D^2V(\mbs{z})
  (\bm{y}-\bm{x})\right\rangle+\mathcal{O}(\|\bm{y}-\bm{x}\|^{4})\\
&=\left\langle \bm{y}-\bm{x},
D^2V(\mbs{x})
  (\bm{y}-\bm{x})\right\rangle+\mathcal{O}(\|\bm{y}-\bm{x}\|^{3})\,.
\end{aligned}
\end{equation}
Now putting everything together, we can rewrite Eq.\eqref{eq-f-cons} as
\begin{equation}
\bm{\mathsf{f}}^{\mathrm{cons}} =
\bm{f}(\mbs{z})
+
\frac{V(\bm{y})-V(\bm{x})-\left\langle \bm{f}(\mbs{z}),\bm{y}-\bm{x}\right\rangle}{
  \left\langle \bm{y}-\bm{x}, D^2V(\mbs{x}) (\bm{y}-\bm{x})\right\rangle }
D^2 V(\mbs{x})(\bm{y}-\bm{x})+\mathcal{O}(\| \bm{y}-\bm{x}\| ^{2})\,.
\end{equation}

Taking a look at this expression, it is apparent that the discrete
force~\eqref{eq-f-cons} is a second order
perturbation of the midpoint approximation and that the metric employed
for the definition of the conserving correction is just
\begin{equation}
  \metric=
  (D^2V (\mbs{x}))^{-1}\ .
\label{eq-metric-d2V}
\end{equation}
We conclude that the an integration scheme based on Eq. (\ref{eq:conservative_part_current_0})
would behave locally in a very similar manner to a method based on
Eq.~\eqref{eq:conservative_part_romero}
with metric~\eqref{eq-metric-d2V}. Their global behavior can, in general, differ.

\subsection{Dissipative algorithmic velocity}

As proposed in \citep{Armero2001_1,Armero2001_2},
the generalized velocity can also be expressed as the
linear combination of a conservative and a dissipative component.
If the mass matrix is configuration-independent, the midpoint
rule provides precisely the conservative part of the velocity.
Following
the ideas adopted for the formulation of the dissipative part of the algorithmic force,
we find next the smallest perturbation of the midpoint velocity that
guarantees dissipation
according to a given dissipation function $\tilde{\mathcal{D}}_{s}(\bm{u},\bm{v})$,
which must be non-negative, at least second order accurate
in $\|\bm{v}-\bm{u}\|$, and tend to $0$ as $\bm{u}$
tends to $\bm{v}$.

According to the midpoint rule, the conservative part
of the algorithmic velocity can be expressed as
\begin{equation}
\bm{\mathsf{s}}^{\mathrm{cons}}=\frac{\bm{u}+\bm{v}}{2} = \bm{w}\:,
\end{equation}
where the velocity $\bm{u}$ corresponds to time instant $t_{n}$,
the velocity $\bm{v}$ corresponds to time instant $t_{n+1}$ and $\bm{w}$ is the averaged velocity.
This particular choice preserves linear and angular momenta. The preservation
of energy is guaranteed when the relation
\begin{equation}
\langle\bm{\mathsf{s}}^{\mathrm{cons}},\bm{M}(\bm{v}-\bm{u})\rangle=T(\bm{v})-T(\bm{u})
\end{equation}
is satisfied, with $T$ being the kinetic energy. Equivalently, this expression can be obtained from
$\langle\delta\bm{s}_{n+1/2},
\bm{\pi}(\bm{\mathsf{s}}^{\mathrm{cons}})\rangle=T(\bm{v})-T(\bm{u})$
when the variation is of the form $\delta\bm{s}_{n+1/2}=\bm{v}-\bm{u}$.
It is apparent that the conservative part of the algorithmic velocity
adopted fulfills this condition without need for further corrections.

Now we follow an idea that is slightly different to the one previously
used to derive the dissipative part of the algorithmic force. Henceforth,
let us assume the existence of a dissipative part of the algorithmic
velocity proportional to $\bm{\mathsf{s}}^{\mathrm{cons}}$,
this is
\begin{equation}
\bm{\mathsf{s}}^{\mathrm{diss}}=\beta^{\text{diss}}\bm{\mathsf{s}}^{\mathrm{cons}}\,,
\end{equation}
where $\beta^{\textrm{diss}}$ is a scalar to be found that
must guaranteed the dissipation of energy according to a given dissipation function
$\tilde{\mathcal{D}}_{s}(\bm{u},\bm{v})$.
This scalar can be obtained as the minimizer of 
\begin{equation}
\begin{array}{ccc}
&  & \frac{1}{2}\|\bm{\mathsf{s}}^{\mathrm{diss}}\|_{\bm{M}}^{2}\\
  \mathrm{subject~to} &
   & \langle\bm{\mathsf{s}}^{\mathrm{diss}},\bm{M}(\bm{v}-\bm{u})
     \rangle-\tilde{\mathcal{D}}_{s}(\bm{u},\bm{v})=0\, ,
\end{array}
\label{eq-beta-program}
\end{equation}
where $\mbs{M}$ is the mass matrix. Eq.~\eqref{eq-beta-program} defines a quadratic program with
linear constraints.  Its solution can be interpreted as the smallest non-conservative velocity
perturbation that satisfies for a discrete variation of the form
$\delta\bm{s}_{n+1/2}=\bm{v}-\bm{u}$, an energy dissipation according to the adopted rule.
The associated Lagrangian function of the optimization problem is
\begin{equation}
\mathcal{L}(\beta^{\textrm{diss}},\mu^{\textrm{diss}})=
\frac{1}{2}\|\bm{\mathsf{s}}^{\mathrm{diss}}\|_{\bm{M}}^{2}+
\mu^{\textrm{diss}}(\langle\bm{\mathsf{s}}^{\mathrm{diss}},\bm{M}(\bm{v}-\bm{u})\rangle-
\tilde{\mathcal{D}}_{s}(\bm{u},\bm{v}))\,,
\end{equation}
where $\mu^{\textrm{diss}}$ is a Lagrange multiplier that enforces
the dissipation constraint. To formulate the stationarity condition,
the variation of the associated Lagrangian function has to be computed.
This procedure yields 
\begin{equation}
\delta\mathcal{L}(\beta^{\textrm{diss}},\mu^{\textrm{diss}})=
\langle\delta\bm{\mathsf{s}}^{\mathrm{diss}},\bm{M}\bm{\mathsf{s}}^{\mathrm{diss}}+
\mu^{\textrm{diss}}\bm{M}(\bm{v}-\bm{u})\rangle+
\delta\mu^{\textrm{diss}}(\langle\bm{\mathsf{s}}^{\mathrm{diss}},\bm{M}(\bm{v}-\bm{u})\rangle-
\tilde{\mathcal{D}}_{s}(\bm{u},\bm{v}))\, .
\end{equation}
Noting that $\delta\bm{\mathsf{s}}^{\mathrm{diss}}=\delta\beta^{\mathrm{diss}}\bm{\mathsf{s}}^{\mathrm{cons}}$, the stationarity condition of the Lagrangian can be written explicitly as
\begin{equation}
\left(\begin{array}{cc}
A_{11}^{s} & A_{12}^{s}\\
A_{12}^{s} & 0
\end{array}\right)\left[\begin{array}{c}
\beta^{\textrm{diss}}\\
\mu^{\textrm{diss}}
\end{array}\right] =\left[\begin{array}{c}
0\\
b_{2}^{s,\:\textrm{diss}}
\end{array}\right] \, ,
\end{equation}
with
\begin{equation}
A_{11}^{s} = 2T(\bm{w})\:,
\end{equation}
\begin{equation}
A_{12}^{s} = T(\bm{v})-T(\bm{u})\:,
\end{equation}
and
\begin{equation}
b_{2}^{s,\:\textrm{diss}}=\tilde{\mathcal{D}}_{s}(\bm{u},\bm{v})\:.
\end{equation}
The solution of this linearly constrained quadratic program is
\begin{equation}
\beta^{\textrm{diss}}=\frac{\tilde{\mathcal{D}}_{s}(\bm{u},\bm{v})}{T(\bm{v})-T(\bm{u})}
\end{equation}
and
\begin{equation}
\mu^{\textrm{diss}}=-\frac{T(\bm{w})\tilde{\mathcal{D}}_{s}(\bm{u},\bm{v})}{(T(\bm{v})-T(\bm{u}))^2}\,.
\end{equation}
Finally, the formula for the dissipative
part of the algorithmic velocity takes the following explicit form:
\begin{equation}
  \bm{\mathsf{s}}^{\mathrm{diss}} =
  \frac{\tilde{\mathcal{D}}_{s}(\bm{u},\bm{v})}{T(\bm{v})-T(\bm{u})}\;
  \frac{\bm{v}+\bm{u}}{2}\,.
\end{equation}
Assuming an additive composition of the algorithmic approximation of the velocity, that is,
\begin{equation}
\bm{\mathsf{s}}=\bm{\mathsf{s}}^{\mathrm{cons}}+\bm{\mathsf{s}}^{\mathrm{diss}}\, ,
\end{equation}
we can write
\begin{equation}
  \bm{\mathsf{s}}(\bm{u},\bm{v})=\left(1+\frac{\tilde{\mathcal{D}}_{s}(\bm{u},\bm{v})}{T(\bm{v})-T(\bm{u})}\right)
  \frac{\bm{v}+\bm{u}}{2}\,.
\end{equation}
This formula is identical to the formula proposed in
\citep{Armero2001_1,Armero2001_2} that was derived by employing only geometric arguments. This time,
it can be clearly interpreted as an optimal approximation.

\subsection{Final equations}

The combination of all ingredients discussed here yields the full
discrete formulation of the dynamic equilibrium for nonlinear
mechanical systems. These consist of two residuals, one for the
generalized velocities and another for the generalized coordinates,
namely,
\begin{equation}
\left[\begin{array}{c}
\bm{r}_{s}\\
\bm{r}_{q}
\end{array}\right]_{n+1/2}
=
\left[ \begin{array}{c}
\bm{\pi}(\bm{s}_{n},\bm{s}_{n+1})-\bm{p}(\bm{q}_{n},\bm{q}_{n+1})\\
\dot{\bm{\pi}}(\bm{s}_{n},\bm{s}_{n+1})+\bm{\mathsf{f}}(\bm{q}_{n},\bm{q}_{n+1})-\bm{f}^{\mathrm{ext}}(\bm{q}_{n+1/2})
\end{array}\right]\, ,
\end{equation}
where both residuals have to be minimized at every time step. This
task is accomplished by means of a Newton-Raphson algorithm.

The generalized-velocity-based momentum term in its algorithmic form
is 
\begin{equation}
\bm{\pi}(\bm{s}_{n},\bm{s}_{n+1})=\bm{M}\left(
  1+\frac{\tilde{\mathcal{D}}_{s}(\bm{s}_{n},\bm{s}_{n+1})}{T(\bm{s}_{n+1})-T(\bm{s}_{n})}\right)
\frac{\bm{s}_{n+1}+\bm{s}_{n}}{2}\, .
\end{equation}
The generalized-coordinate-based momentum term in its discrete version
becomes 
\begin{equation}
  \bm{p}(\bm{q}_{n},\bm{q}_{n+1})
  =
  \bm{M}
  \frac{\bm{q}_{n+1}-\bm{q}_{n}}{\Delta t_n}\, .
\end{equation}
The generalized-coordinate-based momentum rate term in its algorithmic
form is 
\begin{equation}
  \dot{\bm{\pi}}(\bm{s}_{n},\bm{s}_{n+1})
  =
  \bm{M}
  \frac{\bm{s}_{n+1}-\bm{s}_{n}}{\Delta t_n}
  \, .
\end{equation}
Finally, the generalized internal force becomes
\begin{equation}
\bm{\mathsf{f}}(\bm{q}_{n},\bm{q}_{n+1})=
\frac{\bm{f}(\bm{q}_{n})+\bm{f}(\bm{q}_{n+1})}{2}+
\frac{\tilde{\mathcal{C}}_{f}(\bm{q}_{n},\bm{q}_{n+1})+\tilde{\mathcal{D}}_{f}(\bm{q}_{n},\bm{q}_{n+1})}{\langle\bm{f}(\bm{q}_{n+1})-\bm{f}(\bm{q}_{n}),\bm{q}_{n+1}-\bm{q}_{n}\rangle}
(\bm{f}(\bm{q}_{n+1})-\bm{f}(\bm{q}_{n}))\, .
\end{equation}
In the case of accommodating
the preservation of linear and angular momenta, we require the $G$-equivariant version given by 
\begin{equation}
\bm{\mathsf{f}}^G(\bm{q}_{n},\bm{q}_{n+1})=
D \mbs{\Pi}^T_{n+1/2}
\!\! \left( \frac{\bm{f}(\bm{\Pi}_{n+1})+\bm{f}(\bm{\Pi}_{n})}{2}+
\alpha
(\bm{\Pi}_{n},\bm{\Pi}_{n+1})
(\bm{f}(\bm{\Pi}_{n+1})-\bm{f}(\bm{\Pi}_{n}))\right)\,
\label{eq-f-objective2}
\end{equation}
where
\begin{equation}
\alpha
(\bm{\Pi}_{n},\bm{\Pi}_{n+1}) = 
\frac{\tilde{\mathcal{C}}
	_{f}(\bm{\Pi}_{n+1},\bm{\Pi}_{n})+\tilde{\mathcal{D}}
	_{f}(\bm{\Pi}_{n+1},\bm{\Pi}_{n})}
{\langle\bm{f}(\bm{\Pi}_{n+1})-\bm{f}(\bm{\Pi}_{n}),\bm{\Pi}_{n+1}-\bm{\Pi}_{n}\rangle}\,
\end{equation}
and $\bm{\Pi}_{n}=\bm{\Pi}(\bm{q}_{n})$. The discrete conservation function is given by
\begin{equation}
\tilde{\mathcal{C}}_{f}(\bm{q}_{n},\bm{q}_{n+1})=
(V(\bm{q}_{n+1})-V(\bm{q}_{n}))-\left\langle \frac{\bm{f}(\bm{q}_{n})+\bm{f}(\bm{q}_{n+1})}{2},\bm{q}_{n+1}-\bm{q}_{n}\right\rangle \, ,
\end{equation}
and
its $G$-equivariant version given by 
\begin{equation}
\tilde{\mathcal{C}}^G_{f}(\bm{q}_{n},\bm{q}_{n+1})=
(V(\bm{\Pi}_{n+1})-V(\bm{\Pi}_{n}))-\left\langle \frac{\bm{f}(\bm{\Pi}_{n})+\bm{f}(\bm{\Pi}_{n+1})}{2},\bm{\Pi}_{n+1}-\bm{\Pi}_{n}\right\rangle \, .
\end{equation}

The most basic discrete dissipation function at the level of the generalized
internal force that can be chosen
is
\begin{equation}
\tilde{\mathcal{D}}_{f}(\bm{q}_{n},\bm{q}_{n+1}, \tilde{\bm{q}}_{n})=\frac{1}{2}\left\langle \bm{q}_{n+1}-\bm{q}_{n}, \bm{D}(\tilde{\bm{q}}_{n}-\bm{q}_{n})\right\rangle\;,
\end{equation}
or 
its $G$-equivariant counterpart expressed as
\begin{equation}
\tilde{\mathcal{D}}^G_{f}(\bm{q}_{n},\bm{q}_{n+1}, \tilde{\bm{q}}_{n})=\frac{1}{2}\left\langle \bm{\Pi}_{n+1}-\bm{\Pi}_{n}, \bm{D}(\tilde{\bm{\Pi}}_{n}-\bm{\Pi}_{n})\right\rangle\;,
\end{equation}
where $\tilde{\bm{q}}_{n}$
could correspond to an intermediate configuration, and $\bm{D}$ is constant, symmetric and positive
semi-definite. The dissipation function for the velocity is of the form
\begin{equation}
\tilde{\mathcal{D}}_{s}(\bm{s}_{n},\bm{s}_{n+1}, \tilde{\bm{s}}_{n})=\frac{1}{2}\left\langle \bm{s}_{n+1}-\bm{s}_{n}, \bm{M}(\tilde{\bm{s}}_{n}-\bm{s}_{n})\right\rangle\;,
\end{equation}
or 
its $G$-equivariant version given by 
\begin{equation}
\tilde{\mathcal{D}}^G_{s}(\bm{s}_{n},\bm{s}_{n+1}, \tilde{\bm{s}}_{n})
=\frac{1}{2}(\|\bm{s}_{n+1}\|_{\bm{M}}-\|\bm{s}_{n}\|_{\bm{M}})( \|\tilde{\bm{s}}_{n}\|_{\bm{M}}-\|\bm{s}_{n}\|_{\bm{M}}))\ ,
\end{equation}
where $\tilde{\bm{s}}_{n}$ could correspond to an intermediate configuration. With this setting, unconditional stability
in the nonlinear sense can be achieved. 
The chosen dissipation functions correspond to those proposed in \citep{Armero2001_1,Armero2001_2} for the EDMC-1/2.
\section{Numerical results}
\label{sec-results}

In this section, we present four numerical examples which were chosen to show the potentialities of
the proposed approach. With these, we do not pretend to test the new approach exhaustively, but at
rather provide some insight on its properties. For this purpose, we study first two examples involving two-mass
systems with potential functions that can arise in the context of reduced-order models,
and then two examples of nonlinear elastic shell structures employing a neo-Hookean material.
\revone{Additionally, we briefly discuss the dissipation properties of the proposed scheme in the high-frequency range.}

\subsection{Reduced-order models}

The first example is a mechanical system with two degrees of freedom whose potential function
possesses polynomial complexity. The second one considers another mechanical system with also two
degrees of freedom, but whose potential function shows non-polynomial complexity. Here, we adopt the
most basic discrete dissipation functions at the level of the generalized internal force and at the
level of the generalized velocity that are given by
\begin{equation}
  \tilde{\mathcal{D}}_{f}(\bm{q}_{n},\bm{q}_{n+1})
  =\frac{\chi_{f}}{2h}\bigl\Vert\bm{q}_{n+1}-\bm{q}_{n}\bigr\Vert_{\bm{D}}^{2}\;,
\end{equation}
in which $\bm{D}$ is a constant symmetric semi-positive definite matrix, and
\begin{equation}
  \tilde{\mathcal{D}}_{s}(\bm{s}_{n},\bm{s}_{n+1})
  =\frac{\chi_{s}}{h}\left(\sqrt{T(\bm{s}_{n+1})}-\sqrt{T(\bm{s}_{n})}\right)^{2}\:,
\end{equation}
where $\chi_{f}$ and $\chi_{s}$ in $\reals_{\geq 0}$ are merely dissipation parameters. 
For both examples, four cases are considered: \emph{i}) fully conservative, \ie\
$\chi_{f}=\chi_{s}=0$; \emph{ii}) dissipative at the level of the generalized internal force, \ie\
$\chi_{f}\neq0$ and $\chi_{s}=0$; \emph{iii}) dissipative at the level of the generalized
velocities, \ie\ $\chi_{f}=0$ and $\chi_{s}\neq0$; and, \emph{iv}) fully dissipative, \ie\
$\chi_{f}\neq0$ and $\chi_{s}\neq0$. Additionally, to numerically gain some insight about the
accuracy of the method, we provide for these two examples and all cases the second quotient of
precision $Q_{\mathrm{II}}(t)$ computed on the basis of the corresponding states
\textbf{$\bm{\xi}\in Q\times S$}, namely \textbf{$\bm{\xi}=(\bm{q},\bm{s})$}. The definition of
$Q_{\mathrm{II}}(t)$ is presented in the Appendix.

\subsubsection{Two-mass system with a polynomial potential}

Here, we consider a nonlinear oscillatory mechanical system with potential function

\begin{equation}
  V(\bm{q})=\frac{1}{2}V_{ab}^{\mathrm{II}}q^{a}q^{b}+\frac{1}{3}V_{abc}^{\mathrm{III}}q^{a}q^{b}q^{c}
  +
  \frac{1}{4}V_{abcd}^{\mathrm{IV}}q^{a}q^{b}q^{c}q^{d}\, .
\end{equation}

\begin{figure}
		\centering{}
		\includegraphics[width=0.5\textwidth]{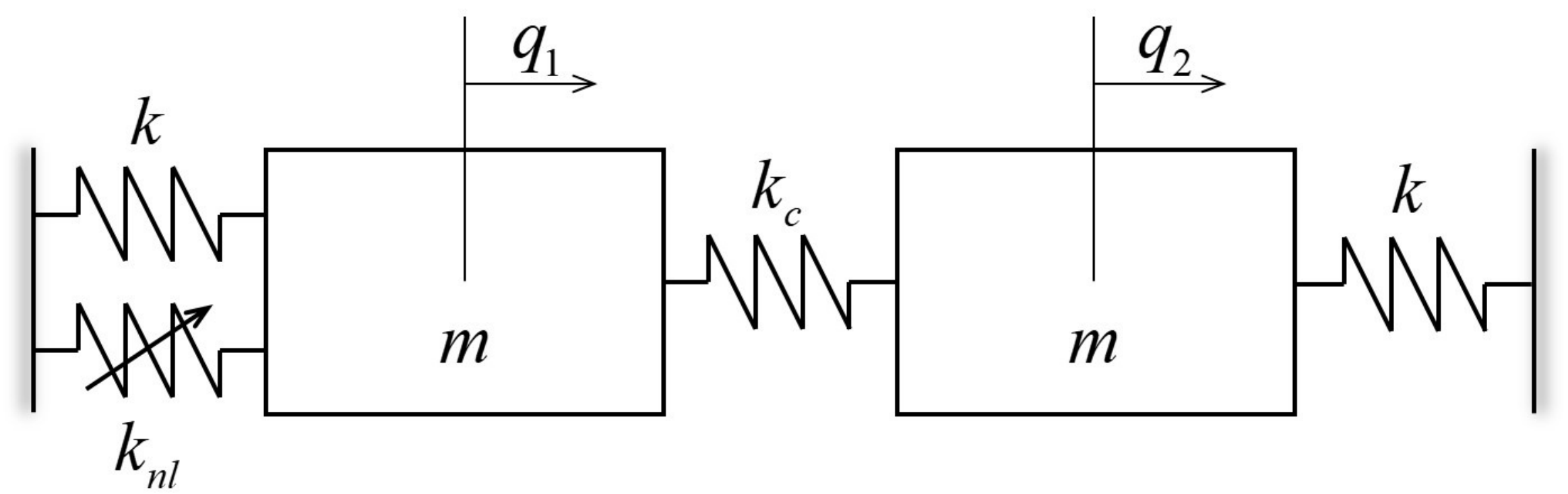}
		\caption{Two masses connected by linear springs between two walls. The first mass is also connected to the left wall through a nonlinear spring.}
		\label{fig_ex1_1}
\end{figure}

\begin{figure}
		\centering{}
		\includegraphics[width=0.5\textwidth]{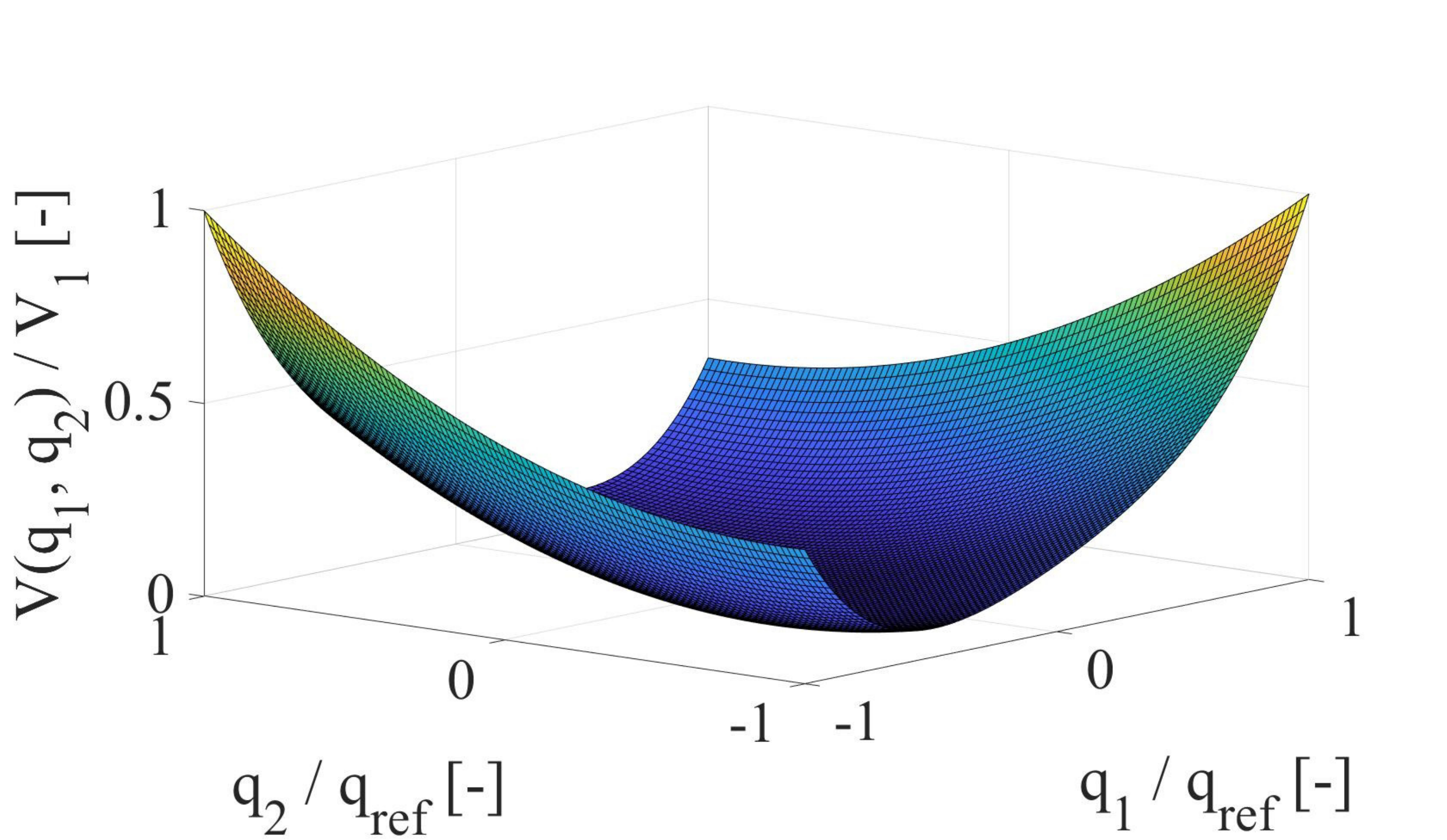}
		\caption{Potential function with polynomial complexity.}
		\label{fig_ex1_2}
\end{figure}

\begin{figure}
	\centering{}
	\begin{tabular}{cc}
		\includegraphics[width=0.45\textwidth]{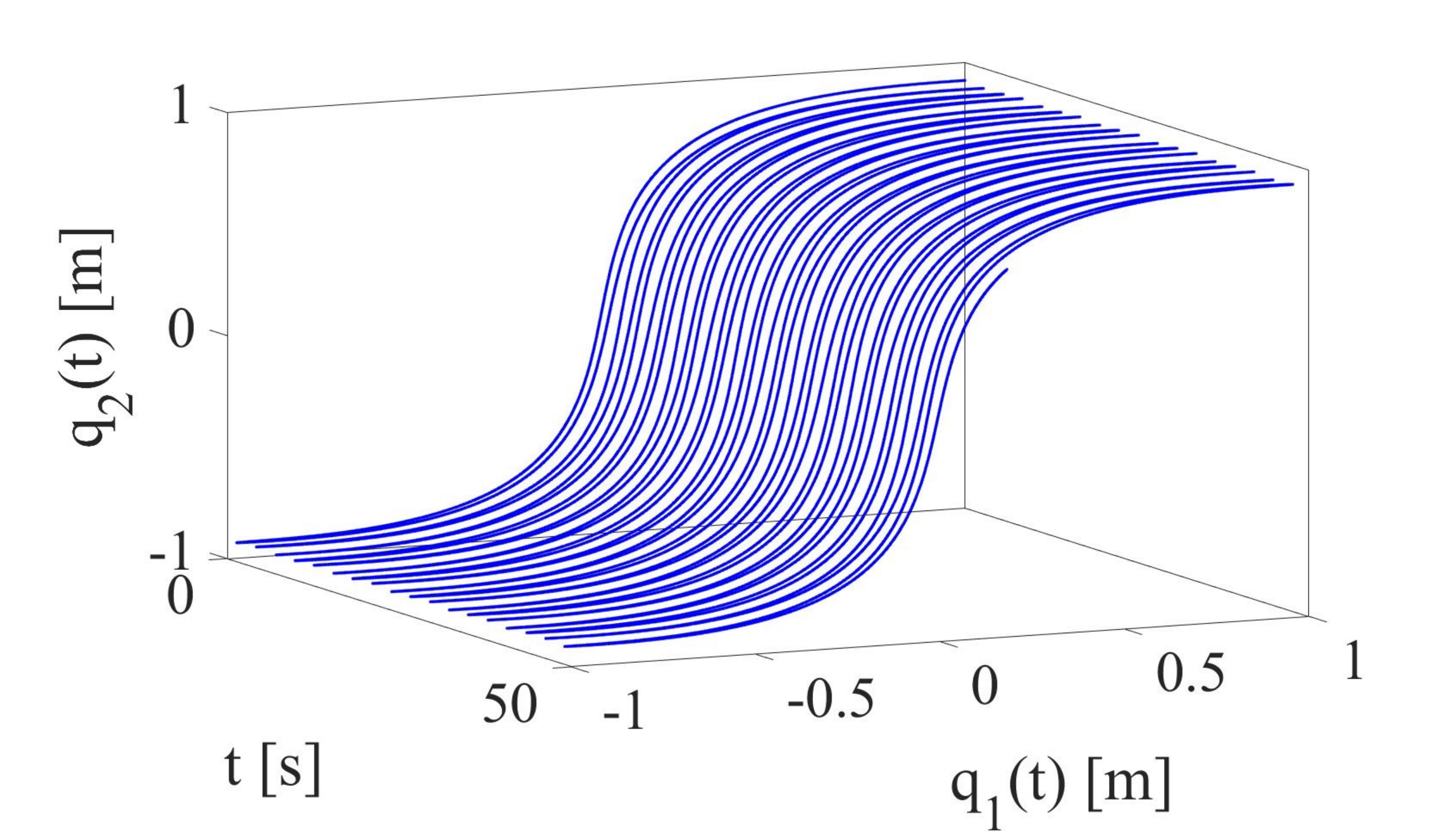} &	\includegraphics[width=0.45\textwidth]{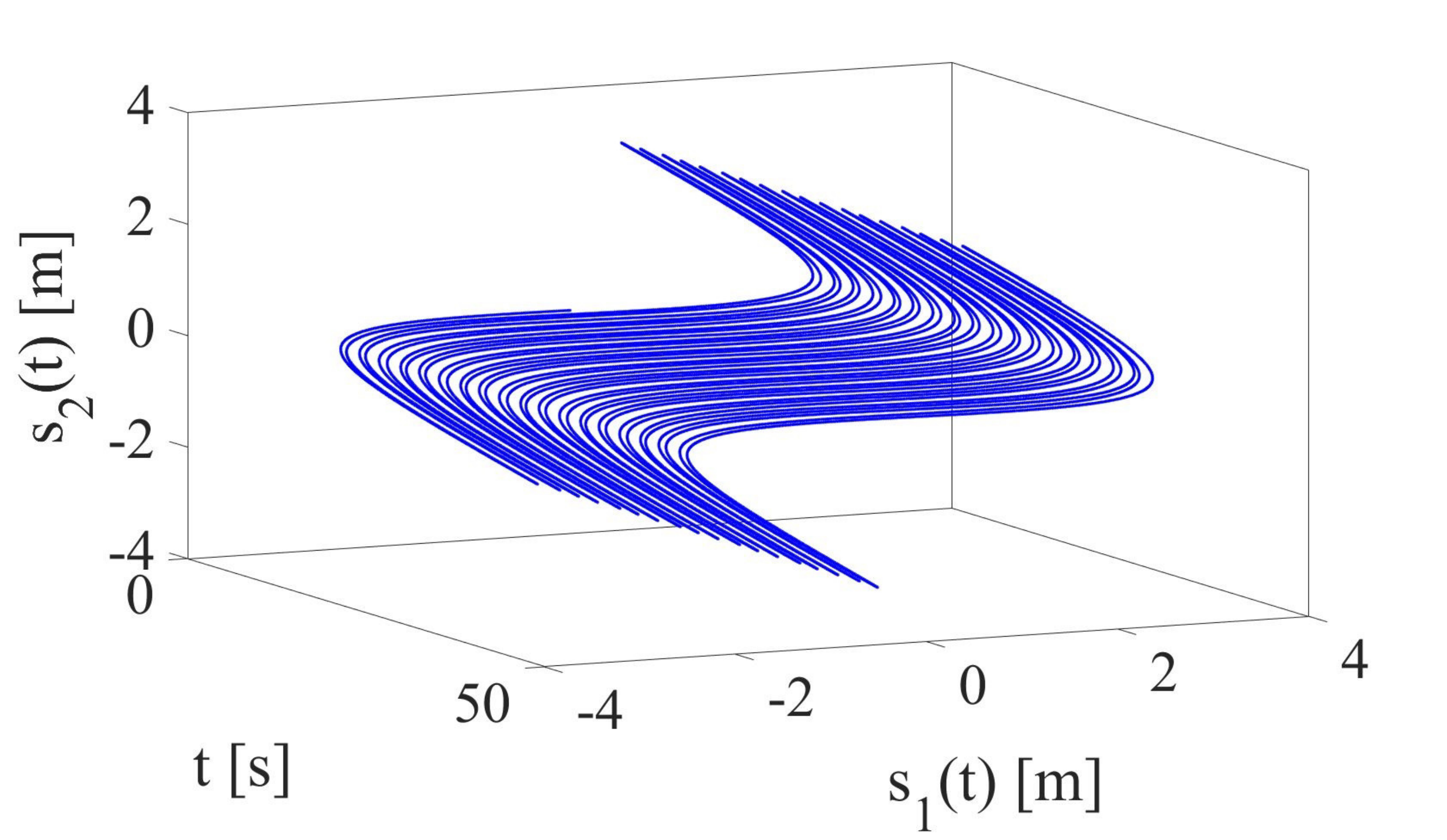}
	\end{tabular}
		\caption{Fully conservative case; extended configuration and velocity diagrams.}
		\label{fig_ex1_3}
\end{figure}

This kind of systems naturally arises in the context of reduced-order models,
see for instance \citep{Kerschen2005, Jansen2005}.
To perform our computations, we adopt a model with two degrees of freedom used as a demonstrator in
\citep{Kerschen2005}. The non-zero mechanical properties are $M_{11}=M_{22}=1\:\mathrm{Kg}$,
$V_{11}^{\mathrm{II}}=V_{22}^{\mathrm{II}}=16\:\mathrm{N/m}$,
$V_{12}^{\mathrm{II}}=V_{21}^{\mathrm{II}}=-15\:\mathrm{N/m}$ and
$V_{1111}^{\mathrm{IV}}=15\:\mathrm{N/m^3}$. The simulation parameters are initial time
$t_{i}=0\:\mathrm{s}$, final time $t_{f}=T\:\mathrm{s}$, simulation time $T=50\:\mathrm{s}$, time step $\Delta
t=0.001\:\mathrm{s}$ and relative iteration tolerance $\varepsilon=10^{-10}$. Additionally, for the
dissipative cases, we set $\chi_{f}=0.0025$ and $\chi_{s}=0.008$ as well as
$\bm{D}=\bm{V}^{\mathrm{II}}$. The initial conditions employed are

%
\begin{equation*}
\bm{q}_{0}=\left[\begin{array}{c} 1.00000 \\ 0.91800\end{array}\right] 
\qquad\textrm{and}\qquad
\bm{s}_{0}=\left[\begin{array}{c} 0.00000 \\ 0.00000\end{array}\right].
\end{equation*}
Fig. \ref{fig_ex1_1} shows the idealized mechanical system
under consideration and Fig. \ref{fig_ex1_2} presents a plot of the potential function, which is clearly convex
within the region where the dynamics of the system takes place. Fig.~\ref{fig_ex1_3}
shows different plots for the solution of the fully conservative case.
We can observe the very complex and nonlinear oscillatory behavior,
which is also in excellent agreement with those results presented in \citep{Kerschen2005}.

\begin{figure}
	\centering{}
	\begin{tabular}{cc}
		\includegraphics[width=0.45\textwidth]{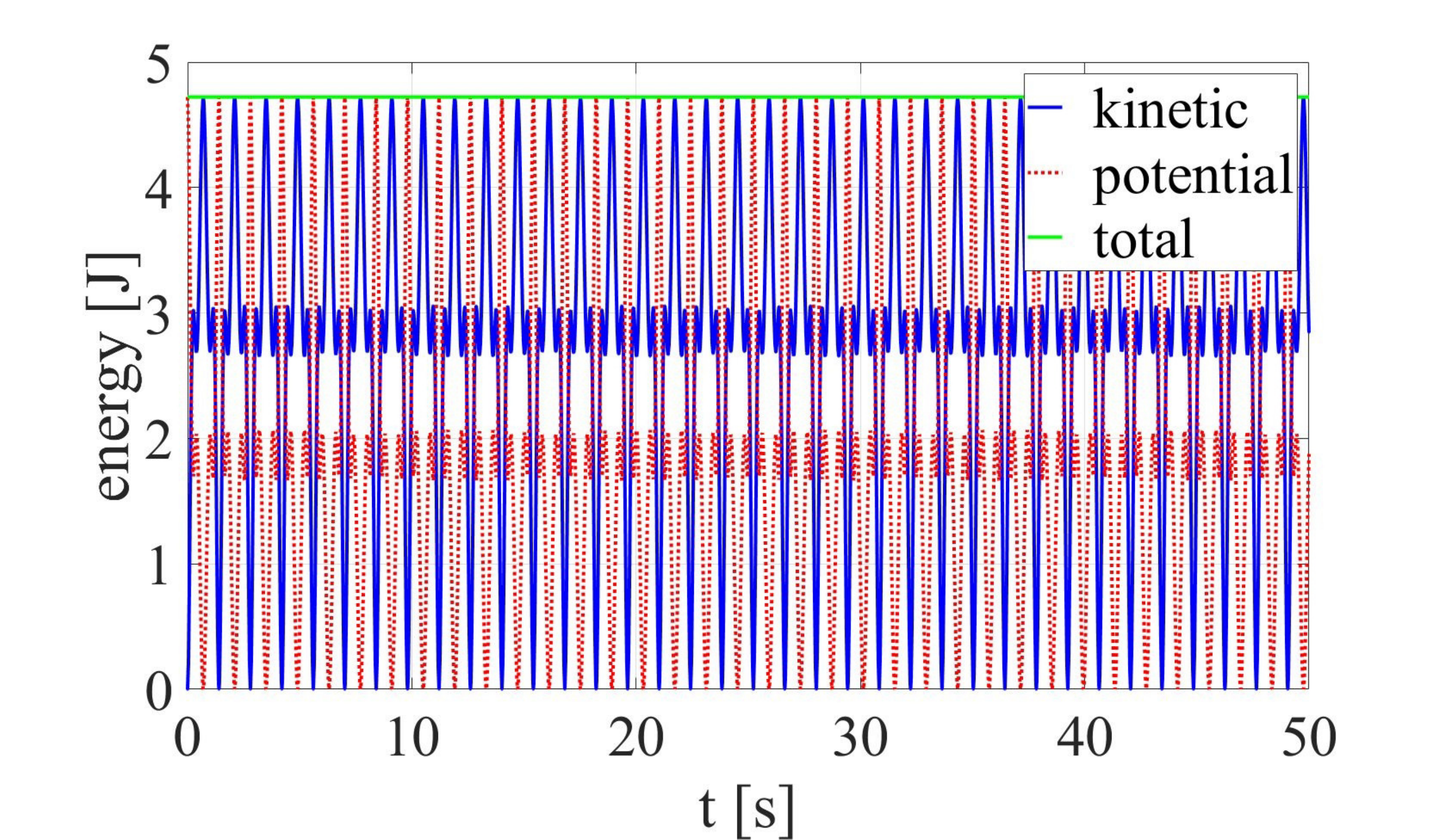} & \includegraphics[width=0.45\textwidth]{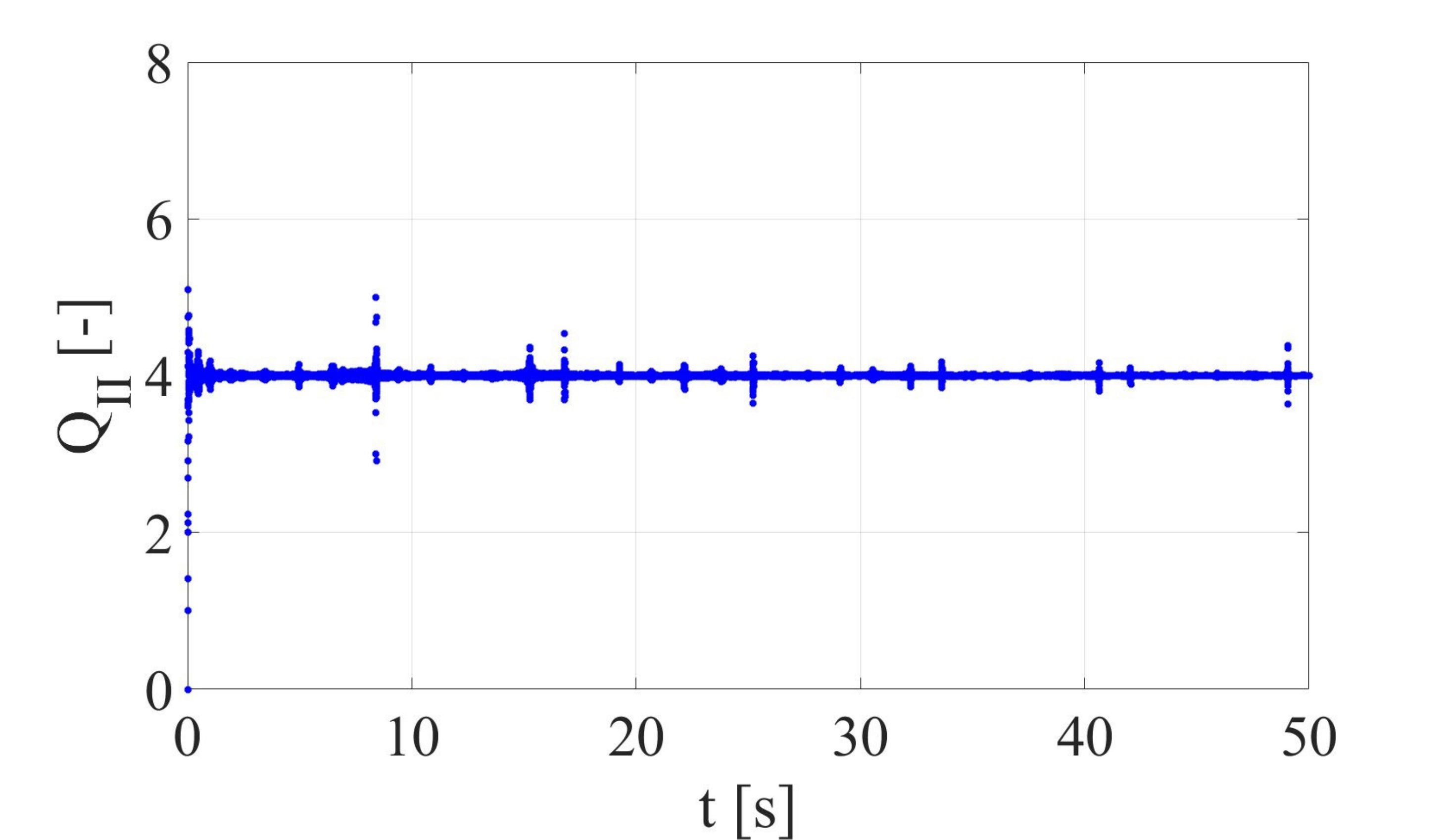}
	\end{tabular}
		\caption{Fully conservative case; energy and $Q_{\mathrm{II}}$.}
		\label{fig_ex1_4}
\end{figure}

\begin{figure}
	\centering{}
	\begin{tabular}{cc}
		\includegraphics[width=0.45\textwidth]{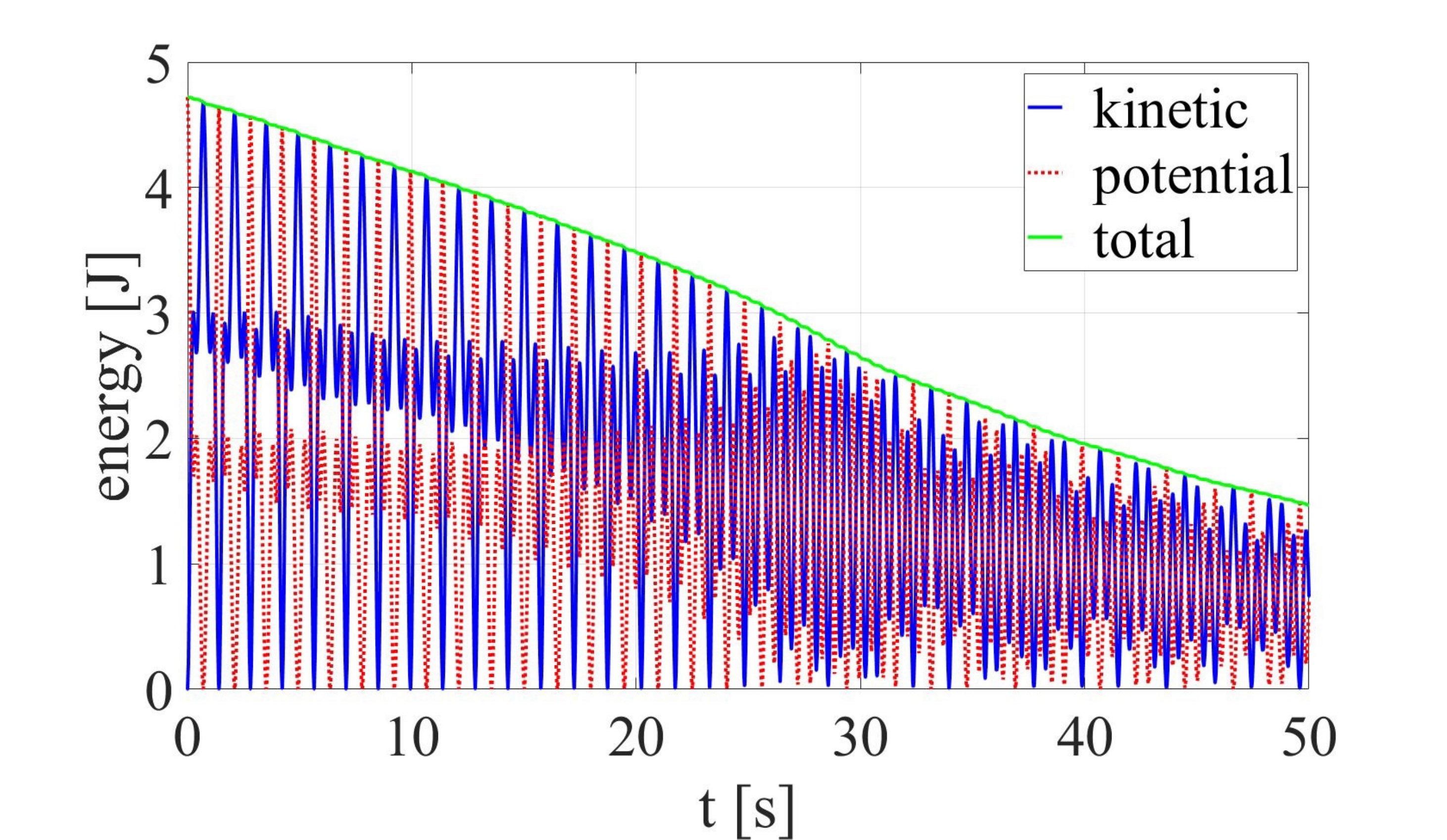} & \includegraphics[width=0.45\textwidth]{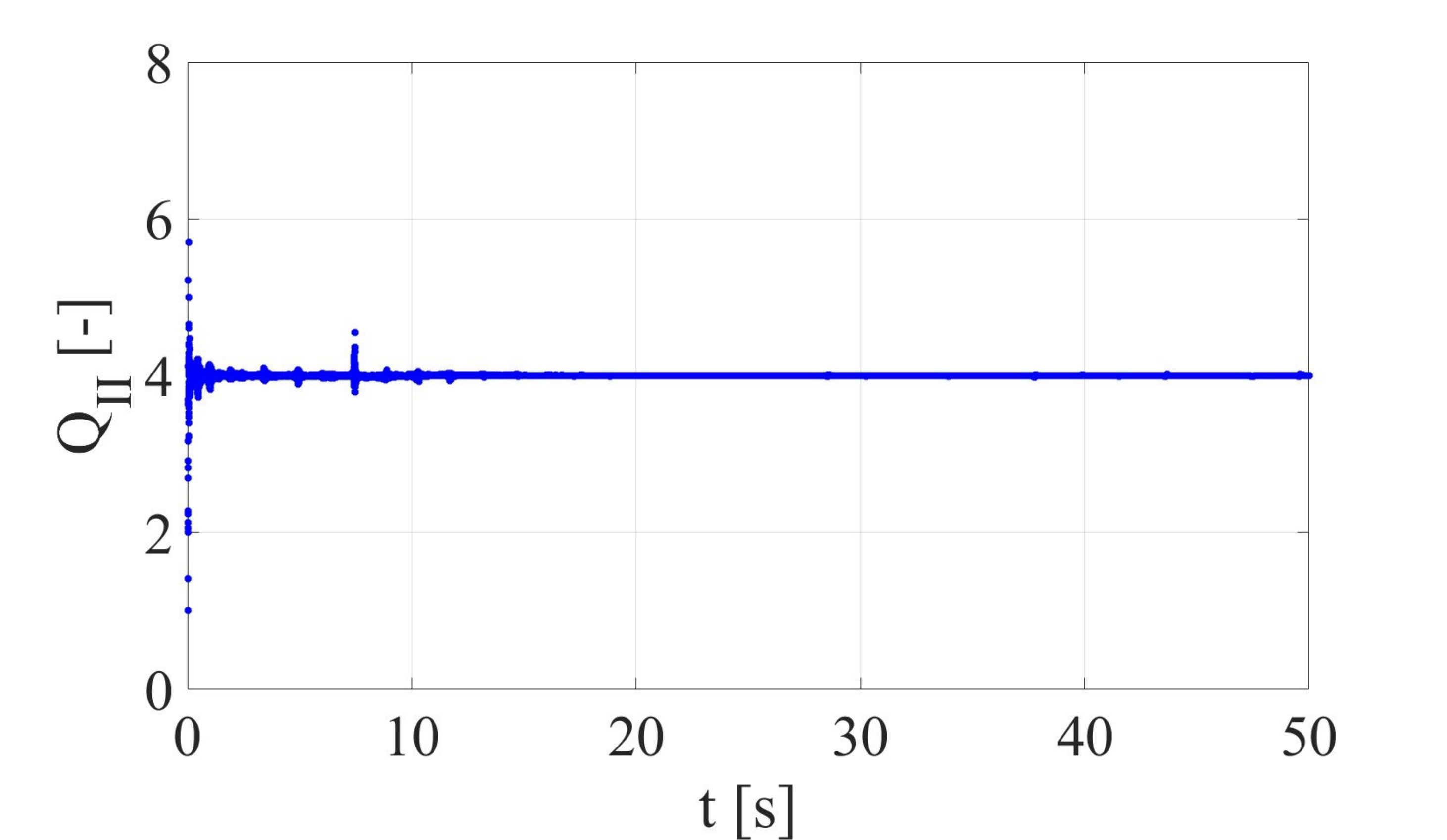}
	\end{tabular}
		\caption{Dissipative case at the level of internal forces; energy and $Q_{\mathrm{II}}$.}
		\label{fig_ex1_5}
\end{figure}

\begin{figure}
	\centering{}
	\begin{tabular}{cc}
		\includegraphics[width=0.45\textwidth]{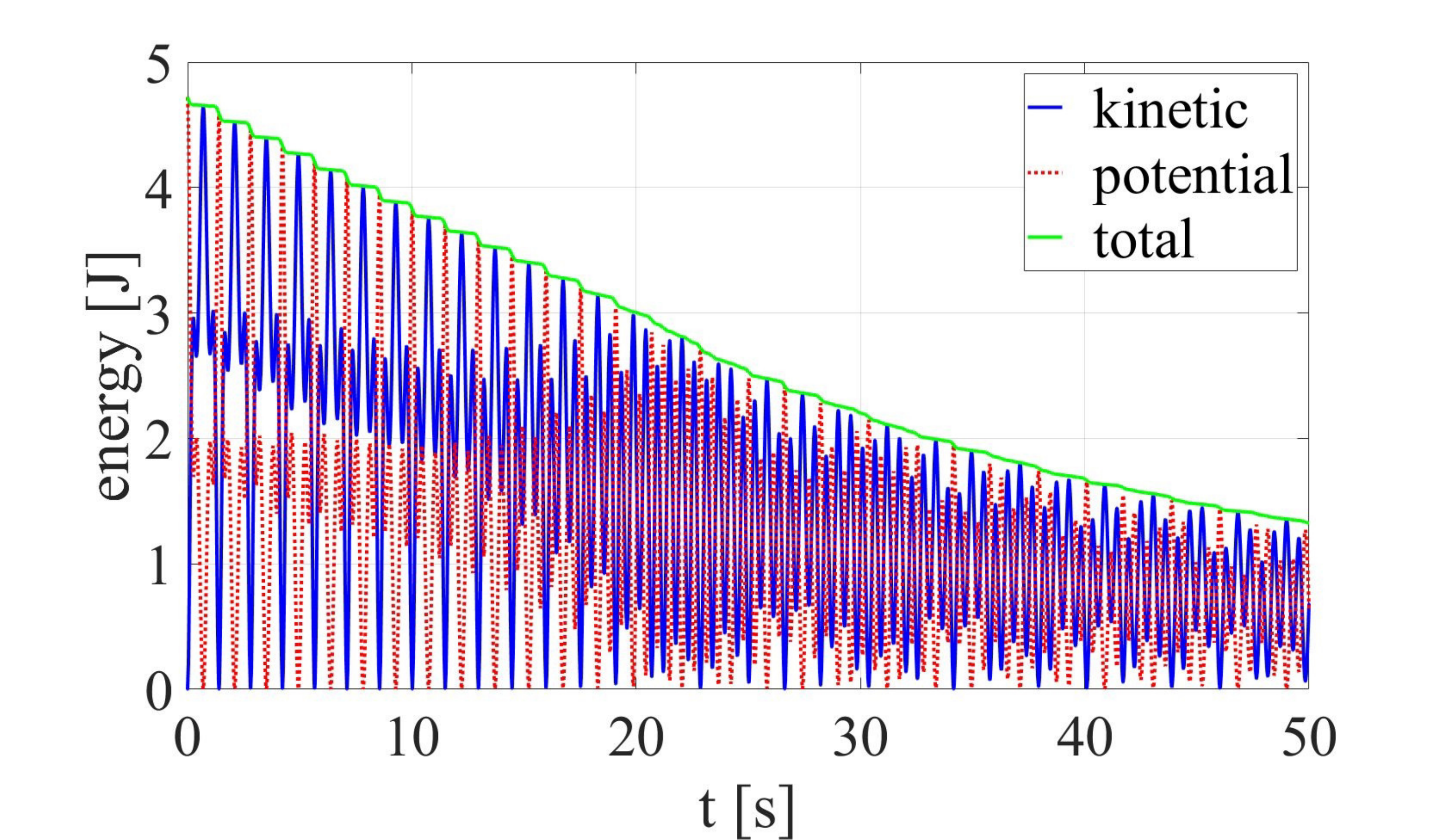} & \includegraphics[width=0.45\textwidth]{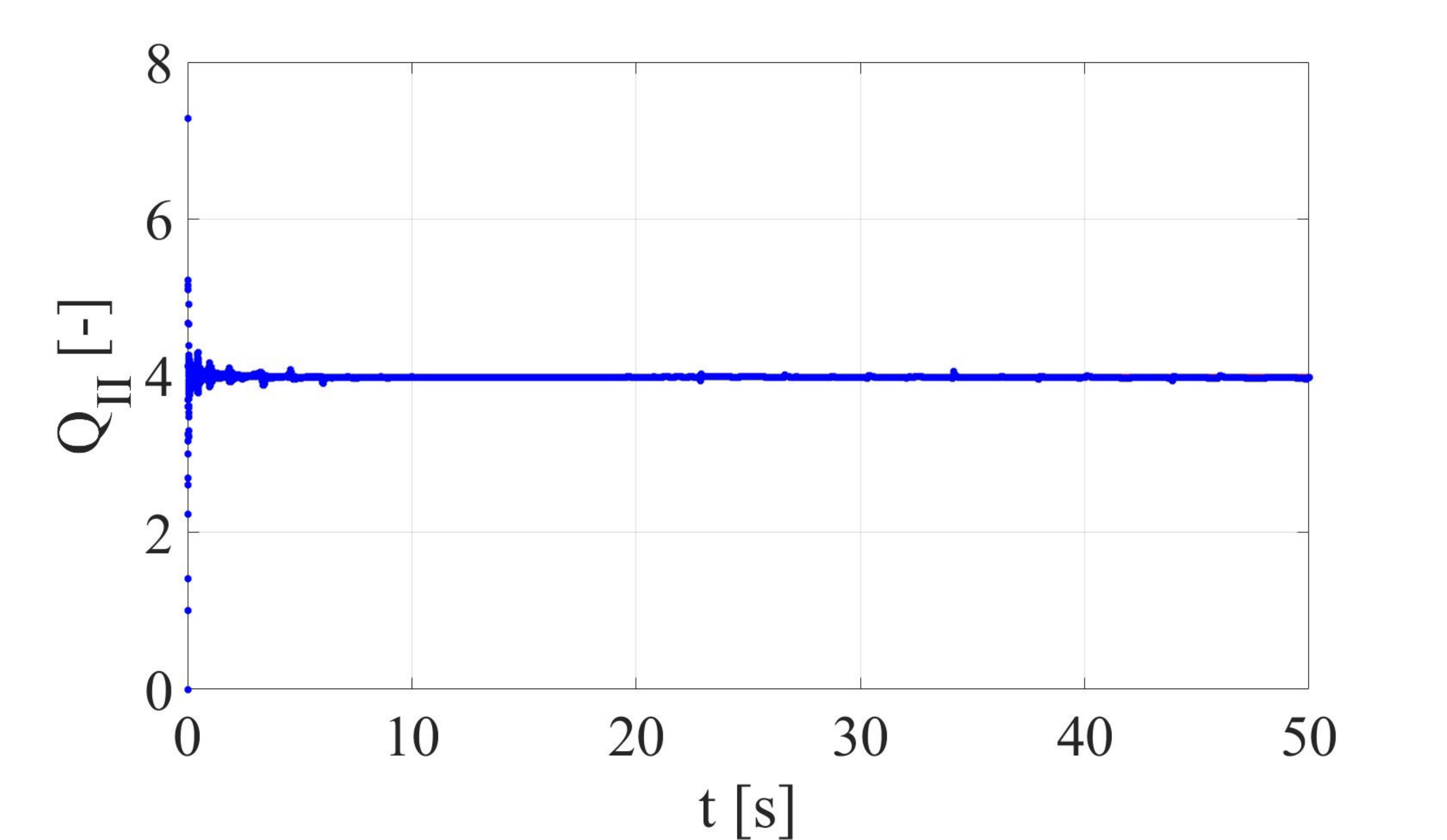}
	\end{tabular}		
		\caption{Dissipative case at the level of generalized velocities; energy and $Q_{\mathrm{II}}$.}
		\label{fig_ex1_6}
\end{figure}

\begin{figure}
	\centering{}
	\begin{tabular}{cc}
		\includegraphics[width=0.45\textwidth]{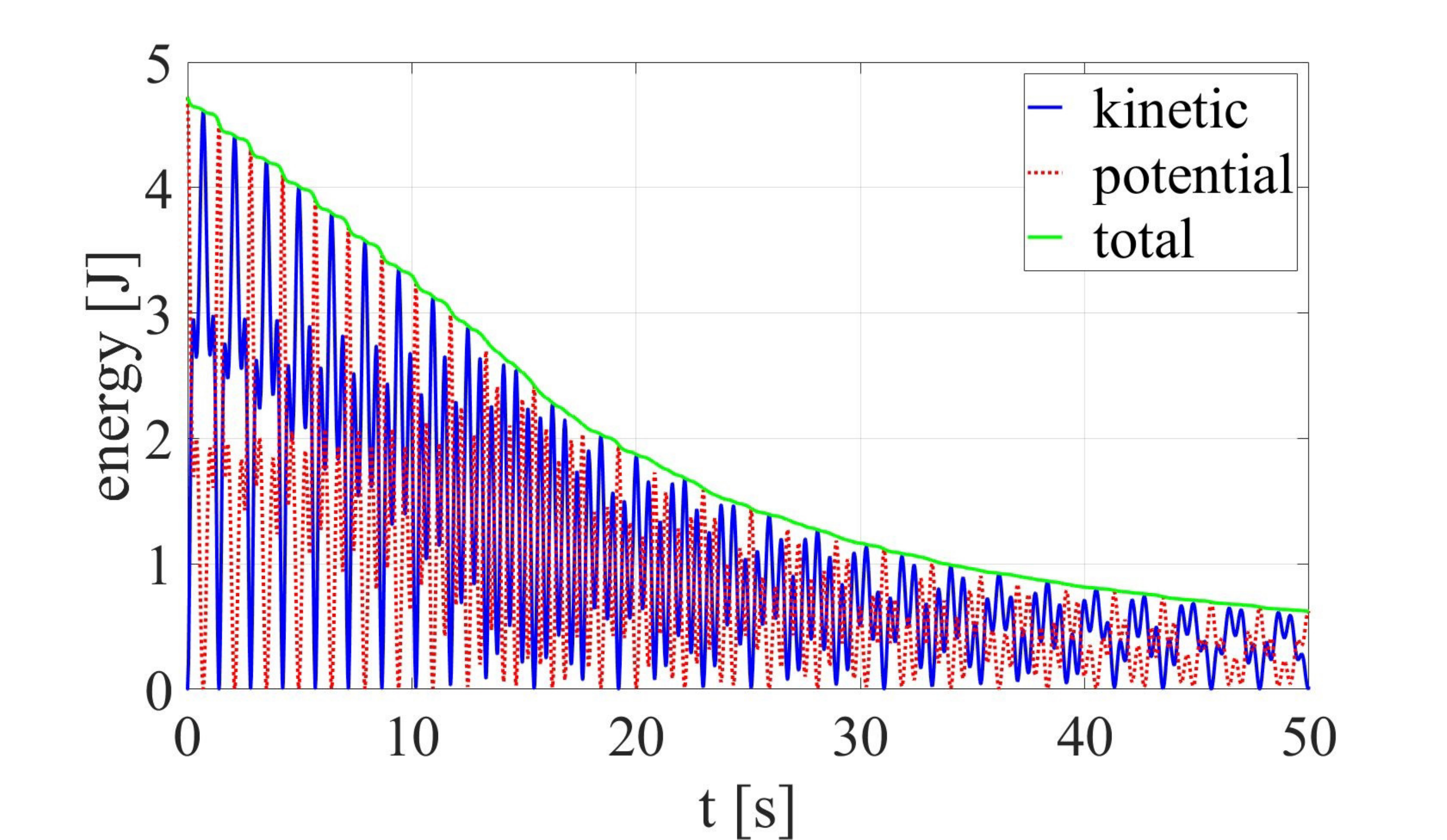} & \includegraphics[width=0.45\textwidth]{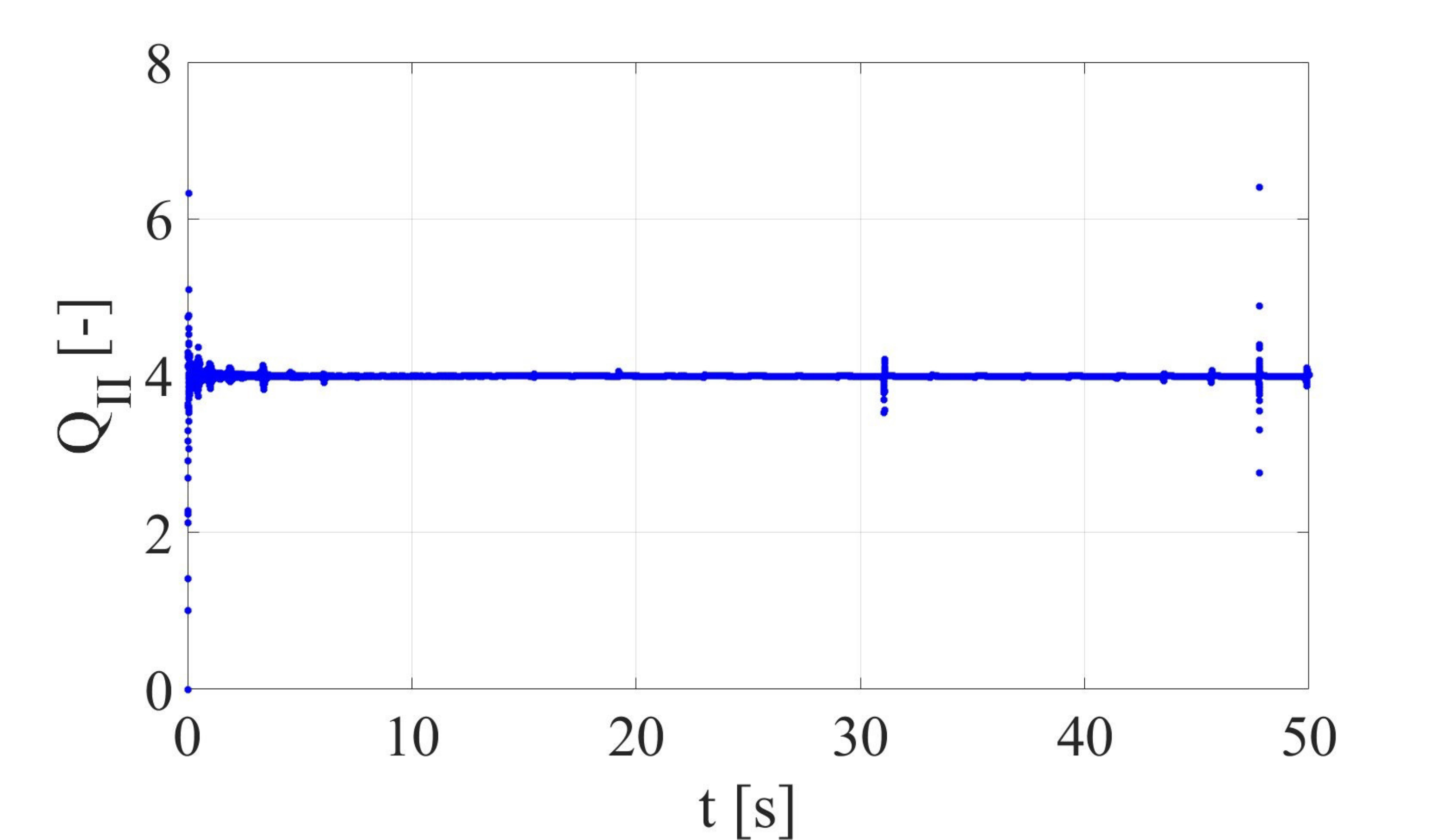}		
	\end{tabular}	
		\caption{Fully dissipative case; energy and $Q_{\mathrm{II}}$.}
		\label{fig_ex1_7}
\end{figure}

On the left of Figs. \ref{fig_ex1_4}, \ref{fig_ex1_5}, \ref{fig_ex1_6} and \ref{fig_ex1_7},
the evolution of the kinetic, potential, and total energies is shown. On the right of
these figure, we show the second precision quotient also as a function of time.
Fig. \ref{fig_ex1_4} evidently corresponds to the fully conservative case.
Fig. \ref{fig_ex1_5} shows the dissipative case at the level of the generalized internal forces.
Fig. \ref{fig_ex1_6} corresponds to the dissipative case at the level of the generalized velocities.
Finally, Fig. \ref{fig_ex1_7} corresponds to the fully dissipative case.
In the latter, the energy decay is larger than the two previous cases.

For all cases the second quotient of precision is almost constant and its value is approximately
$4$. Therefore, as expected, the numerical method is second-order accurate.
According to Eq. \eqref{eq:second_quotient}, a method of a given order is unable to produce
solutions with higher quotients of precision.
In \citep{Kreiss2014}, it is stated that even if
the method is correctly implemented, it is not trivial to find the right set of parameters in order
to numerically obtain precision quotients of a high quality like the one presented herein.

\subsubsection{Two-mass system with a non-polynomial potential}

Next, we consider a nonlinear oscillatory mechanical system whose potential function
is

\begin{equation}
  V(\bm{q})=
  \frac{1}{2}\left[V_{ab}^{\mathrm{II}}
    +\frac{V_{ab}^{\mathrm{N}}}{\left(1+V_{ab}^{\mathrm{D}}q^{a}q^{b}\right)^{n}}\right]q^{a}q^{b},
\end{equation}

\noindent in which the conditions $V_{ab}^{\mathrm{II}}=V_{ba}^{\mathrm{II}}$,
$V_{ab}^{\mathrm{N}}=V_{ba}^{\mathrm{N}}$ and $V_{ab}^{\mathrm{D}}=V_{ba}^{\mathrm{D}}$
are used for the sake of simplicity. This kind of systems could
arise in the context of Euler-Bernoulli beams \citep{Kopmaz2003,Nayfeh2008} or for structures
with softening behavior \citep{Ye2018}.

To perform our computations, we adopt a model
with two degrees of freedom. The non-zero constants of the model are
are $M_{11}=M_{22}=1\:\mathrm{Kg}$, $V_{11}^{\mathrm{II}}=V_{22}^{\mathrm{II}}=10\:\mathrm{N/m}$, $V_{11}^{\mathrm{N}}=V_{22}^{\mathrm{N}}=-V_{12}^{\mathrm{N}}=-V_{21}^{\mathrm{N}}=300\:\mathrm{N/m}$,
$V_{11}^{\mathrm{D}}=V_{22}^{\mathrm{D}}=-V_{12}^{\mathrm{D}}=-V_{21}^{\mathrm{D}}=5\:\mathrm{N/m^2}$
and $n=3$.

\begin{figure}
		\centering{}
		\includegraphics[width=0.5\textwidth]{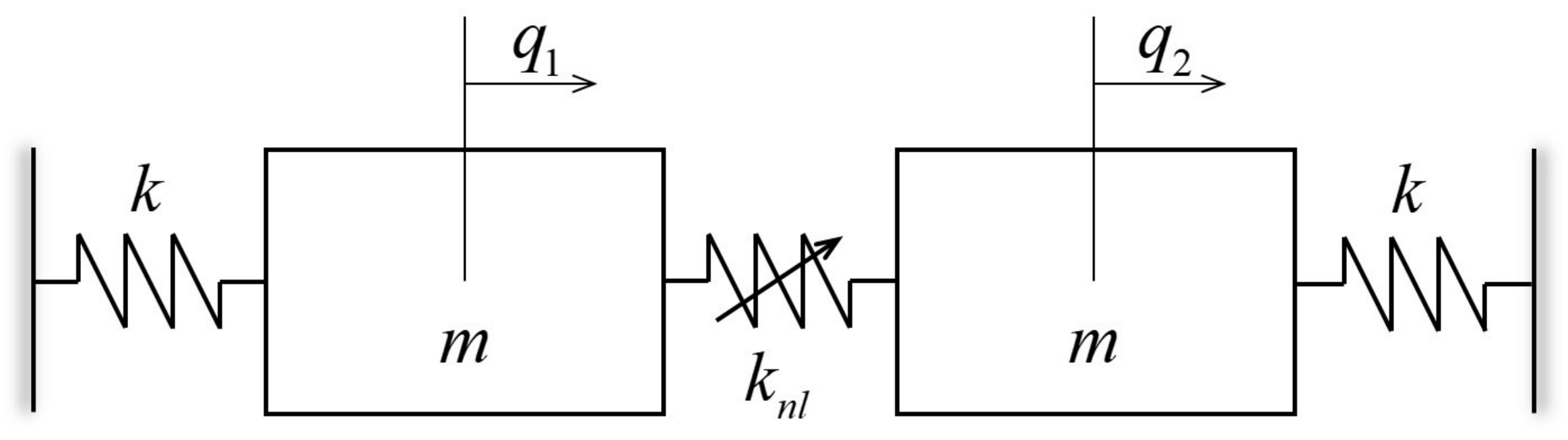}
		\caption{Nonlinear system with two masses connected by linear springs to two walls. Both masses are connected to each other by a nonlinear spring.}
		\label{fig_ex2_1}
\end{figure}

\begin{figure}
		\centering{}
		\includegraphics[width=0.5\textwidth]{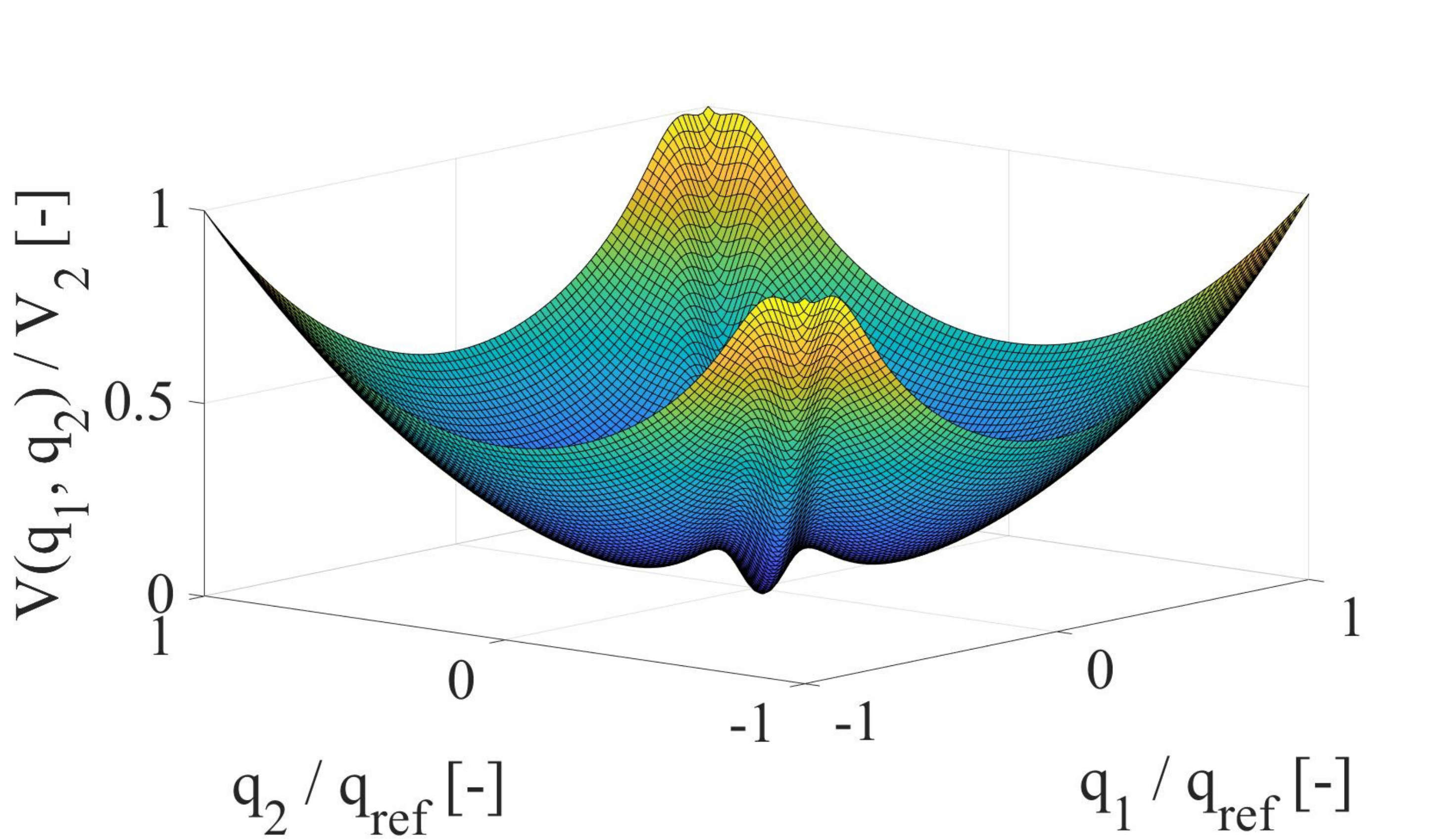}
		\caption{Potential function with non-polynomial complexity.}
		\label{fig_ex2_2}
\end{figure}

\begin{figure}
	\centering{}
	\begin{tabular}{cc}
		\includegraphics[width=0.45\textwidth]{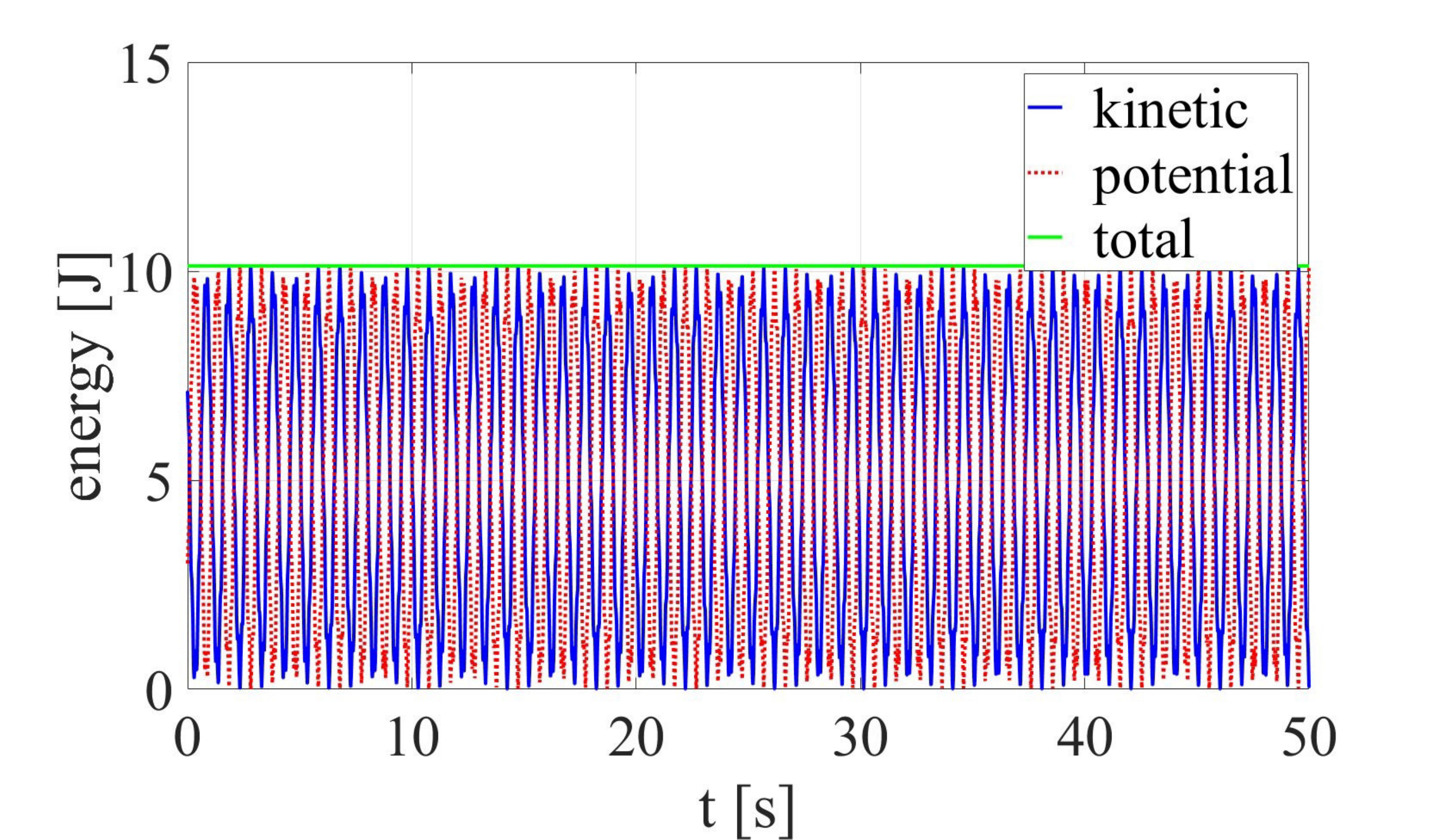} & \includegraphics[width=0.45\textwidth]{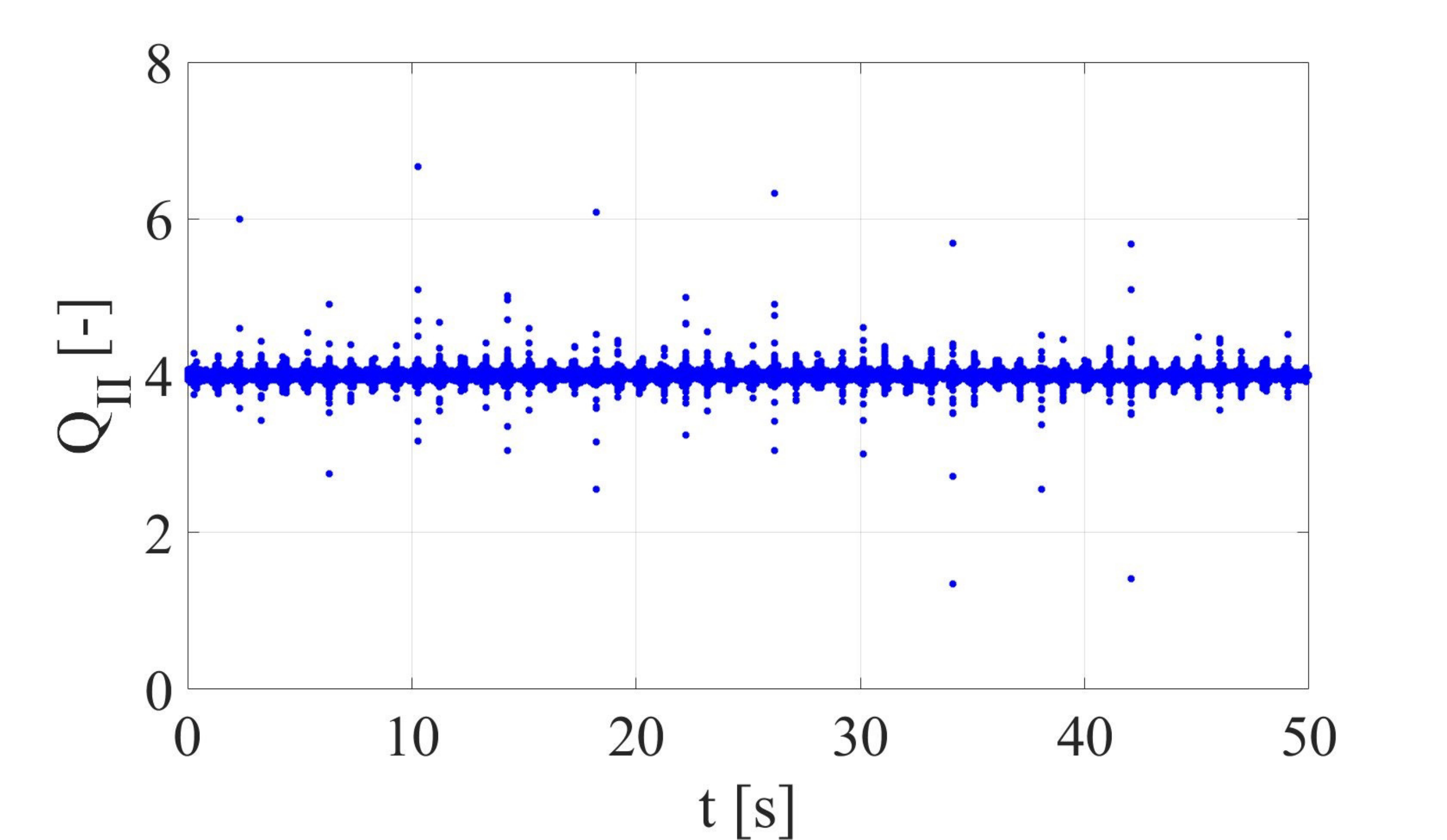}
	\end{tabular}
		\caption{Fully conservative case; energy and $Q_{\mathrm{II}}$.}
		\label{fig_ex2_4}
\end{figure}

\begin{figure}
	\centering{}
	\begin{tabular}{cc}
		\includegraphics[width=0.45\textwidth]{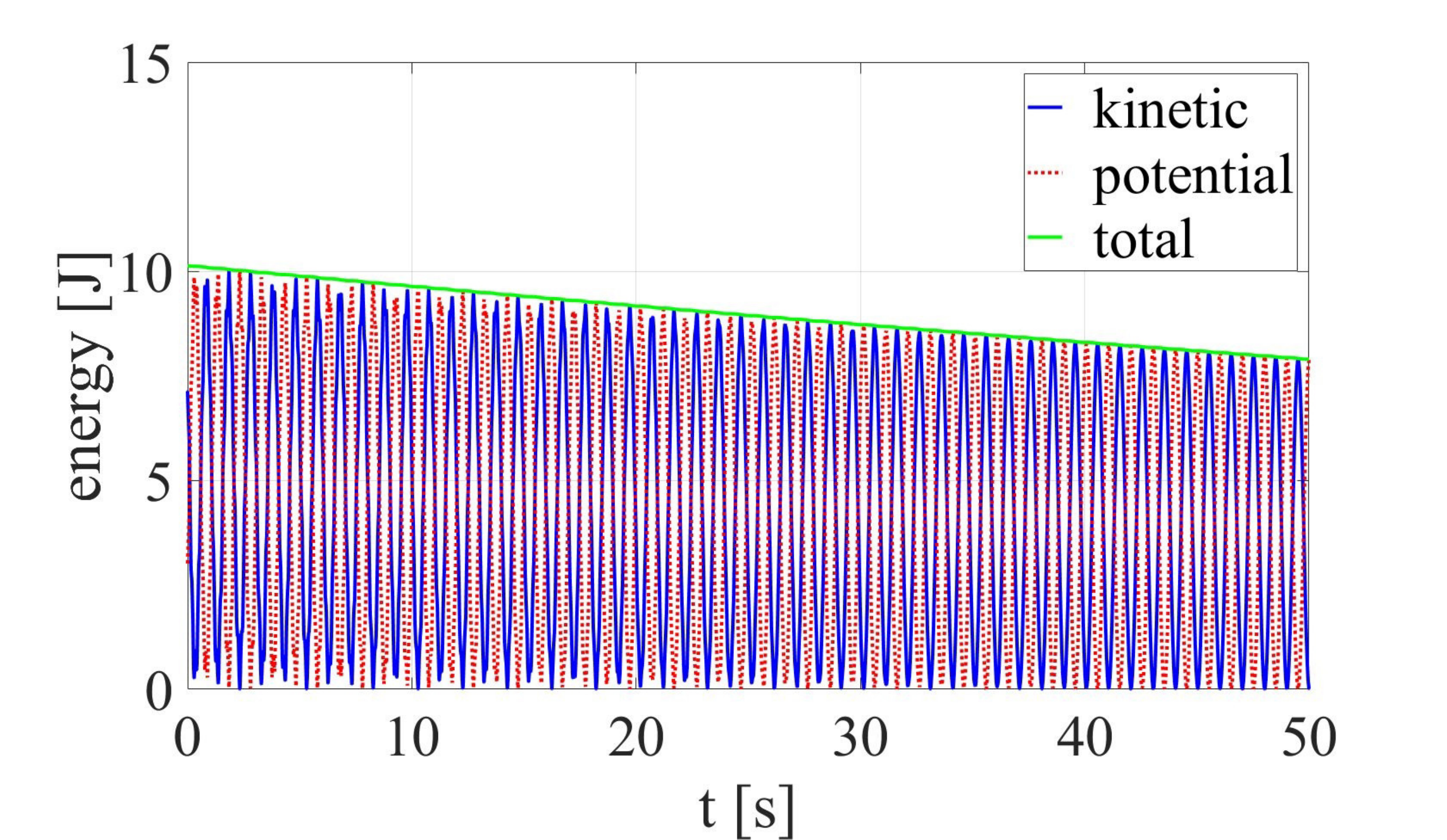} & \includegraphics[width=0.45\textwidth]{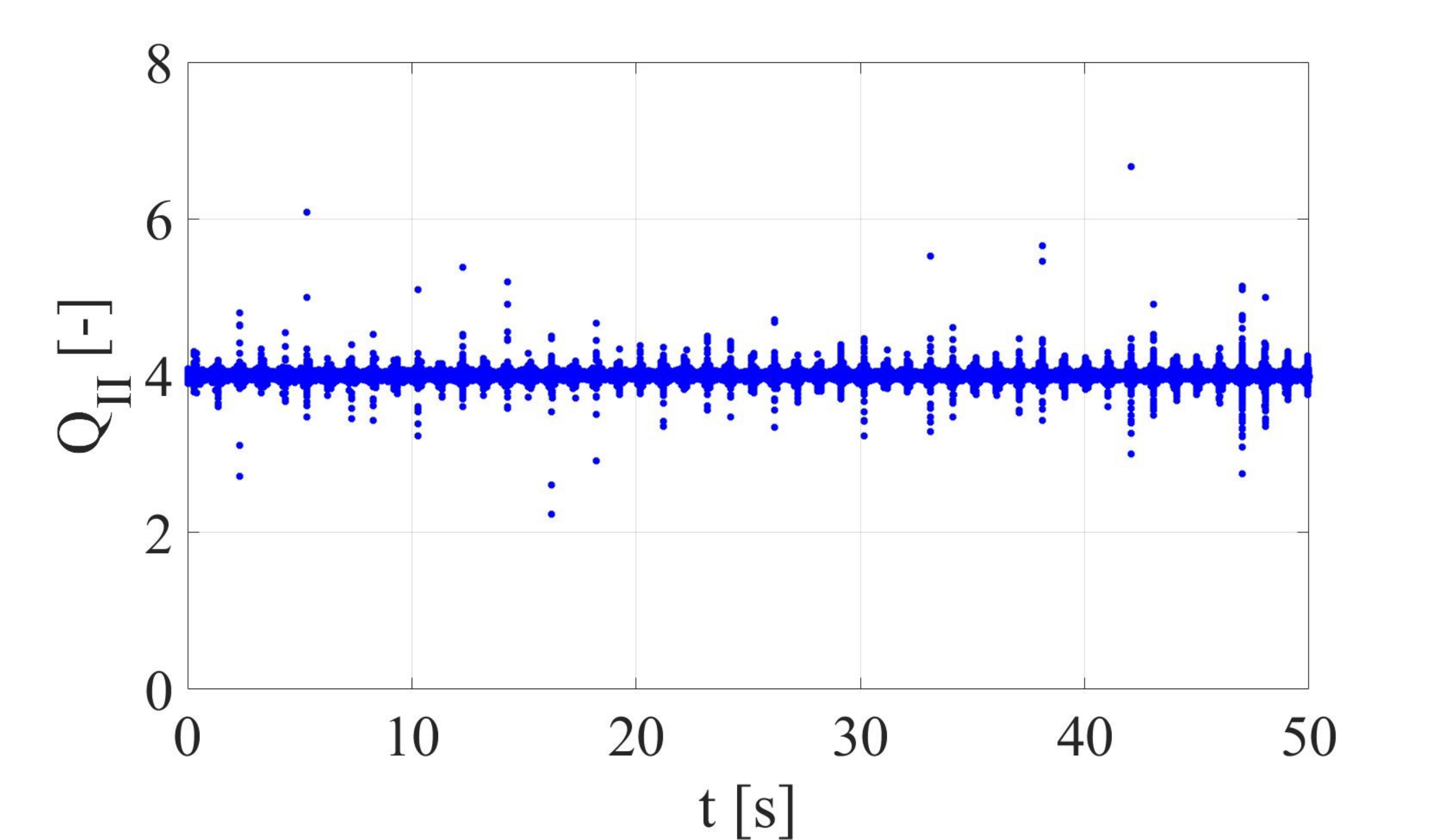}
	\end{tabular}
		\caption{Dissipative case at the level of internal forces; energy and $Q_{\mathrm{II}}$.}
		\label{fig_ex2_5}
\end{figure}

\begin{figure}
	\centering{}
	\begin{tabular}{cc}
		\includegraphics[width=0.45\textwidth]{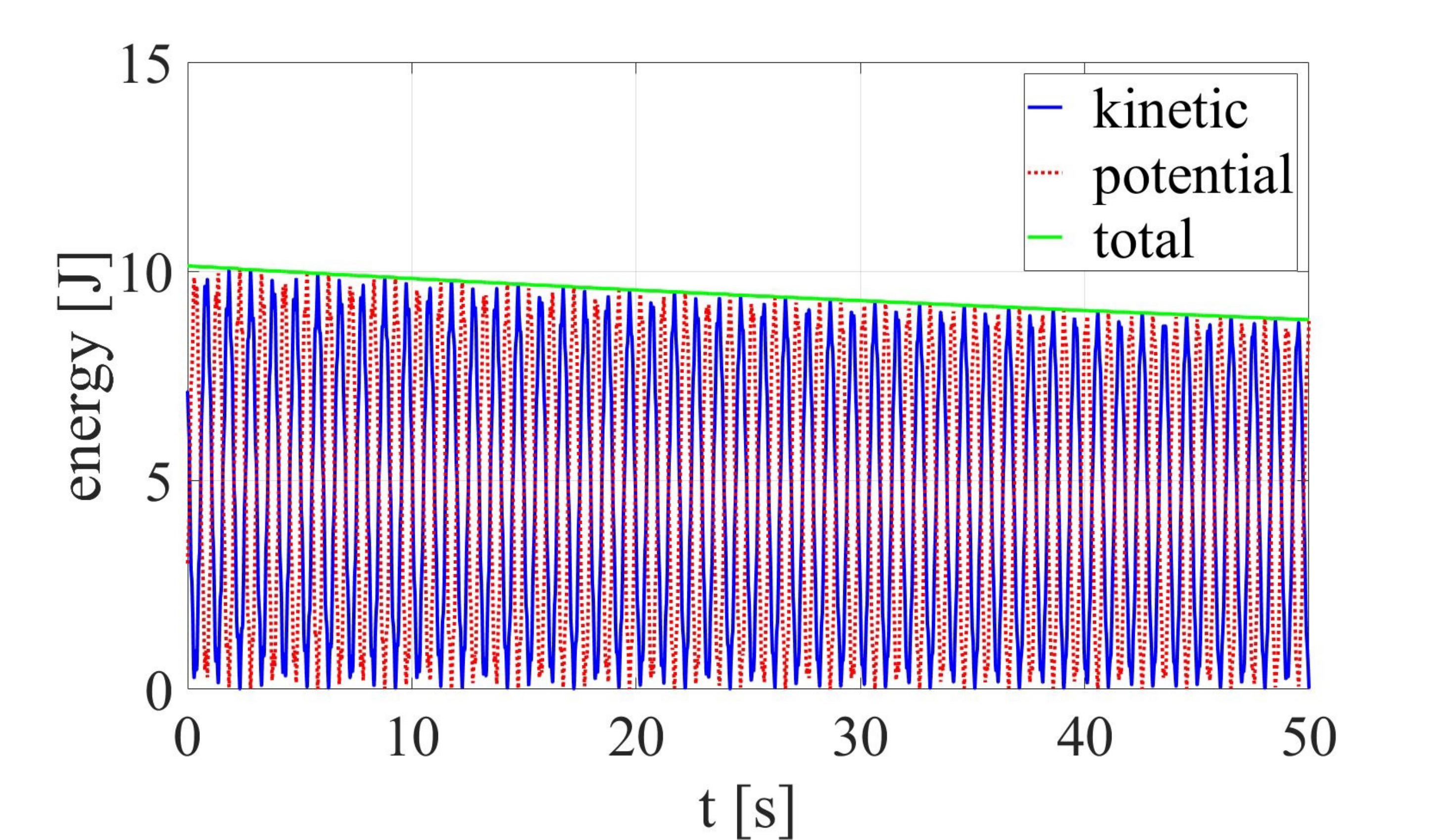} & \includegraphics[width=0.45\textwidth]{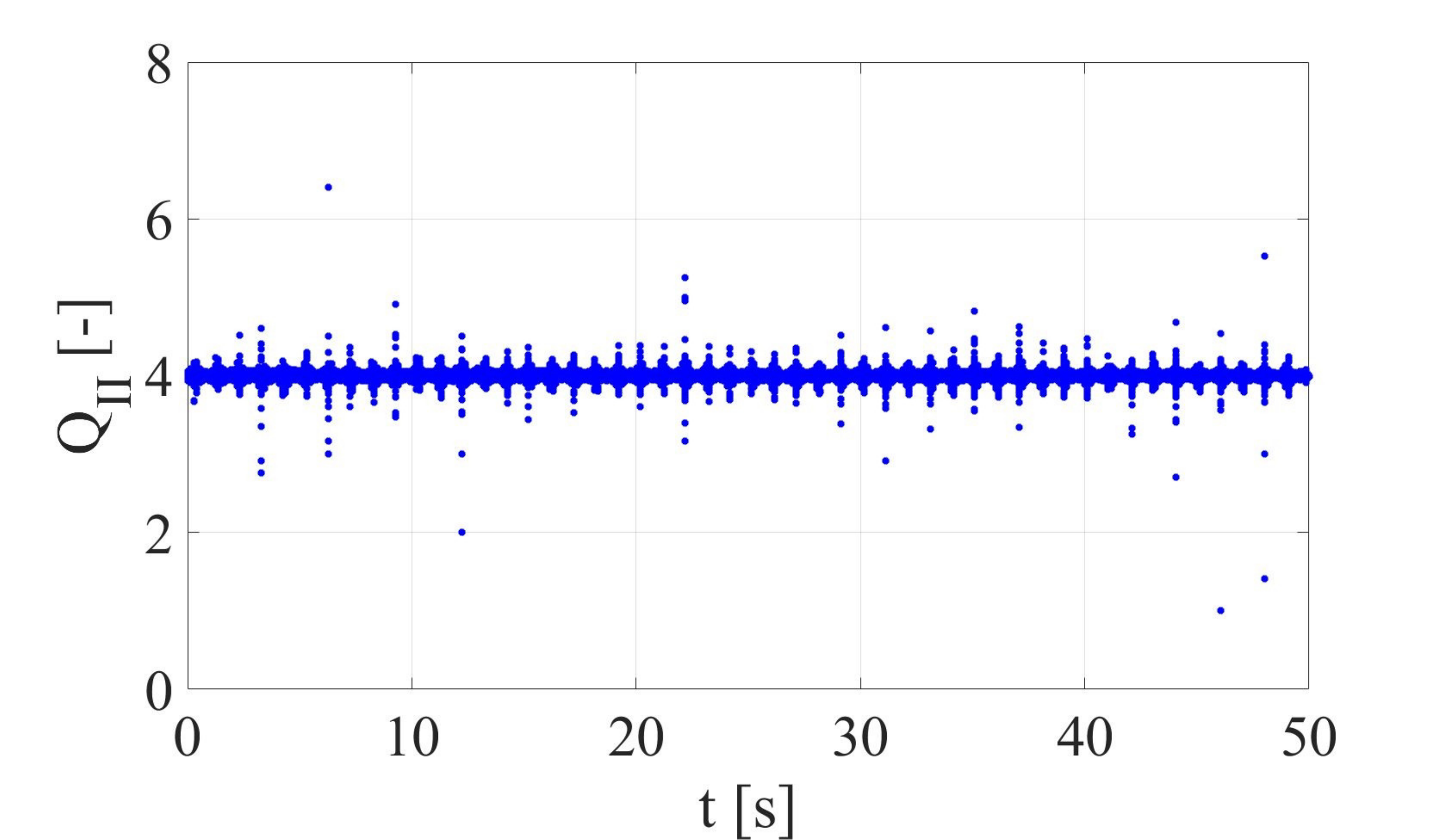}
	\end{tabular}		
		\caption{Dissipative case at the level of generalized velocities; energy and $Q_{\mathrm{II}}$.}
		\label{fig_ex2_6}
\end{figure}

\begin{figure}
	\centering{}
	\begin{tabular}{cc}
		\includegraphics[width=0.45\textwidth]{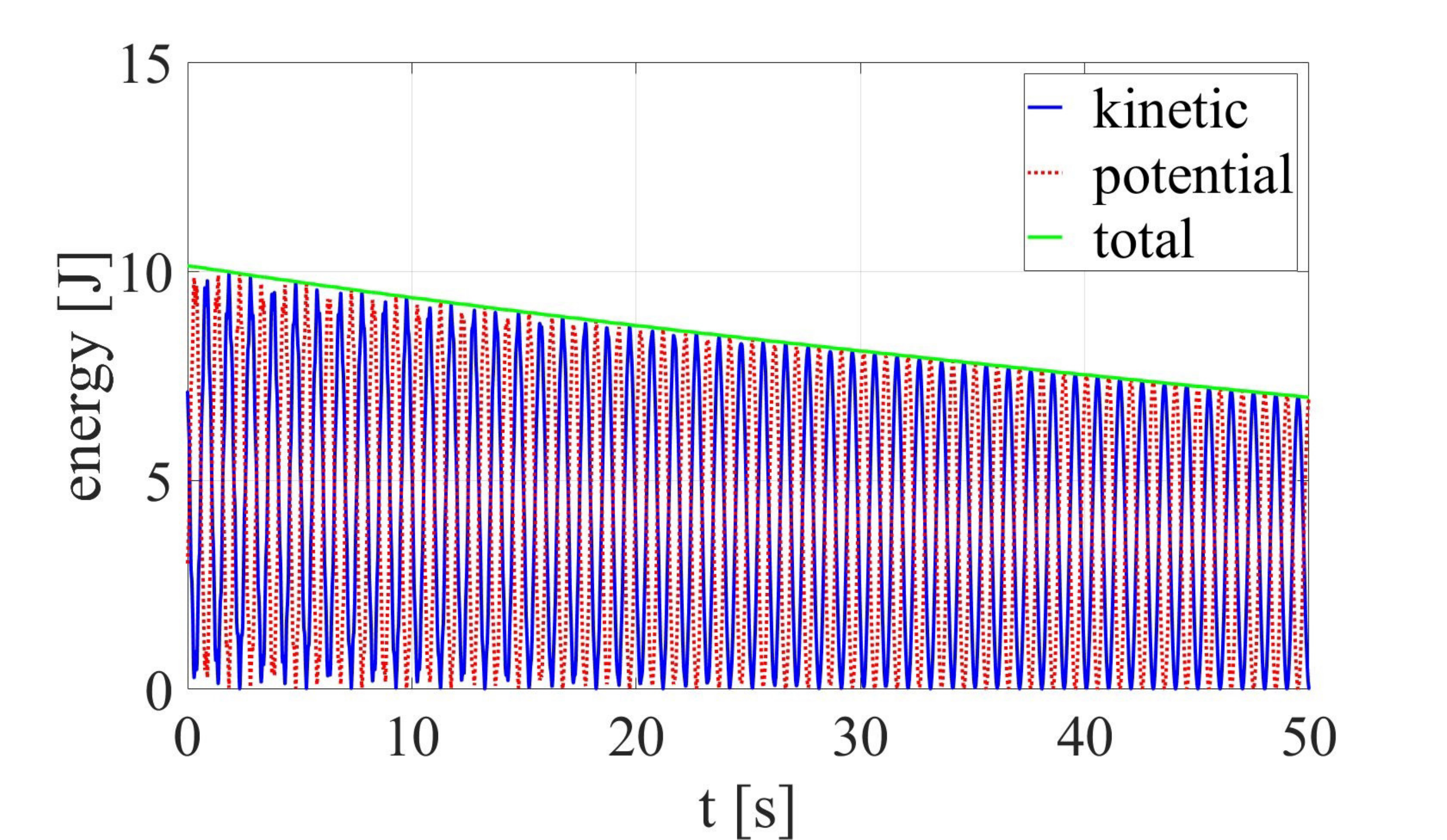} & \includegraphics[width=0.45\textwidth]{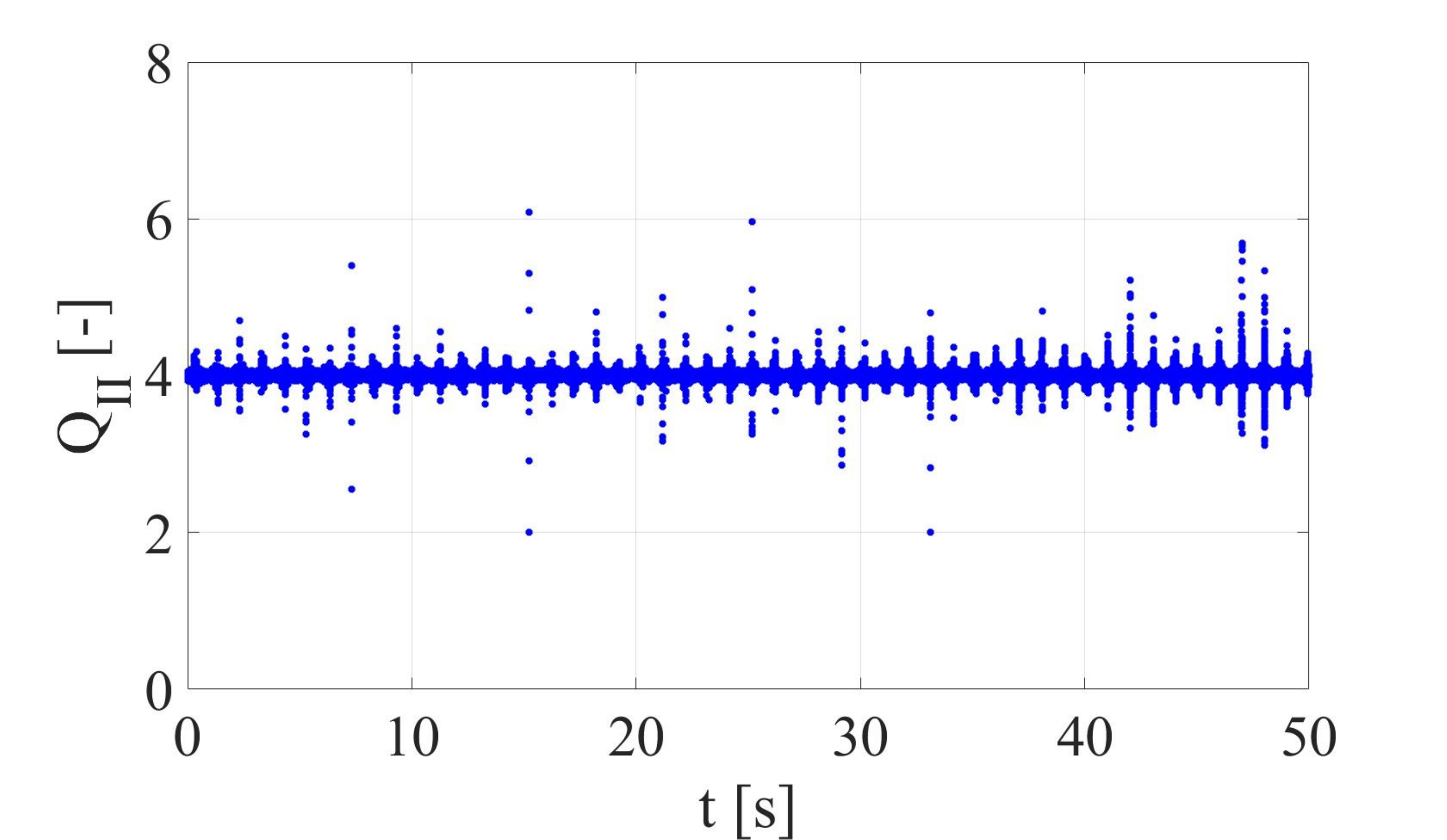}		
	\end{tabular}	
		\caption{Fully dissipative case; energy and $Q_{\mathrm{II}}$.}
		\label{fig_ex2_7}
\end{figure}

The simulation parameters are initial time $t_{i}=0\:\mathrm{s}$,
final time $t_{f}=T\:\mathrm{s}$, simulation time $T=50\:\mathrm{s}$,
time step $\Delta t=0.0001\:\mathrm{s}$, and relative iteration tolerance $\varepsilon=10^{-10}$.
Additionally, for the dissipative cases, we set $\chi_{f}=0.001$ and $\chi_{s}=0.001$ as well as $\bm{D}=\bm{V}^{\mathrm{II}}$. The initial conditions employed are

%
\begin{equation*}
\bm{q}_{0}=\left[\begin{array}{c} -0.41726 \\ -0.49840\end{array}\right] 
\qquad\textrm{and}\qquad
\bm{s}_{0}=\left[\begin{array}{c} -2.53182 \\ -2.79761\end{array}\right].
\end{equation*}
Fig. \ref{fig_ex2_1} shows the idealized mechanical system under consideration and
Fig. \ref{fig_ex2_2} depicts the potential function, which is clearly non-convex within
the region where the dynamics of the system takes place. This feature pushes the numerical method to
its limits. On the left of Figs. \ref{fig_ex2_4}, \ref{fig_ex2_5}, \ref{fig_ex2_6} and \ref{fig_ex2_7}
the kinetic, potential and total energies and plotted. On their right, these show
the second precision quotient also as a function of time.  Fig. \ref{fig_ex2_4}
depicts the energies in the conserving solution.  Figs. \ref{fig_ex2_5} and \ref{fig_ex2_6}
plot the energies when dissipation is introduced in the internal forces and generalized
velocities, respectively. Last, Fig. \ref{fig_ex2_7}, provides the results obtained
when both dissipation functions are employed, resulting in a larger dissipation of energy.
In all cases, the second quotient of precision is very close to 4 for all time, confirming
the second order accuracy of the method in all the simulations.

\begin{figure}
	\centering{}
	\begin{tabular}{cc}
		\includegraphics[width=0.45\textwidth]{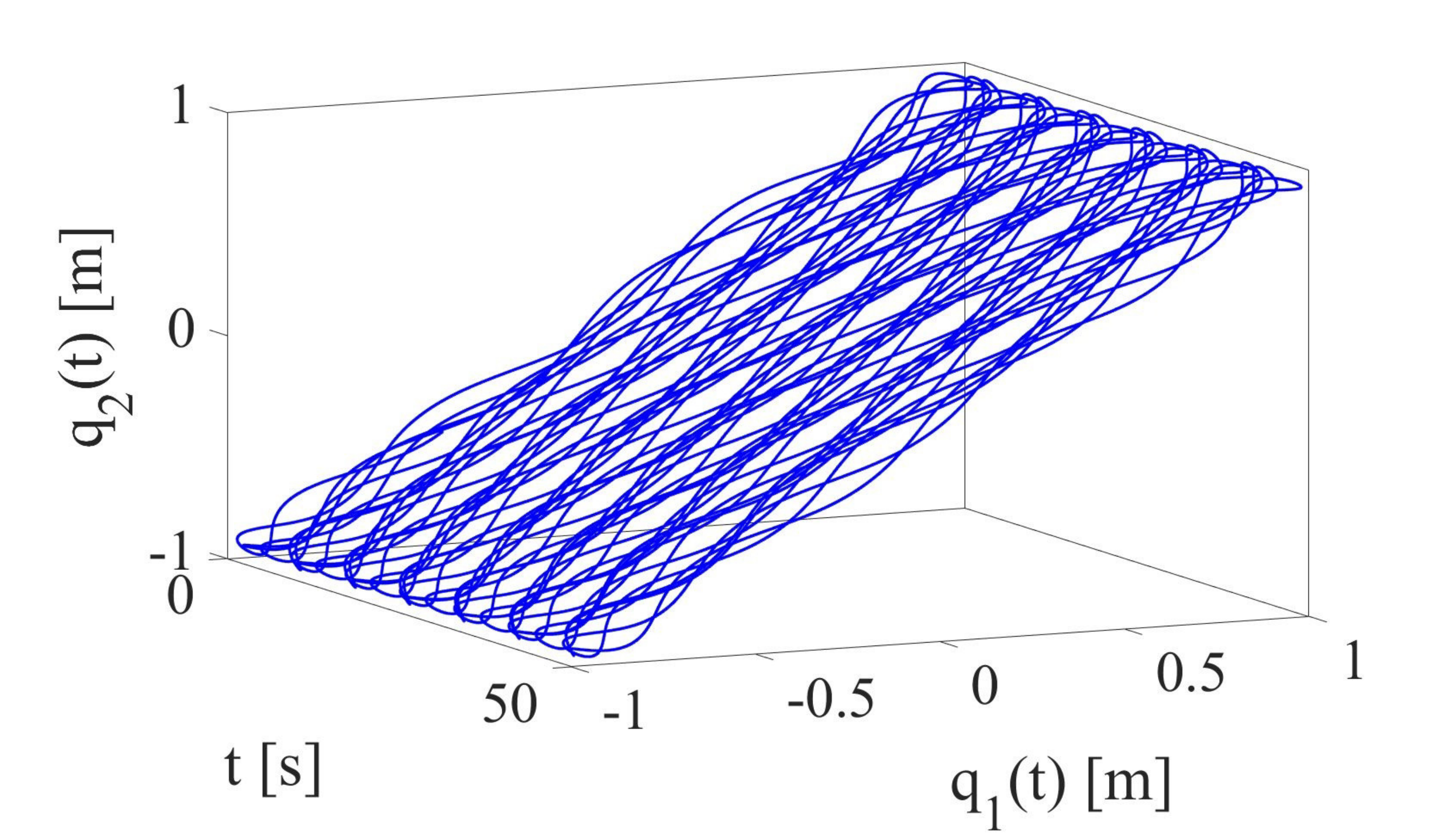} &	\includegraphics[width=0.45\textwidth]{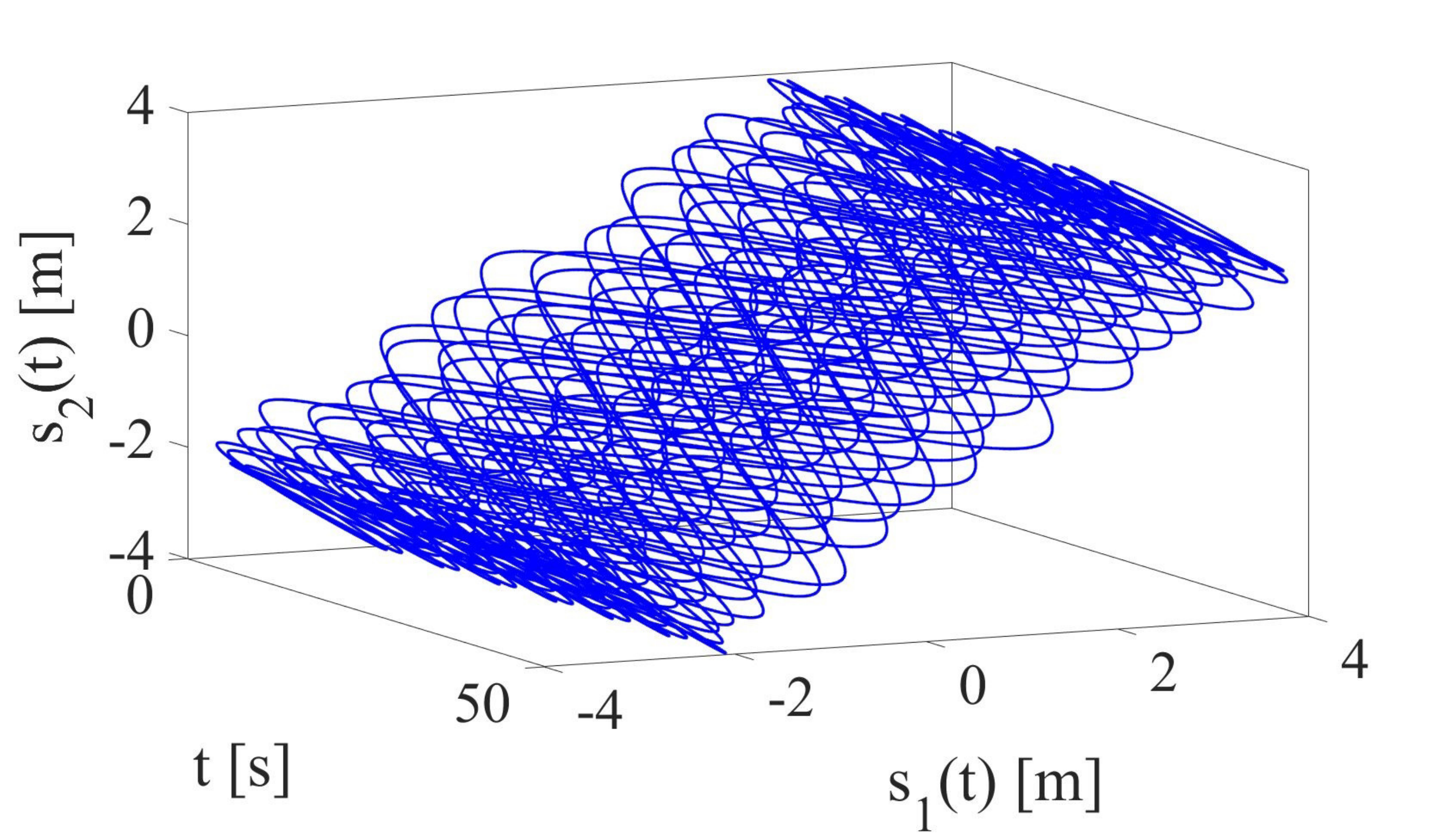}
	\end{tabular}
		\caption{Fully conservative case; extended configuration and velocity diagrams.}
		\label{fig_ex2_3}
\end{figure} 

\subsection{Finite elasticity models}

Here we analyze two finite elasticity models. The first one is a tumbling cylinder and the second
one is a free-flying, single-layer, shell. In both cases, the spatial discretizations
are based on a four-node shell
element \cite{Gebhardt2017, Gebhardt2019a}, \ie~ an extensible-director-based solid-degenerate shell
model, in which the shear locking and the artificial thickness strains are controlled by means of the
assumed natural strain method. Also, the enhancement of the strain field in the thickness direction
and the cure of the membrane locking are achieved by means of the enhanced assumed strain
method. Such element allows to consider unmodified three-dimensional constitutives laws. For the
current study, we adopt the neo-Hookean hyperelastic material model, whose strain energy density is
given by
\begin{equation}
  \tilde{W}(\bm{C}) = \frac{\lambda}{2}\log^2(J)+\frac{\mu}{2}(I_1-3)-\mu\log(J)\, ,
\end{equation}
with $\bm{C}$, the right Cauchy-Green deformation tensor, $J = \sqrt{\det(\bm{C})}$, $I_1 =
\mathrm{trace}(\bm{C})$, and $\lambda$ and $\mu$ are the first and second Lam\'e parameters,
respectively.

\subsubsection{Tumbling cylinder}
\begin{figure}[!ht]
	\centering{}
	\begin{tabular}{ccc}
		$t = 0.00$ s & $t = 1.36$ s & $t = 2.72$ s \tabularnewline
		\tabularnewline
		\includegraphics[align=c, scale=0.5]{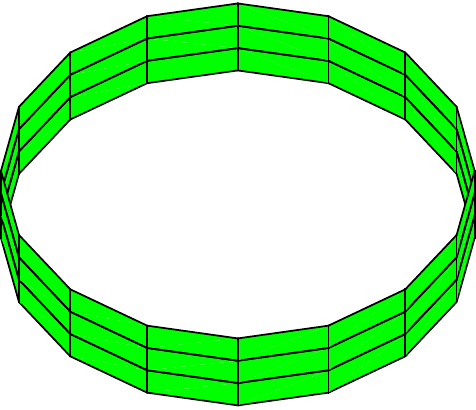} &
		\includegraphics[align=c, scale=0.5]{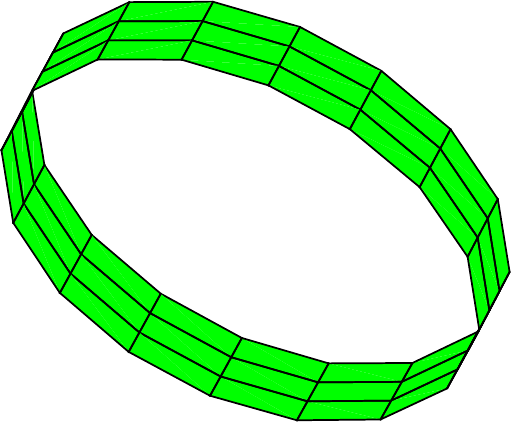} &
		\includegraphics[align=c, scale=0.5]{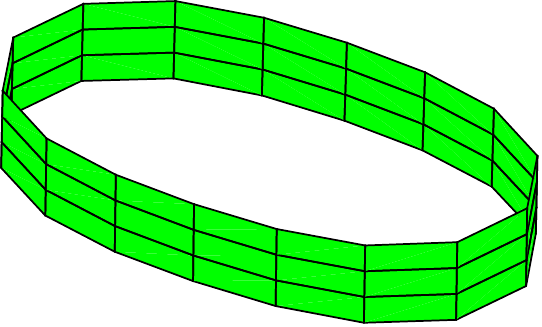} \tabularnewline
		\tabularnewline
		$t = 4.08$ s & $t = 5.44$ s & $t = 6.80$ s \tabularnewline
		\tabularnewline
		\includegraphics[align=c, scale=0.5]{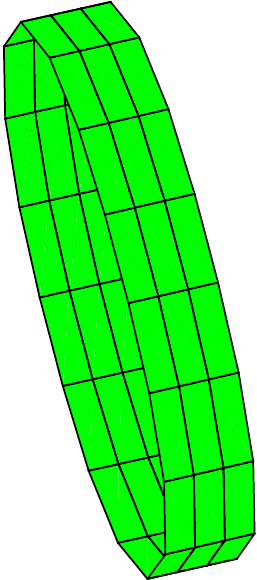} &
		\includegraphics[align=c, scale=0.5]{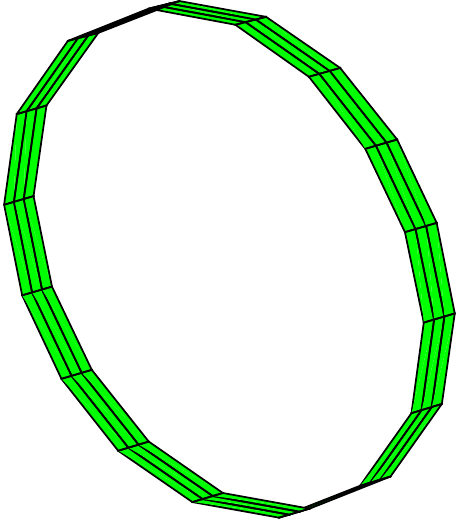} &
		\includegraphics[align=c, scale=0.5]{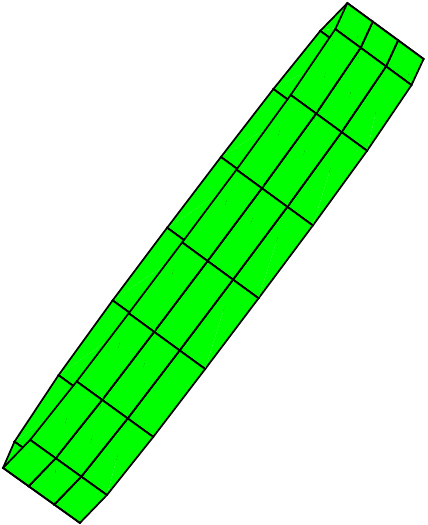} \tabularnewline
		\tabularnewline
		$t = 8.16$ s & $t = 9.52$ s & $t = 10.88$ s \tabularnewline
		\tabularnewline
		\includegraphics[align=c, scale=0.5]{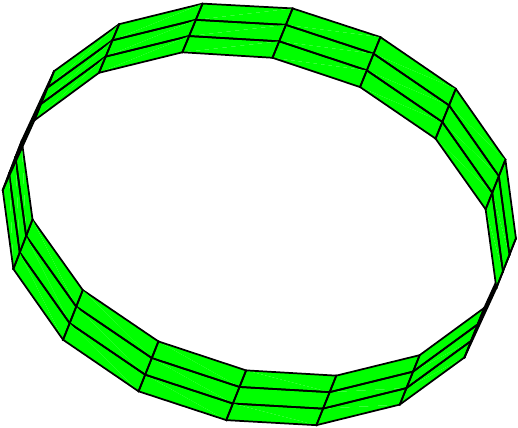} &
		\includegraphics[align=c, scale=0.5]{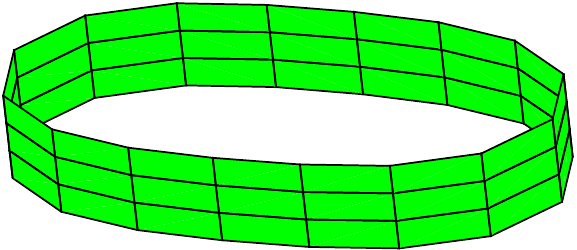} &
		\includegraphics[align=c, scale=0.5]{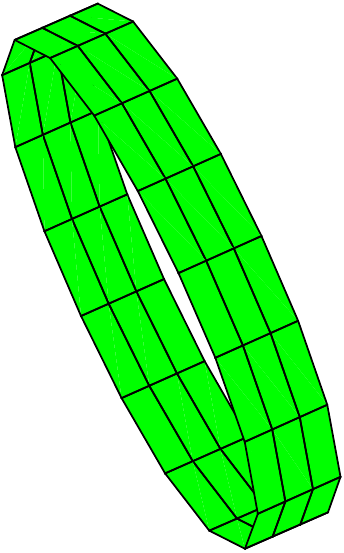} \tabularnewline
		\tabularnewline
		$t = 12.24$ s & $t = 13.60$ s & $t = 14.96$ s \tabularnewline
		\tabularnewline
		\includegraphics[align=c, scale=0.5]{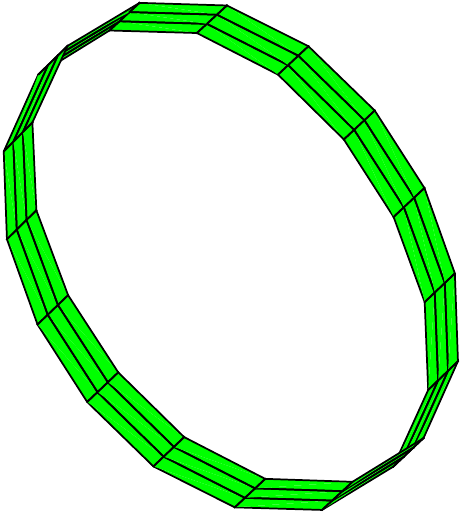} &
		\includegraphics[align=c, scale=0.5]{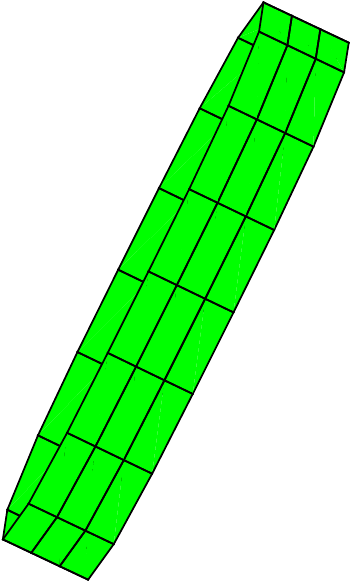} &
		\includegraphics[align=c, scale=0.5]{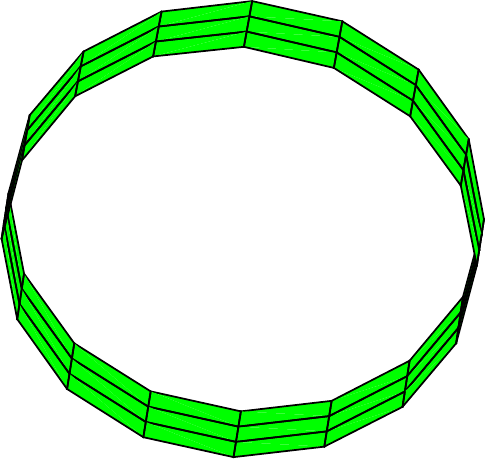}
	\end{tabular}
	\caption{Tumbling cylinder (conservative) - sequence of motion.}
	\label{fig:tumbling cylinder conservative - sequence of motion}	
\end{figure}

\begin{figure}[!ht]
	\centering{}
	\begin{tabular}{r}
		\includegraphics[scale=0.3]{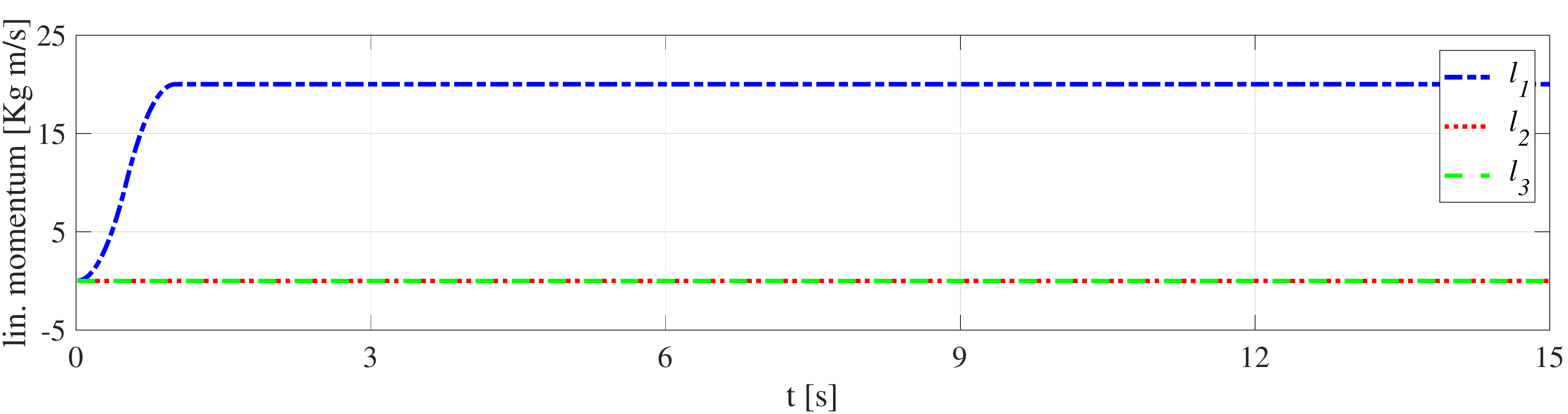}\tabularnewline
		\includegraphics[scale=0.3]{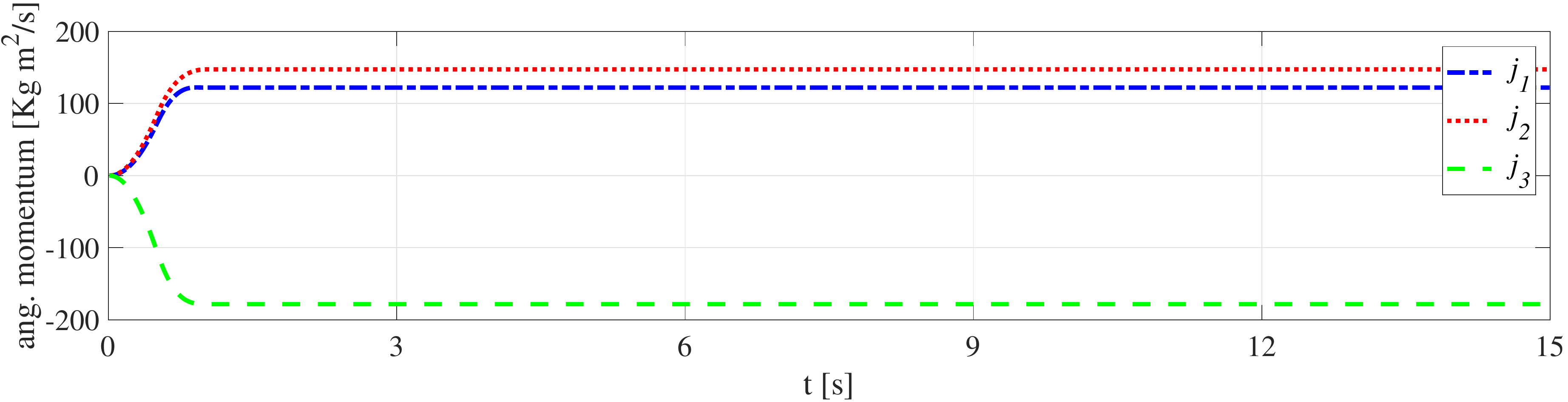}\tabularnewline
		\includegraphics[scale=0.3]{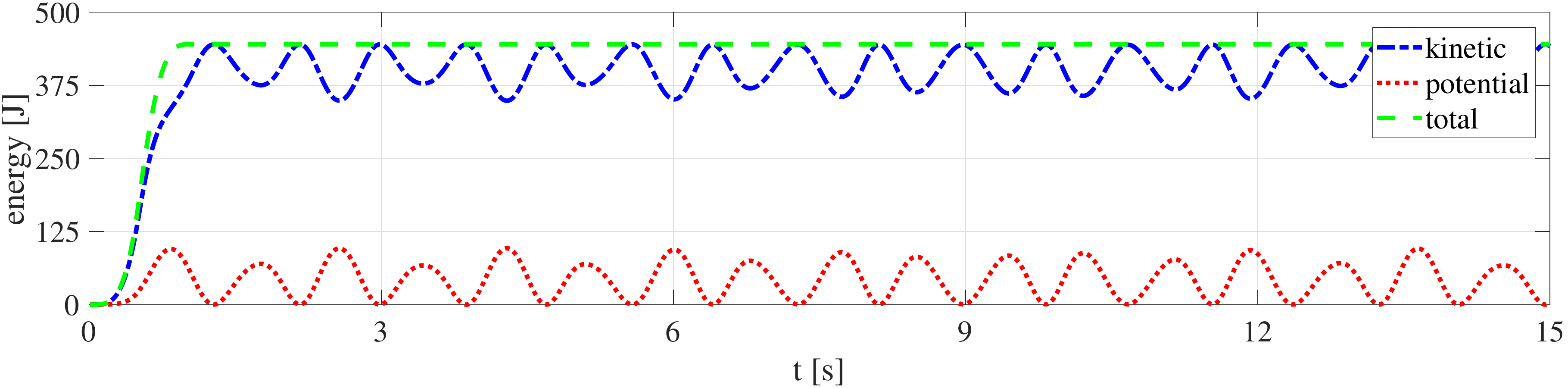}\tabularnewline
	\end{tabular}\caption{Tumbling cylinder (conservative) - momenta and energy.}
	\label{fig:tumbling cylinder conservative - momenta and energy}
\end{figure}

\begin{figure}[!ht]
	\centering{}
	\begin{tabular}{r}
		\includegraphics[scale=0.3]{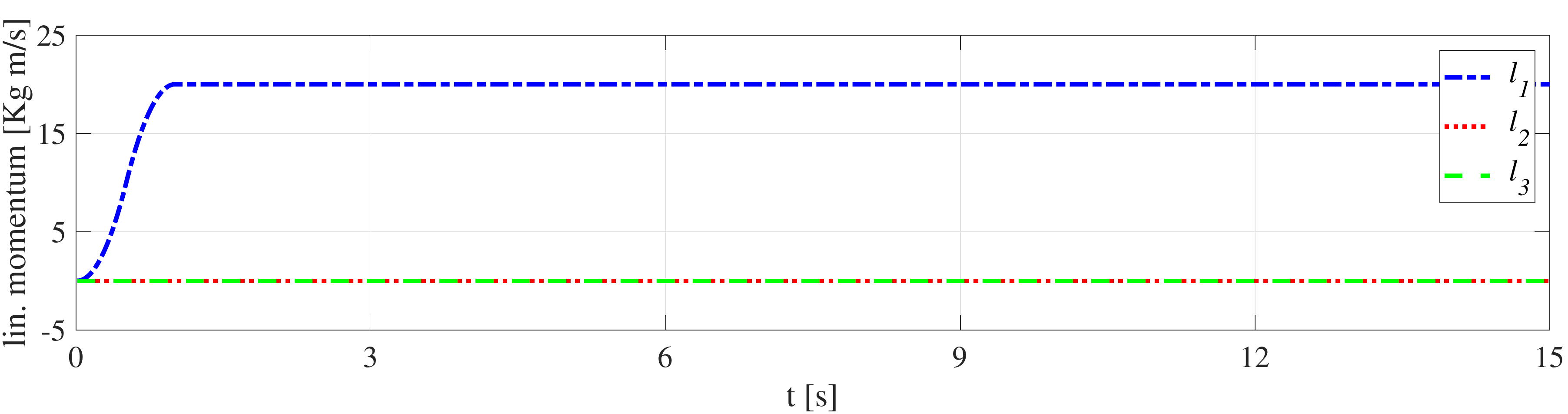}\tabularnewline
		\includegraphics[scale=0.3]{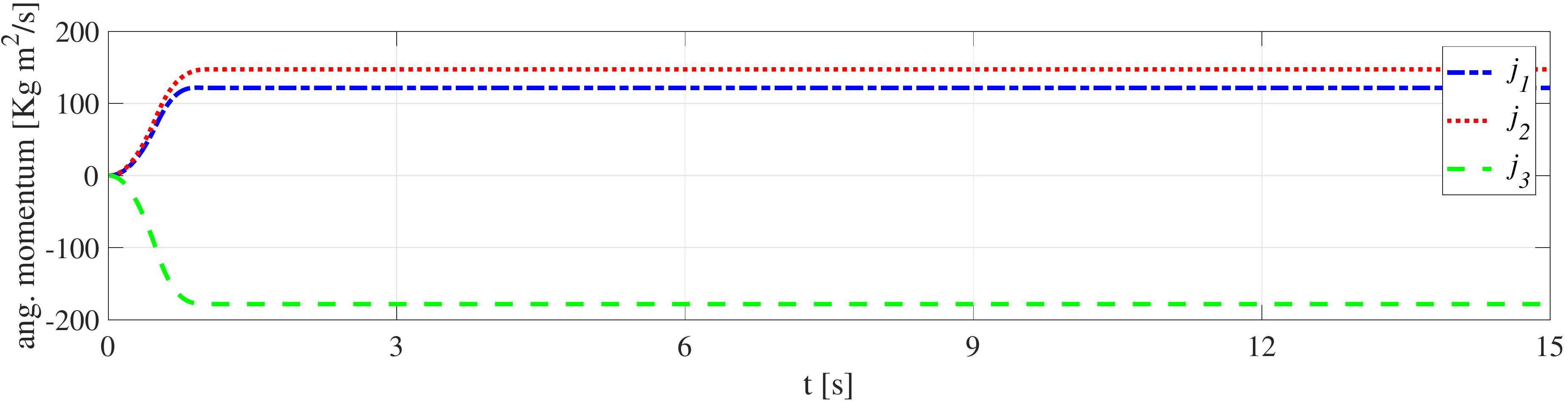}\tabularnewline
		\includegraphics[scale=0.3]{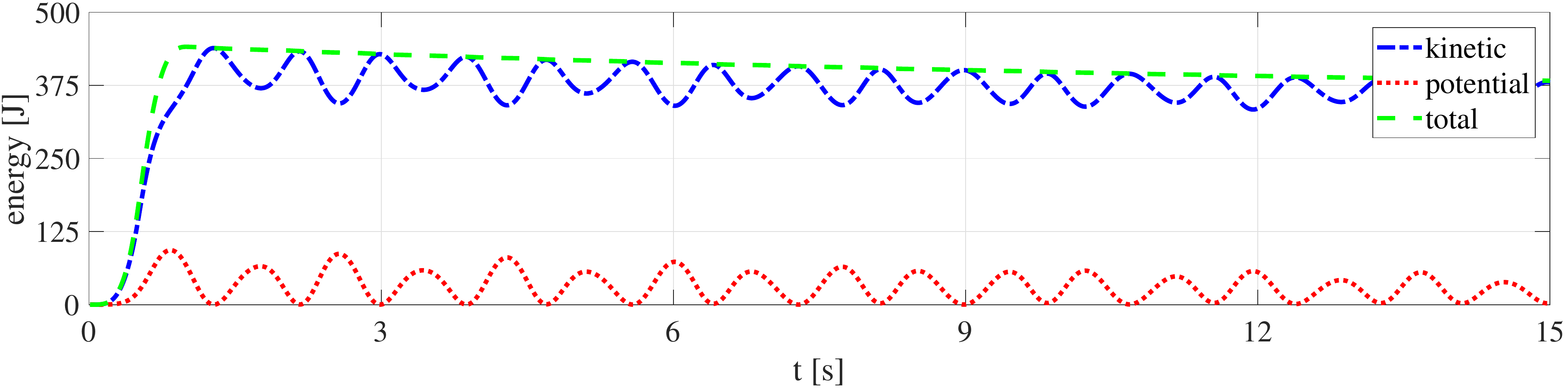}\tabularnewline
	\end{tabular}\caption{Tumbling cylinder (dissipative) - momenta and energy.}
	\label{fig:tumbling cylinder dissipative - momenta and energy}
\end{figure}

This structure is a cylindrical shell subject to body loads with a prescribed time variation
\revone{and was already investigated, for instance, in \cite{Simo1994, Betsch2009, Gebhardt2017} and in many other works}. The
geometrical and material properties are the following: mean radius $7.5\,\mathrm{m}$, height
$3.0\,\mathrm{m}$, thickness $0.02\,\text{m}$, first Lam\'e parameter $80\,\mathrm{MPa}$, second
Lam\'e parameter $80\,\mathrm{MPa}$ and mass density per volume unit
$1.0\,\mathrm{Kg/}\mathrm{m}^{3}$. The cylinder is discretized with $48$ elements, in which $16$
elements are located along the circumference and $3$ elements along the height. The total number of
nodes is $60$. Moreover, no kinematic boundary conditions are enforced. For the dissipative case we set $\chi = 0.25$.
\begin{figure}[!ht]
	\centering{}
	\includegraphics[scale=0.25]{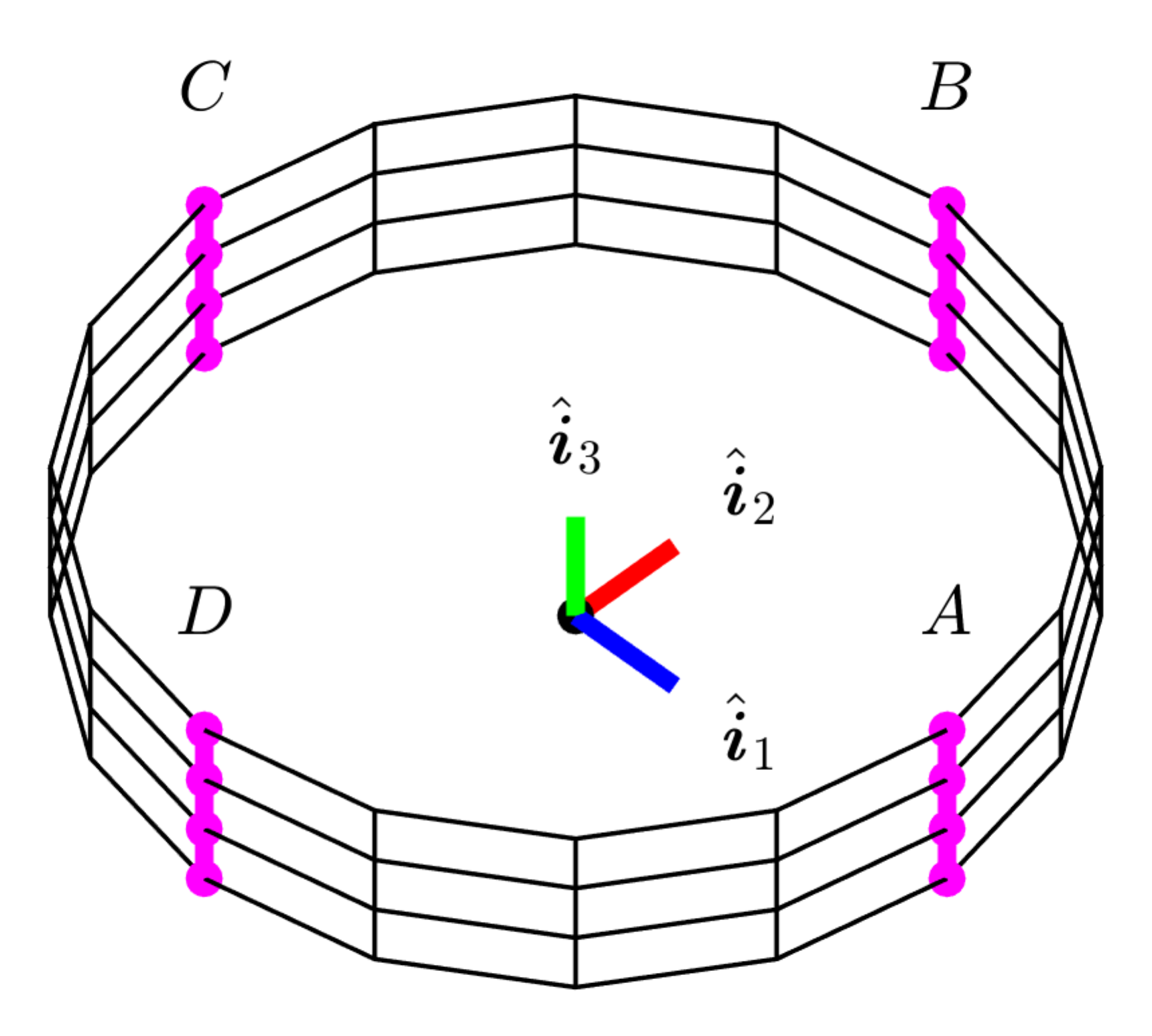}\caption{Tumbling cylinder - finite element representation.}
	\label{fig:tumbling cylinder - finite element representation}
\end{figure}
Fig. \ref{fig:tumbling cylinder - finite element representation} shows the finite element
model of the tumbling cylinder. Additionally, the line segments $A$, $B$, $C$ and $D$, to
which the spatial loads are applied, are indicated in magenta. Tab. \ref{tab:Tumbling cylinder -
spatial loads per length unit} presents the values for the loads that are applied to the
structure. The loads are then multiplied with a function that describes the variation of the applied
force over the time, which is defined in Eq. \eqref{eq:tumbling cylinder - loads time
variation}. Then
$\bm{f}_{0}^{\mathrm{ext}}=f_{1}\hat{\bm{i}}^{1}+f_{2}\hat{\bm{i}}^{2}+f_{3}\hat{\bm{i}}^{3}$ and
$\bm{f}^{\mathrm{ext}}(t)=f(t)\bm{f}_{0}^{\mathrm{ext}}$, in which the last expression is the applied
load.
\begin{table}[!ht]
	\centering{}
	\begin{tabular}{|c|r|r|r|r|}
		\cline{2-5} 
		\multicolumn{1}{c|}{} & $A$ & $B$ & $C$ & $D$\tabularnewline
		\hline 
		$f_{1}$ & $0$ & $1$ & 1 & $0$\tabularnewline
		\hline 
		$f_{2}$ & $-1$ & $1$ & $1$ & $-1$\tabularnewline
		\hline 
		$f_{3}$ & $-1$ & $1$ & $1$ & $-1$\tabularnewline
		\hline 
	\end{tabular}\caption{Tumbling cylinder - spatial loads per length unit in $\text{N}/\text{m}$.}
	\label{tab:Tumbling cylinder - spatial loads per length unit}
\end{table}
\begin{equation}
f(t)=\left\{ \begin{array}{ccc}
10t & \mathrm{for} & 0\leq t<0.5\\
5-10t & \mathrm{for} & 0.5\leq t<1\\
0 & \mathrm{for} & t\geq1
\end{array}\right.\label{eq:tumbling cylinder - loads time variation}
\end{equation}
Fig. \ref{fig:tumbling cylinder conservative - sequence of motion} shows a motion sequence for the
conservative case, where the original configuration is located at the upper-left corner of the plot,
and some deformed configurations are sequentially shown from left to right and from the top to the
bottom.  Tab. \ref{fig:tumbling cylinder - stationary values} provides the stationary values for
momenta and energy computed with the current method for both the conservative and the dissipative
cases. Fig. \ref{fig:tumbling cylinder conservative - momenta and energy} shows the time history of
momenta and energy for the conservative case. It can be observed that the linear momentum, angular
momentum and total energy vary during the time in which the external load is active, \ie~ the first
$1\,\mathrm{s}$. After the external loads vanish, these three quantities are identically
preserved through the time. These results confirm that the newly proposed integration scheme
preserves momenta and energy. Although the total energy remains constant, the potential and
kinematic energies vary in time, complementing each other in such a way that the total energy is
perfectly constant.  Fig. \ref{fig:tumbling cylinder dissipative - momenta and energy} shows the
time history of momenta and energy for the dissipative case. Clearly, the momenta is identically
preserved and energy is dissipated.
\begin{table}[!ht]
	\centering{}
	\begin{tabular}{|c|c|c|c|c|c|c|c|}
		\hline
        $t>t_{load}$ & $l_{1}$ & $l_{2}$ & $l_{3}$ & $j_{1}$ & $j_{2}$ & $j_{3}$  & $T+V$ \tabularnewline
        ($1.0$ s) & $[$Kg$\,$m/s$]$ & $[$Kg$\,$m/s$]$ & $[$Kg$\,$m/s$]$ & $[$Kg$\,$m$^2$/s$]$ & $[$Kg$\,$m$^2$/s$]$ & $[$Kg$\,$m$^2$/s$]$ & $[$J$]$ \tabularnewline
		\hline
        cons. & $20.00000$ & $0.00000$ & $0.00000$ & $122.00960$ & $147.20774$ & $-178.26475$ & $445.22767$ \tabularnewline
		\hline
		diss.  & $20.00000$ & $0.00000$ & $0.00000$ & $121.62771$ & $147.28455$ & $-178.14972$ & -- \tabularnewline
		\hline		
	\end{tabular}\caption{Tumbling cylinder - stationary values.}
	\label{fig:tumbling cylinder - stationary values}
\end{table}

\subsubsection{Free-flying single-layer plate}
\begin{figure}[!ht]
	\centering{}
	\begin{tabular}{ccc}
		$t = 0.00000$ s & $t = 0.00905$ s & $t = 0.01810$ s \tabularnewline
		\tabularnewline
		\includegraphics[align=c, scale=0.5]{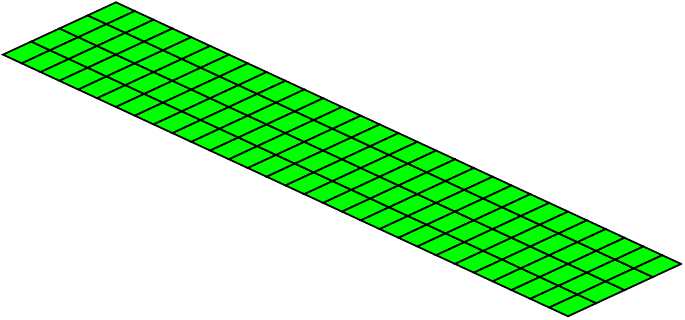} &
		\includegraphics[align=c, scale=0.5]{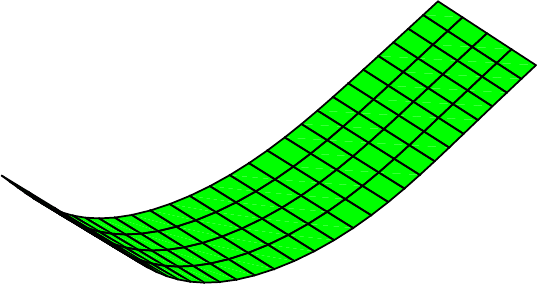} &
		\includegraphics[align=c, scale=0.5]{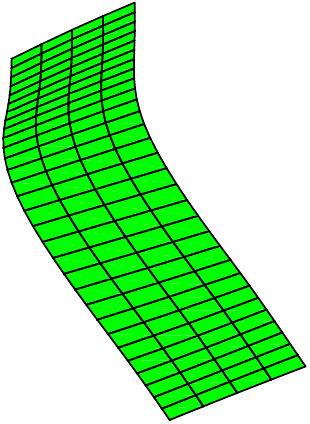} \tabularnewline
		\tabularnewline
		$t = 0.02715$ s & $t = 0.03620$ s & $t = 0.04525$ s \tabularnewline
		\tabularnewline
		\includegraphics[align=c, scale=0.5]{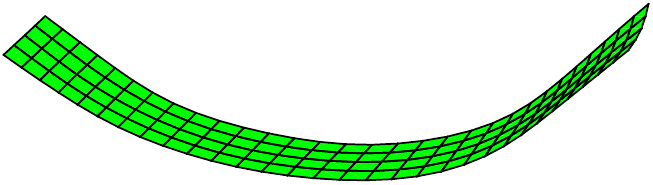} &
		\includegraphics[align=c, scale=0.5]{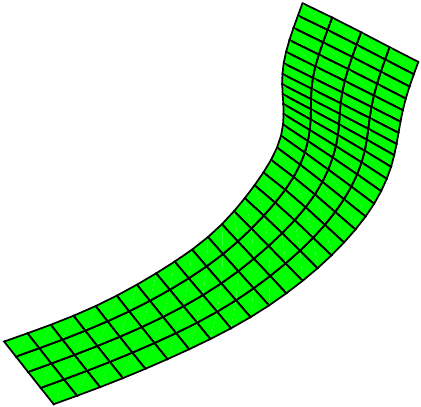} &
		\includegraphics[align=c, scale=0.5]{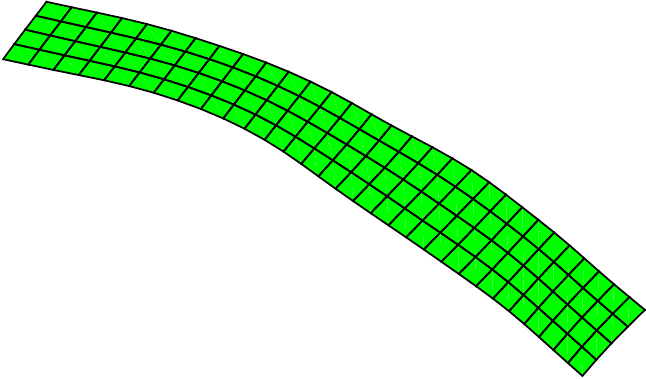} \tabularnewline
		\tabularnewline
		$t = 0.05430$ s & $t = 0.06335$ s & $t = 0.07240$ s \tabularnewline
		\tabularnewline
		\includegraphics[align=c, scale=0.5]{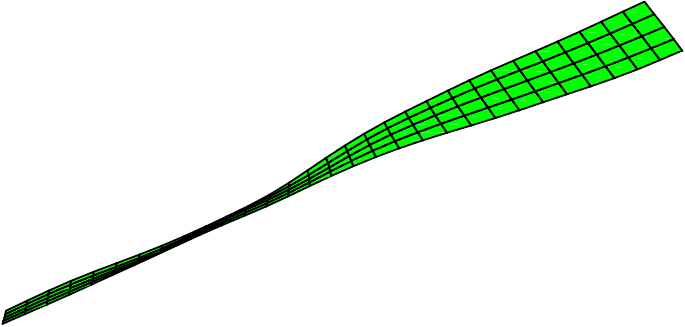} &
		\includegraphics[align=c, scale=0.5]{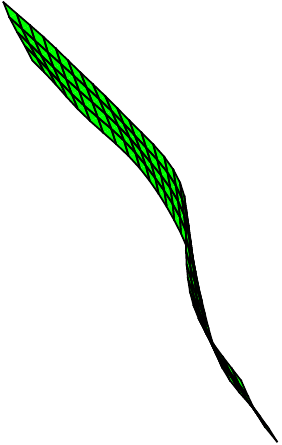} &
		\includegraphics[align=c, scale=0.5]{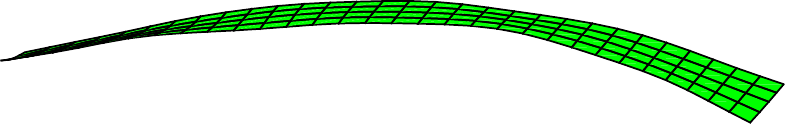} \tabularnewline
		\tabularnewline
		$t = 0.08145$ s & $t = 0.09050$ s & $t = 0.09955$ s \tabularnewline
		\tabularnewline
		\includegraphics[align=c, scale=0.5]{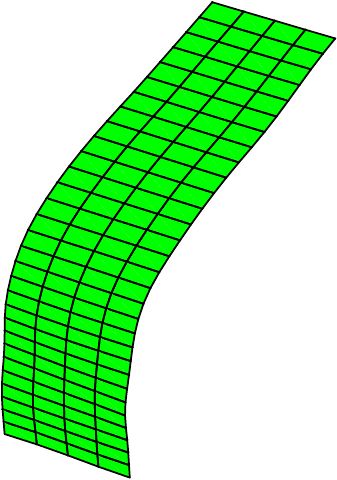} &
		\includegraphics[align=c, scale=0.5]{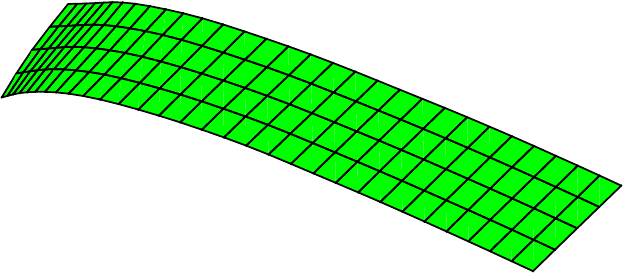} &
		\includegraphics[align=c, scale=0.5]{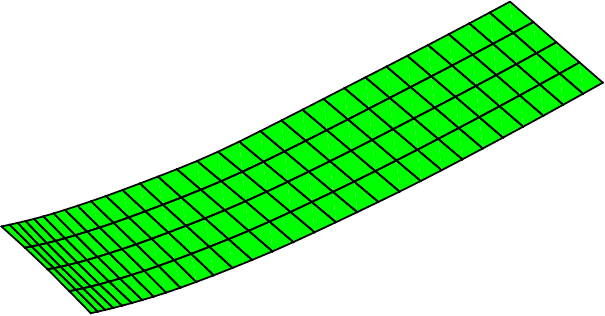}
	\end{tabular}
	\caption{Free-flying single-layer plate (conservative) - sequence of motion.}
	\label{fig:free-flying single-layer plate conservative - sequence of motion}	
\end{figure}

The structure considered in this last example
is a rectangular flat plate, which consisting of a single material layer, subject to
spatial loads with a prescribed time variation
\revone{and was considered, for example,
  in \cite{Kuhl:1996vu, VuQuoc2003, Gebhardt2017} and in many other works}. The geometrical and material properties are the
following: length $0.3\,\mathrm{m}$, width $0.06\,\mathrm{m}$, thickness $0.002\,\mathrm{m}$, first
Lam\'e parameter $0.0\,\mathrm{Pa}$, second Lam\'e parameter $103.0\,\mathrm{GPa}$ and mass density
per volume unit $7.3\times10^{3}\,\mathrm{Kg}/\mathrm{m}^{3}$. The plate is then discretized
with $120$ elements, $30$ elements being located along the largest dimension and $4$ elements along
the smallest dimension. The total amount of nodes is $155$ and for the dissipative case we set $\chi = 0.5$.
\begin{figure}[!ht]
	\centering{}
	\includegraphics[scale=0.25]{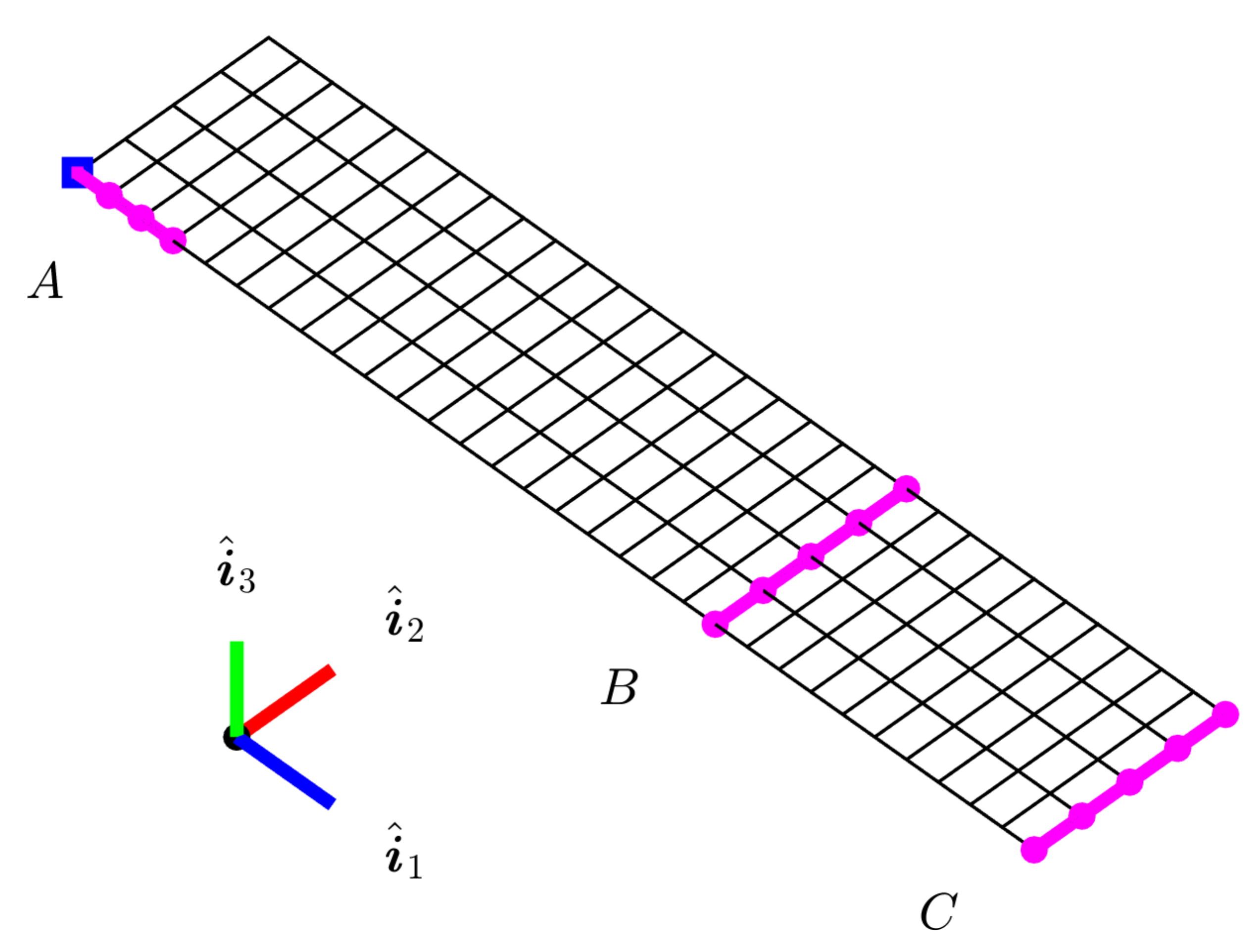}\caption{Free-flying single-layer plate - finite element representation.}
	\label{fig:free-flying single-layer plate - finite element representation}
\end{figure}
Fig. \ref{fig:free-flying single-layer plate - finite element representation}
depicts the finite element discretization of this structure. The loads are applied
over the line segments $A$, $B$ and $C$, indicated in magenta on the figure.
The reference point for the angular momentum is indicated with the symbol $\square$.
Tab. \ref{tab:free-flying single-layer plate - force density per length
  unit} gathers the values for the loads that are applied to the structure
and Eq.~\eqref{eq:free-flying plate - loads time variation} defines their scaling factor.
\begin{table}[!ht]
	\centering{}
	\begin{tabular}{|c|c|c|c|}
		\cline{2-4} 
		\multicolumn{1}{c|}{} & $A$ & $B$ & $C$\tabularnewline
		\hline 
		$f_{1}$ & $0$ & $0$ & $40000$\tabularnewline
		\hline 
		$f_{2}$ & $40000$ & $0$ & $0$\tabularnewline
		\hline 
		$f_{3}$ & $40000$ & $-40000$ & $40000$\tabularnewline
		\hline 
	\end{tabular}\caption{Free-flying singe-layer plate - force density per length unit in $\text{N}/\text{m}$.}
	\label{tab:free-flying single-layer plate - force density per length unit}
\end{table}
\begin{equation}
f(t)=\left\{ \begin{array}{ccc}
500t & \mathrm{for} & 0\leq t<0.002\\
2-500t & \mathrm{for} & 0.002\leq t<0.004\\
0 & \mathrm{for} & t\geq0.004
\end{array}\right.\label{eq:free-flying plate - loads time variation}
\end{equation}
Fig. \ref{fig:free-flying single-layer plate conservative - sequence of motion} shows a motion
sequence for the conservative case.  The linear momentum, angular momentum and energy during the
simulation are constant once reached the stationary state, and their values are provided in
Tab. \ref{tab:free-flying single-layer plate - stationary values} for both the conservative and the
dissipative cases. Momenta and energy values in time are plotted in Fig. \ref{fig:free-flying
single-layer plate conservative - momenta and energy} for the conservative case, proving that after
the removal of the force, they all remain constant. Fig. \ref{fig:free-flying single-layer plate
dissipative - momenta and energy} shows the values in time of momenta and energy for the dissipative case. Once again, momenta is perfectly preserved and energy is artificially dissipated. 
\begin{table}[!ht]
	\centering{}
	\begin{tabular}{|c|c|c|c|c|c|c|c|}
		\hline
		$t>t_{load}$ & $l_{1}$ & $l_{2}$ & $l_{3}$ & $j_{1}$ & $j_{2}$ & $j_{3}$  & $T+V$
                                                                                            \tabularnewline
          ($0.004$ s) & $[$Kg$\,$m/s$]$ & $[$Kg$\,$m/s$]$ & $[$Kg$\,$m/s$]$
                      & $[$Kg$\,$m$^2$/s$]$ & $[$Kg$\,$m$^2$/s$]$ & $[$Kg$\,$m$^2$/s$]$
                      & $[$J$]$ \tabularnewline
		\hline
		cons. & $4.80000$ & $3.20000$ & $3.20000$ & $0.02880$ & $-0.38690$ & $-0.03596$ & $246.53283$ \tabularnewline
		\hline
		diss. & $4.80000$ & $3.20000$ & $3.20000$ & $0.02896$ & $-0.38701$ & $-0.03597$ & -- \tabularnewline
		\hline
	\end{tabular}\caption{Free-flying single-layer plate - stationary values. }
	\label{tab:free-flying single-layer plate - stationary values}
\end{table}
\begin{figure}[!ht]
	\centering{}
	\begin{tabular}{r}
		\includegraphics[scale=0.3]{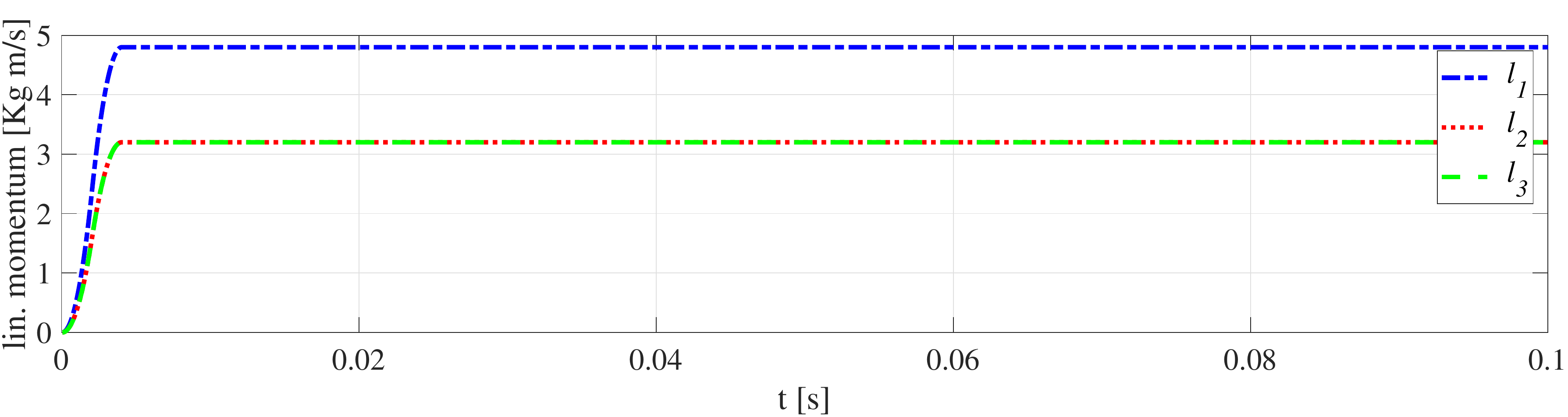}\tabularnewline
		\includegraphics[scale=0.3]{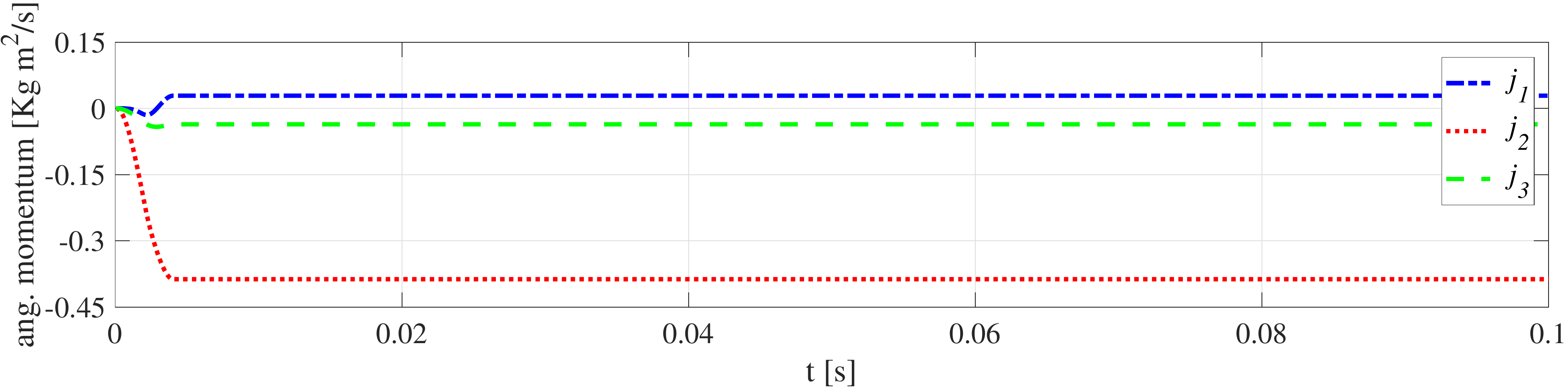}\tabularnewline
		\includegraphics[scale=0.3]{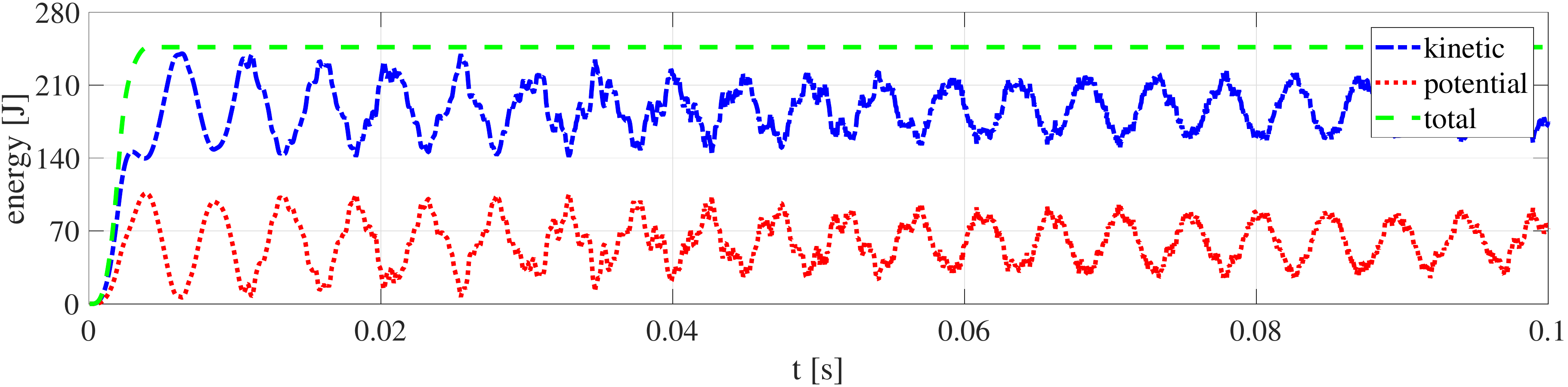}\tabularnewline
	\end{tabular}
	\caption{Free-flying single-layer plate (conservative) - momenta and energy.}
	\label{fig:free-flying single-layer plate conservative - momenta and energy}
\end{figure}
\begin{figure}[!ht]
	\centering{}
	\begin{tabular}{r}
		\includegraphics[scale=0.3]{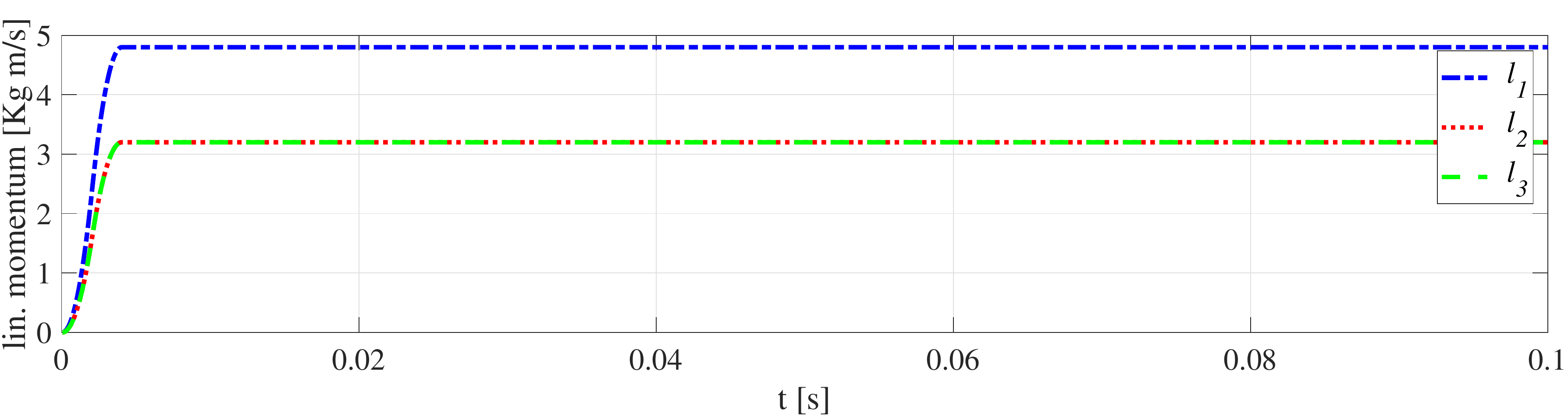}\tabularnewline
		\includegraphics[scale=0.3]{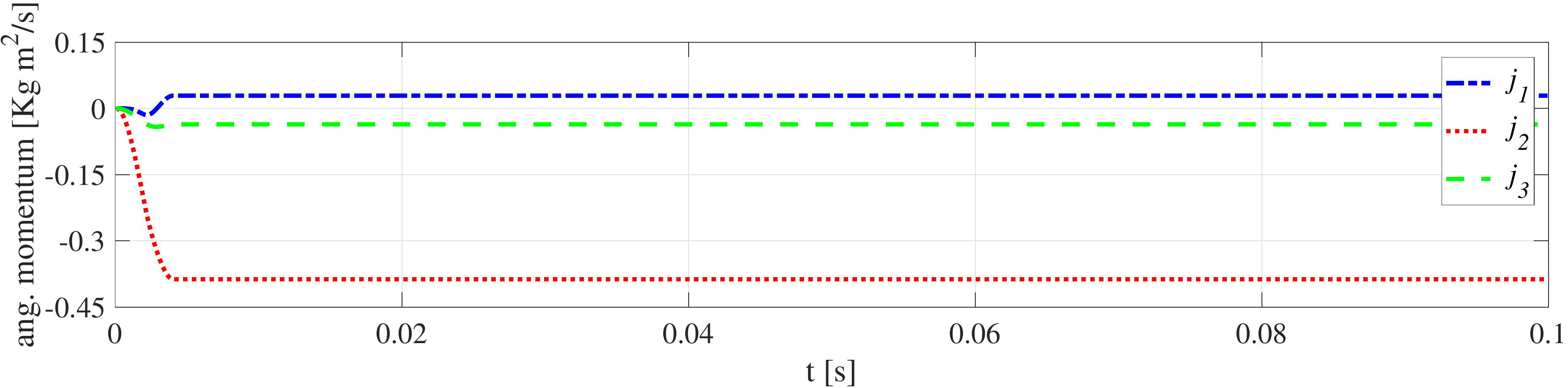}\tabularnewline
		\includegraphics[scale=0.3]{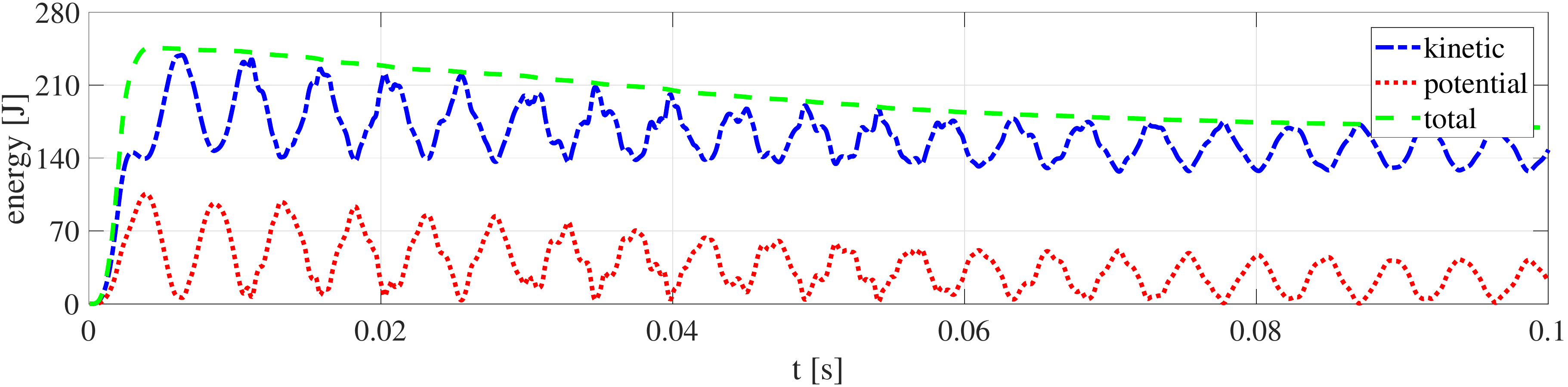}\tabularnewline
	\end{tabular}\caption{Free-flying single-layer plate (dissipative) - momenta and energy.}
	\label{fig:free-flying single-layer plate dissipative - momenta and energy}
\end{figure}

\subsubsection{\revone{On the dissipation properties}}
\revone{Dissipative schemes would be of little interest if the energy dissipation did not take place
  mostly in the high-frequency range. It is well-known that the dissipation
  of the high frequencies results in enhanced stability for the integration of stiff differential
  equations. Therefore, we could conclude that the only useful dissipative schemes are those that can
annihilate the high-frequency content of the response without radically affecting the low frequency
content of the response. In a nonlinear mechanical context, and to the best of our knowledge, there
exists no dissipation function that only eliminates the high-frequency content and leaves untouched
the low-frequency content. Dissipation always takes place along the whole frequency range. Moreover,
there is no formal proof that the dissipation can be split in that sense for nonlinear problems
and thus, we can only claim
that some dissipation functions seem to be effective to address the high-frequency problem, fact
that is mainly justified by experience. Further detailed analysis
regarding intrinsic features of dissipation functions would fall outside the scope of the current
work that addresses the derivation of a new structure preserving schema that is enriched with the
inclusion of numerical dissipation. The choice of a particular dissipation function is left to the
structural analyst based on the special demands of the problem to be solved.}

\begin{figure}[ht!]
	\centering{}
	\begin{tabular}{c c}
		\includegraphics[width=0.5\textwidth]{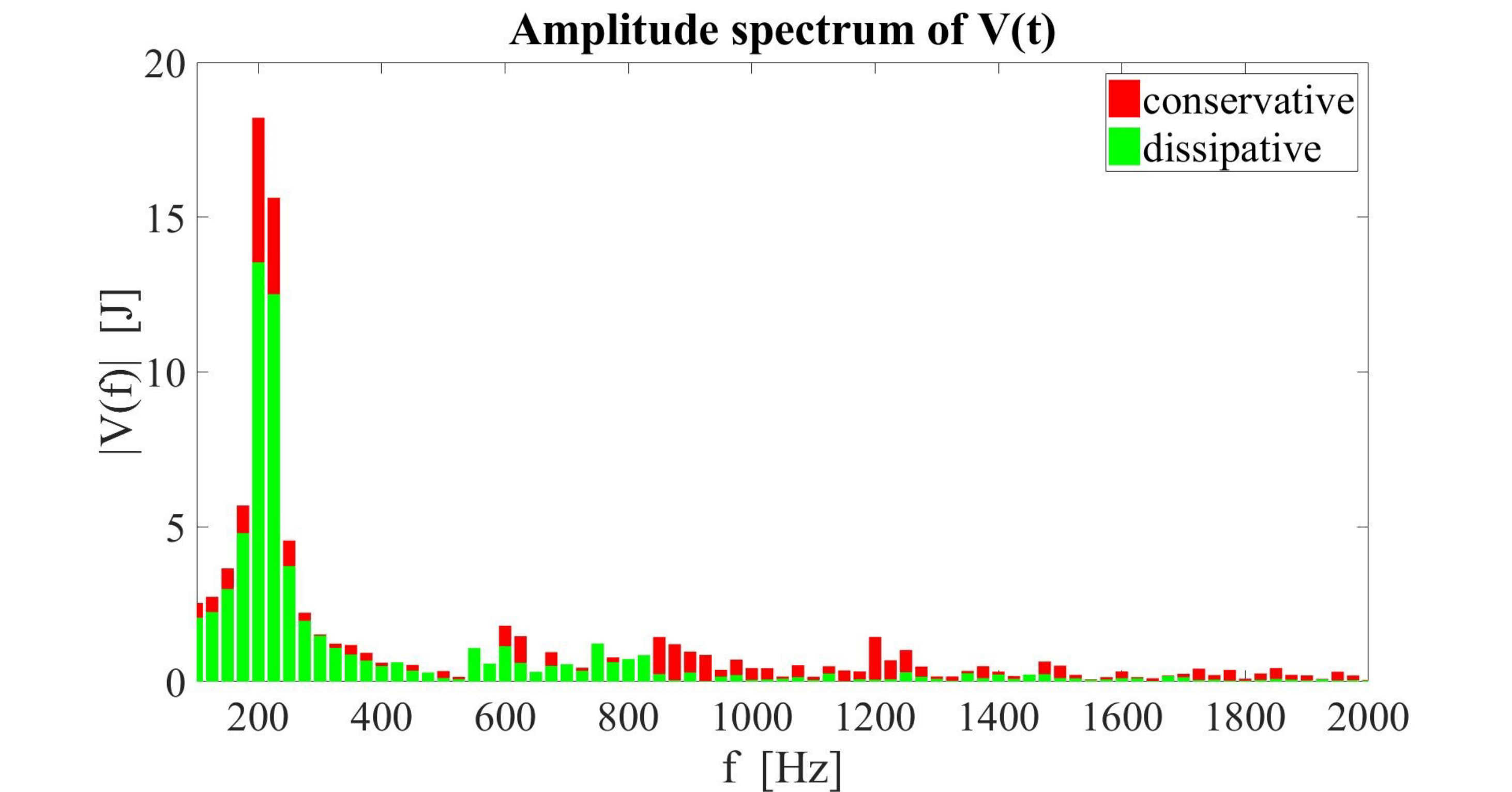} &	\includegraphics[width=0.5\textwidth]{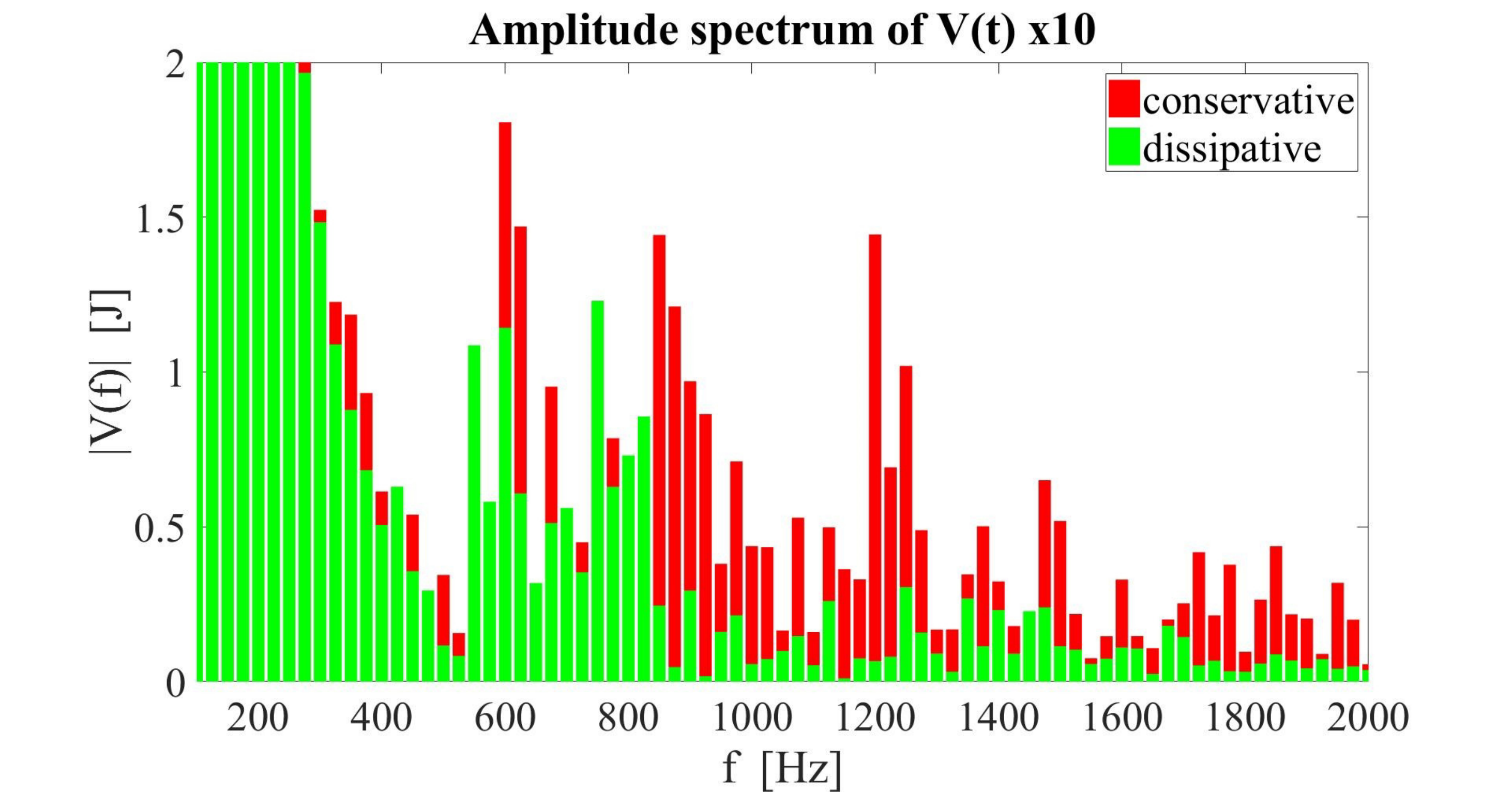}
	\end{tabular}
	\caption{\revone{Effectiveness of the dissipative scheme at the potential energy level.}}
	\label{fig:dissproppotential}
\end{figure}

\begin{figure}[ht!]
	\centering{}
	\begin{tabular}{c c}
		\includegraphics[width=0.5\textwidth]{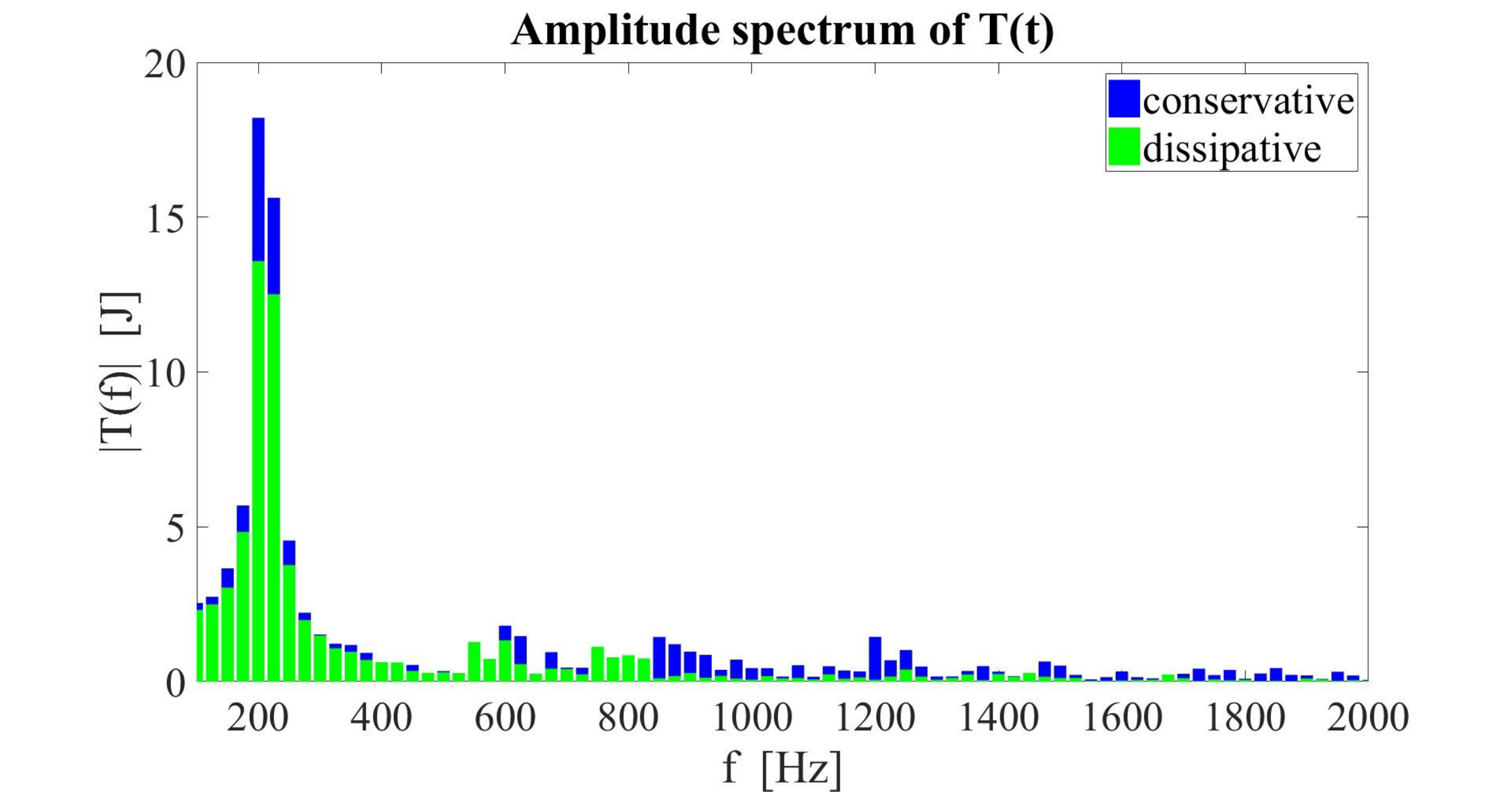} &	\includegraphics[width=0.5\textwidth]{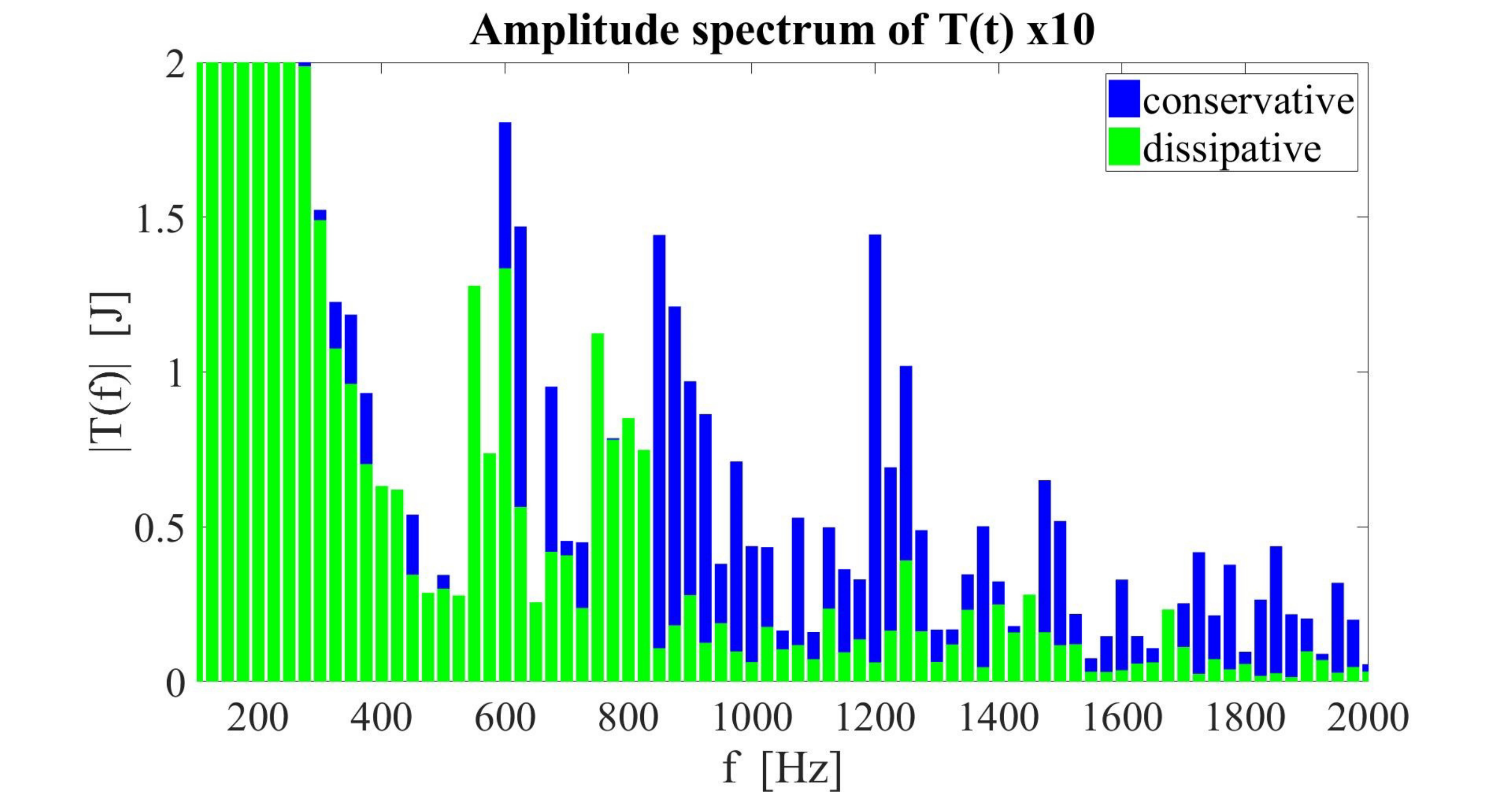}
	\end{tabular}
	\caption{\revone{Effectiveness of the dissipative scheme at the kinetic energy level.}}
	\label{fig:disspropkinetic}
\end{figure}

\revone{Keeping these limitations in mind, the free-flying single-layer plate turns to be a suitable example to show the
good dissipation properties of the new proposed scheme. Fig. \ref{fig:dissproppotential}, to
the left, presents the amplitude spectrum based on the fast Fourier transform of the potential
energy for both, the conservative and dissipative cases within the time range $0.06-0.1\textrm{ s}$
such that the direct influence of the initial transient is avoided. The subsequent analysis
corresponds to the frequency range $100-2000\textrm{ Hz}$ and to the energy amplitude range
$0-20\textrm{ J}$. Fig. \ref{fig:dissproppotential}, to the right, presents the same information,
but for the energy amplitude range $0-2\textrm{ J}$. Clearly, the dissipative algorithm works very
effectively beyond $600\textrm{ Hz}$. For the kinetic energy, Fig. \ref{fig:disspropkinetic}, to the
left and to the right, shows almost identical dissipative properties. Up to $600\textrm{ Hz}$, even
if slightly different due to some dissipation within $200-210\textrm{ Hz}$, the behavior for the
conservative and dissipative cases looks similar. Thus, we can claim that the proposed scheme seems
to have very interesting dissipation properties.}

\section{Concluding remarks}
\label{sec-closure}

We considered the conservative/dissipative time integration in
the context of nonlinear mechanical systems. A systematic
approach to derive algorithmic internal forces and generalized velocities
that ensure the preservation or the controlled dissipation of energy
was presented. As a main concrete result, we proposed a new second-order
formula for the algorithmic internal forces. Moreover, this formula
was investigated from a geometric point of view and also
interpreted for the fully conservative case in terms of a general
approach available in the literature.

In contrast with conservative/dissipative methods available
in the literature and based on the midpoint rule, the proposed formulas
are perturbations of averaged evaluations, and thus not equivalent to
existing ones.

The proposed methods are able to preserve the total energy of conservative
equations, or add artificial dissipation in a controllable fashion, while
preserving, in both cases, the linear and angular momenta of the system.
Numerical tests verify all the previous assertions.


The proposed methods could be extended to integrate differential-algebraic equations or to include
consistently dissipation functions involving derivatives with fractional orders, among others. The
reformulation in the context of \revtwo{polyconvex large strain elasticity as well as of} Lie Groups may yield interesting results. Beyond that, rigorous
mathematical proofs on the robustness are still necessary.

\section*{Acknowledgments}

Cristian Guillermo Gebhardt and Raimund Rolfes acknowledge the financial support of the
Lower Saxony Ministry of Science and Culture (research project \emph{ventus efficiens}, FKZ ZN3024)
and the German Federal Ministry for Economic Affairs and Energy (research project \emph{Deutsche
Forschungsplattform f\"{u}r Windenergie}, FKZ 0325936E) that enabled this work.


\appendix
 
\section{Precision quotient}

It is very useful to have means for checking the correctness
of integration algorithms during their development and implementation. Therefore, we introduce
here two tests that can be applied once an integration
scheme has been numerically implemented. According to Kreiss
and Ortiz \citep{Kreiss2014}, the numerical solution of an initial
value problem can be expanded as
\begin{equation}
\bm{\xi}(t,h,k)=
\bm{\xi}(t)+
\left(\frac{h}{k}\right)\bm{\psi}_{1}(t)+
\left(\frac{h}{k}\right)^{2}\bm{\psi}_{2}(t)+
\ldots+
\left(\frac{h}{k}\right)^{n}\bm{\psi}_{n}(t)+
\mathcal{O}(h^{n+1})\,,
\end{equation}
where $\bm{\xi}(t)$ is the exact solution of the given initial
value problem and $\bm{\psi}_{i}$ for $i=1,\ldots,n$ are
smooth functions of the time $t$ that do not depend on the reference
time step $h$. A positive
integer number $k$ allows to define finer solutions based on the
original resolution given by the time step $h$ that are necessary
to compute precision coefficients, a tool that may be very effective
to check the correctness of a running program. 

A first precision quotient can be defined as
\begin{equation}
  Q_{\textrm{I}}(t)=\frac{\left\| \bm{\xi}(t,h,1)-\bm{\xi}(t)\right\| }{\left\|
      \bm{\xi}(t,h,2)-\bm{\xi}(t)\right\| }\,,
\label{eq:first-quotient}
\end{equation}
where the numerator is computed as
\begin{equation}
  \left\| \bm{\xi}(t,h,1)-\bm{\xi}(t)\right\|
  =\left(\frac{h}{1}\right)^{n}\left\| \bm{\psi}_{n}(t)\right\| +\mathcal{O}(h^{n+1}) ,
\end{equation}
and the denominator is given by
\begin{equation}
\left\| \bm{\xi}(t,h,2)-\bm{\xi}(t)\right\| =\left(\frac{h}{2}\right)^{n}\left\| \bm{\psi}_{n}(t)\right\| +\mathcal{O}(h^{n+1})\,.
\end{equation}
It is possible to show that for sufficiently small time steps, the
first precision quotient can be directly approximated by $2^{n}$,
where $n$ denotes the order of accuracy of the integration method,
namely
\begin{equation}
Q_{\textrm{I}}(t)=\frac{\left(\frac{h}{1}\right)^{n}\left\| \bm{\psi}_{n}(t)\right\| +\mathcal{O}(h^{n+1})}{\left(\frac{h}{2}\right)^{n}\left\| \bm{\psi}_{n}(t)\right\| +\mathcal{O}(h^{n+1})}=2^{n}+\mathcal{O}(h^{n+1})\approx2^{n}\,.
\end{equation}
The main issue with this definition is that the exact solution of
the initial value problem is required and, in general, is not available,
especially in the context of mechanical systems involving nonlinear
constitutive relations. To circumvent this drawback, it is possible
to define a second precision quotient as
\begin{equation}
Q_{\textrm{II}}(t)=\frac{\left\| \bm{\xi}(t,h,1)-\bm{\xi}(t,h,2)\right\| }{\left\| \bm{\xi}(t,h,2)-\bm{\xi}(t,h,4)\right\| }\,,
\end{equation}
where the numerator is computed as
\begin{equation}
\begin{aligned}
\left\| \bm{\xi}(t,h,1)-\bm{\xi}(t,h,2)\right\| 
& = \left\| \left(\frac{h}{1}\right)^{n}\bm{\psi}_{n}(t)-\left(\frac{h}{2}\right)^{n}\bm{\psi}_{n}(t)+\mathcal{O}(h^{n+1})\right\| \\
& = \left(\frac{2^{n}-1}{2^{n}}\right)h^{n}\left\| \bm{\psi}_{n}(t)\right\| +\mathcal{O}(h^{n+1})
\end{aligned}
\end{equation}
and the denominator is given by
\begin{equation}
\begin{aligned}
\left\| \bm{\xi}(t,h,2)-\bm{\xi}(t,h,4)\right\|
& = \left\| \left(\frac{h}{2}\right)^{n}\bm{\psi}_{n}(t)-\left(\frac{h}{4}\right)^{n}\bm{\psi}_{n}(t)+\mathcal{O}(h^{n+1})\right\| \\
& = \left(\frac{2^{n}-1}{4^{n}}\right)h^{n}\left\| \bm{\psi}_{n}(t)\right\| +\mathcal{O}(h^{n+1})\,.
\end{aligned}
\end{equation}
Notice that this concept removes intrinsically the need for the
exact solution of the considered initial value problem. Once again,
it is possible to show that for sufficiently small time steps, the
second precision quotient can be approximated by $2^{n}$ as well
as in the case of the first precision quotient, namely
\begin{equation}
Q_{\textrm{II}}(t)=
\frac{\left(\frac{2^{n}-1}{2^{n}}\right)h^{n}\left\| \bm{\psi}_{n}(t)\right\| +\mathcal{O}(h^{n+1})}{\left(\frac{2^{n}-1}{4^{n}}\right)h^{n}\left\| \bm{\psi}_{n}(t)\right\| +\mathcal{O}(h^{n+1})}=
2^{n}+\mathcal{O}(h^{n+1})
\approx
2^{n}\,.
\label{eq:second_quotient}
\end{equation}
For the integration scheme considered in this work (an energy-conservative/dissipative
method), accuracy of second order can be guaranteed, meaning that $\log_{2}[Q_{\textrm{I}}(t)]\approx2$
and $\log_{2}[Q_{\textrm{II}}(t)]\approx2$. Let us note that for the calculation of
precision quotients, $h$ has to be chosen small enough, and the
choice may vary from case to case.
In addition, if $\|\bm{\psi}_{n}(t)\|$
is very small, both tests may fail even if the implementation
is right. For this reason it is sometime necessary to experiment with
several initial conditions and time step sizes in order to achieve
correct pictures. As a general rule, the quotients
of accuracy show better performance when the trajectories are periodic
or quasi-periodic.

\bibliographystyle{ieeetr}
\bibliography{main}

\end{document}